\documentclass[
psfig,11pt]{article}
\usepackage{amsmath,epsfig,amssymb}
\title{A Smooth Compactification of Moduli Space of Instantons
and Its Application}
\author{Bohui Chen, \\
 Dept. of Math., MIT\\
bchen@math.mit.edu}

\usepackage{amsmath,amsthm}
\begin{document}
\maketitle

\newtheorem{theorem}{Theorem}[section]
\newtheorem{assertion}[theorem]{Assertion}
\newtheorem{claim}[theorem]{Claim}
\newtheorem{conjecture}[theorem]{Conjecture}
\newtheorem{corollary}[theorem]{Corollary}
\newtheorem{definition}[theorem]{Definition}
\newtheorem{example}[theorem]{Example}
\newtheorem{figger}[theorem]{Figure}
\newtheorem{lemma}[theorem]{Lemma}
\newtheorem{prop}[theorem]{Proposition}
\newtheorem{remark}[theorem]{Remark}

\section{Introduction}
\subsection{Compactification of Moduli of Instantons}
This paper consists of two parts.
The purpose of part I is to introduce a smooth
compactification of moduli of instantons. Here, we say a space is 
smooth if it is a smooth orbifold.
It is well known that
moduli spaces of instantons play very important roles in the
study of differential four-manifolds. The most famous example is
the Donaldson theory (\cite{DK}). In general, the moduli 
spaces are not
compact. One of the fundamental problems is to give a
compactification. The most commonly used compactification is the
so-called  Uhlenbeck compactification (\cite{DK}, \cite{FU}). 
The Uhlenbeck
compactification works well for many purposes in the Donaldson
theory. But the Uhlenbeck type compactification
 is less successful  in tackling the problems in
which studies of lower strata are involved.  Examples are the
Kotschick-Morgan conjecture (\cite{KM}, \cite{G})
 on the wall-crossing
formula and the Witten conjecture (\cite{W}, \cite{CLRT},\cite{FT})
 on the
relationship between the Donaldson invariants and the
Seiberg-Witten invariants. 
 The problem with the 
Uhlenbeck compactification is that the spaces are 
only statified spaces
and very singular in lower
strata. Although the neighborhoods of lower strata can be described
(\cite{T}) 
and some applications are worked out (\cite{KM}, \cite{L}), they are intractable
in general. For example, the Kotschick-Morgan conjecture  has been
around for almost ten years and a great deal of efforts has been
made to solve it. So far, it remains 
unsolved. The author believes that the most important reason for
the lack of closure is the fact that the  Uhlenbeck compactification
is not good enough.

    For anyone who knows  algebraic geometry, it is not difficult to
     appreciate the importance of a good compactification,
 which is often the make-or-break
     point of a problem. A tremendous amount of 
literature has been devoted to
     construct a good compactification in various 
geometric problems. An important notion
     emerging from these studies is the notion of
 ``stability''. Let us take the
     example of quantum cohomologies, 
which is similar to the gauge theory in many ways.
     A pseudo-holomorphic map from a 
Riemann surface (smooth or singular) has a natural
     automorphism given by the automorphism of Riemann surface. 
Such a pseudo-holomorphic map
     is called ``stable'' if the stabilizer is {\it finite}.
 Once the stability is achieved, the
     moduli space behaves as good as a smooth orbifold 
via so-called ``virtual cycle'' or ``virtual
     integration'' techniques. The moduli space 
of stable maps was constructed by  Parker-Wolfson-Ye-Kontsevich
     and  has been proved to be extremely important for 
quantum cohomologies. Indeed, the computations
in the theory of quantum cohomologies has gone very far based on  the
moduli spaces of stable maps. Using a similar analysis,
Parker-Wolfson (\cite{PW}) pointed out
  a similar compactification should work for moduli spaces
of instantons. 
We refer these type compactification as the bubble tree compactification.
The author was pointed out by T. Mrowka that the bubble tree compactifcation
had been known by many people in  80's.

However,
there are instances where the stability fails
when the bubble tree compactification 
is used for the moduli of instanton. 
Even in the generic situation,
the corresponding moduli spaces will have singularities. As far as
the author knows, there is no general method to deal with the
singularities caused by the failure of stability. Unfortunately, a
close examination of bubble tree compactification shows that the
stability fails at ``ghost strata'' (see \S 3.2 ). 
This is very disappointed! 
To search a smooth and reasonible compactification, 
a new idea is
needed!

     The main purpose of this paper is
    to introduce a method to resolve these singularities. 
Technically, it is similar to the
    ``flip'' in algebraic geometry
 where we first 
perform a blow-up and then a blow-down.
    It is easy to do a blow-up
 to obtain a manifold with boundary. The key is that
    a blow-down can follow to obtain a closed space
 again. Topologically, it is equivalent to
the    well-known interchange of handles. 
Such a blow-down follows from our key observation that ``ghost
strata'' are closely related to the Fulton-McPherson
compactification (\cite{FuM})
 of configuration of points. This main technique is given in 
\S 4.

\subsection{Wall-Crossing Formula for 4-manifolds with $b^+=1$}
In \S 5, we give an
affirmative answer to Kotschick-Morgan conjecture. 
Our new compactifation is one of our main ingredients.
With this smooth moduli space, we are able to apply the 
localization technique which is our other ingredient
to solve the problem.
The idea can be also applied to the moduli of $U(2)$-monopoles,
which 
is more complicated than the case considered here(\cite{C}):
we also have 
a new type  compactification for this moduli space; with
this space,
 the ``virtual'' localization
technique and the method
 developed
in \S 5,
this leads to an alternative proof of
the Witten conjecture  mentioned earlier with certain condition
(\cite{CLRT}).

Let $X$ be a smooth, simply connected, closed compact 4-manifold.
Denote by $b^+$ the maximal dimension of a maximal subspace of
$H^2(X,\mathbb{R)}$ on which the cup product is positive definite.
Let $P$ be an $SO(3)$-bundle over $X$.
It is well known that the Donaldson invariants (or Donaldson 
polynomials)
based on the moduli 
space $\mathcal{M}_P(X)$, are well defined 
when $b^+(X)\geq 1$. Moreover, when $b^+\geq 2$
they are  unique in the sense that the invariants are independent
of metrics as long as they are generic. 
For $b^+=1$, the Donaldson invariants depend on the metric.
 In fact, the invariants depend on
 chambers in $H^2(X,\mathbb{R})$ determined by $P$. Let $\Delta_P$
be the set of chambers. Donaldson invariants are expressed as a map
$$
\Psi^X_{P,c}: \Delta_P\to Sym^d(H_2(X,\mathbb{Z})),
$$
where $c\in H^2(X,\mathbb{Z})$ is an integral lifting of $w_2(P)$.
$\Psi^X_{P,c}(C)$ is the Donaldson polynomial in terms of the chamber
$C$. The invariants vary with the chambers
according to a wall-crossing formula. Let $C_1, C_{-1}$ be two 
chambers.
Suppose $\Lambda$ is the collection of walls between them. Each wall
corresponds to an element of $H^2(X,\mathbb{Z})$. So we can
 identify
$\Lambda$ as a subset of $H^2(X,\mathbb{Z})$.  The following
 is conjectured 
\cite{KM} 
\vskip 0.1in
\noindent
{\bf Conjecture (Kotschick-Morgan): }{\it Suppose
$$
\Psi_{P,c}^X(C_1)-\Psi_{P,c}^X(C_{-1})
=\sum_{\alpha\in \Lambda} \epsilon(c,\alpha)\delta_P(\alpha).
$$
where $\epsilon(c,\alpha)=(c-\alpha)^2/2$, then 
$$
\delta_P(\alpha)= \sum_{i=0}^r a_i(r,d,X) q^{r-i}\alpha^{d-2r-2i},
$$
where $q$ is the intersection form of $X$ and the coefficients
$a_i$ depend only on $r,d$ and homotopy type of $X$.
}
\vskip 0.1in
\noindent
Assuming the conjecture, Gottsche \cite{G}
has derived an elegant formula
for $\delta_P(\alpha)$ in terms of modular forms. 
So the remaining question is the conjecture itself.
The core of the problem is to prove the dependence on the 
``homotopy
type'' of $X$. For example,  we show in theorem 5.12
that $a_i$'s depend
on the Euler number $\chi$ and signature $\sigma$ of $X$ which 
are well known homotopy invariants. One of the main ingredients
is  our smooth compactification.
$\delta(\alpha)$
can be computed via a localization formula of the equivariant theory
for models with $S^1$-action. But 
when we use the Uhlenbeck compactification, these models
become
very singular because of the structure of lower strata. 
In practice, computations are unavailable. It turns out
that our smooth compactification fits the equivariant structure 
perfectly. This enables us to  apply the ordinary localization 
theory.  Beginning with this observation, we reduce 
the proof of the conjecture
to computing $p_1$ of some bundle over a space, which itself
is a fibration over $X$ with fiber 
$\underline{\mathcal{M}}_K^b$. Fortunately, this
problem can also be dealt with by using the equivariant theory again. 
Therefore, we prove the conjecture. Details are given in \S 5.3.

\subsection{Acknowledgements}
This project would not be possible without 
my adviser Yongbin
    Ruan's guidance and encouragement. In fact, many ideas of this
    paper were first generated in our joint work on other projects.
    I would like to extend to him my heartfelt appreciation. 
My special thanks are due to T. Mrowka for his patience and
 invaluable suggestions, to G. Tian for helpful discussions
and comments on the paper.
    It is my pleasure to record here my deep gratitude to
An-min Li for his
advice, encouragement and perspective,
to Prof. Guosong Zhao
for his generous help, 
Without whom I would not  do mathematics. 
I  am also indebted to  
Dr. Shengda Hu, Dr. Wanchuan Zhang and Dr. Quan Zheng for useful comments.

\vskip 0.1in
\noindent
{\bf Updates: }
This is an updated  version of  the author's thesis (\cite{C}). 
Recemtly, Feehan-Leness (\cite{FL1})
announced a proof to the Witten conjecture along their cobordism approach.
In the same paper, they also mentioned that their approach 
can also work for the Kotschick-Morgan conjecture.

\section{Background}

We give a quick review of some  basic facts about instantons and their
moduli spaces for four-manifolds. Readers are
referred to \cite{DK} for most of the details.

Let $(X, g)$ be a simply connected 4-manifold with metric $g$.
There are two kinds of  manifolds considered in this paper. One is
manifolds with $b^+_2 > 0$; the other is
$S^4$ with standard $g_0$.
For $S^4$,
using the  stereographic projection $s$ from the south pole
to the equatorial plane, we
identify $s: S^4\setminus\{\mbox{south pole}\}\to R^4$.
 Alternatively,  $0$ and  $\infty$ are used for the north 
and south poles
of $S^4$ respectively.
Via $s$ we pull back the standard coordinate
functions on $R^4$ to functions $x_i, 1\leq i\leq 4,$
on $S^4$. 

Suppose 
that $P$ 
is a  principal  $G$-bundle over $X$, where
$G$ is  
either $SO(3)$ or $SU(2)$.
When $G=SU(2)$, let $E$ be the associated rank 2 vector
bundle of $P$ and $K=c_2(E)$. 
When $G= SO(3)$, choose  $c\in H^2(X,\mathbb{Z})$ to be an integral lifting 
of $w_2(P)$. 
There exists a $U(2)$-bundle $E$ such that 
$p_1(E)=p_1(P)$ and $c_1(E)=c$. 
We say that $E$ is associated to $(P,c)$.
Let $n_0$ be the biggest integer 
satisfying $d(n_0)\geq 0$ and $n_0\cong p_1(P)(\mathrm{mod} 4)$. Set
$\bar{0}= (c^2-n_0)/4$ and let
$K= c_2(E)-\bar{0}$.
Now fix a connection $a_E$ on the determinant line bundle $\det(E)$.
For $G=SU(2)$, choose $a_E=0$. For
 a connection $A$ of $E$ we write $\det(A)$
to be the induced connection on $\det(E)$.
For both cases
define
$$
\mathcal{A}_E = \{ A| \det(A)=a_E\}.
$$
 $\mathcal{G}_E$, {\it the Gauge group} of $E$,
 consists of gauge transformations 
of 
$E$ preserving $\det(E)$. Define
$
\mathcal{B}_E = \mathcal{A}_E/\mathcal{G}_E. $
$\mathcal{B}_E$
is called the {\it configuration space}.
Let  $\mathcal{A}^\ast_E\subset \mathcal{A}_E $
be the set of non-reducible connections and
$\mathcal{B}^\ast_E=\mathcal{A}^\ast_E/\mathcal{G}_E$.
All definitions can be given directly by using $P$ instead of $E$.
A connection 
$A$ is  called {\it anti-self-dual} (or ASD) if 
$$
F^+(A)=0.
$$
Class $[A]\in \mathcal{B}_E$ is called an {\it (anti-)instanton}.
$\tilde{\mathcal{M}}_E(X)$ is the collection of
ASD-connections of $E$.
The {\it moduli space} of instantons is
$$\mathcal{M}_E(X)
=\tilde{\mathcal{M}}_E(X)/\mathcal{G}.
$$
Let $\mathfrak{g}_E$ be the adjoint bundle of $E$.
There is a complex
$$
0\to
\Omega^0(\mathfrak{g}_E)
\xrightarrow{d_A}
\Omega^1(\mathfrak{g}_E)
\xrightarrow{d_A^+}
\Omega^+(\mathfrak{g}_E)
\to 0
$$
%\node{\Omega^0(\mathfrak{g}_E)}
%\arrow{e,t}{d_A}
%\node{\Omega^1(\mathfrak{g}_E)}
%\arrow{e,t}{d_A^+}
%\node{\Omega^+(\mathfrak{g}_E)}\arrow{e}
%\node{0}
%\end{diagram}
%\end{equation}
%\end{figure}
for any ASD-connection $A$.  The complex
gives the Kuranishi model for the ``virtual'' tangent space
of $\mathcal{M}_E(X)$ at $[A]$. Let $H^i_A, i=0,1,2$ 
be the cohomologies
defined by this complex. We say that
$[A]$ is {\it regular} if $H^0_A=H^2_A=0$
 and
{\it semi-regular} if $H^2_A=0$. 
We always assume that $H^2_A=0$ in this paper.
When $H^0_A\not= 0$, $A$ is called {\it reducible}, otherwise A is called
{\it irreducible}.
$\mathcal{M}^\ast_E(X)$
is the moduli space of irreducible instantons.
$[A]$ is regular iff it is irreducible. 
When $[A]$ is regular, $H^1_A$ is the tangent space 
of $[A]$ in the moduli space. 
When $[A]$ is semi-regular, first of all, $[A]$ is
reducible: its isotropy group $\Gamma_A$ is either
$S^1$  or $G$; secondly the local model of $[A]$ in the moduli space
is then given by $H^1_A/\Gamma_A$.
When $b^+_2(X)>0$ and $c_2(E)\not= 0$,
the transversality theorem (\cite{DK}, \cite{FU}) says that 
$\mathcal{M}_E(X)= \mathcal{M}_E^\ast(X)$ is smooth
for generic metrics $g$.
For $X=S^4$ with standard metric $g_0$, instantons are regular
if $c_2(E)\not= 0$. 
When $c_2(E)=0$, the trivial instanton is semi-regular and its 
isotropy group is $SO(3)$.
For  general $X$, when $G$ is $SU(2)$ and $c_2(E)=0$,
$H^2$ is nontrivial for
the trivial connections.
We exclude this case
in our paper. In general, in order to consider the trivial connection, one 
should introduce the ``thicked moduli space''(\cite{T},\cite{FM}).   
For generic metrics,
the dimension
of the moduli space is computed by the Atiyah-Singer index theory:
$$
\dim \mathcal{M}_E(X) = d(p_1(E)),
$$
where $d$ is defined by
\begin{equation}
d(n) = -2n -
\frac{3}{2}(\chi+\sigma),
\end{equation}
where $\chi$ is the Euler number of $X$ and $\sigma$
is the signature of $X$.

One would see that when $X= S^4$, the dimensions are $5(\mathrm{mod }8)$.
In fact  when $X= S^4$, we are  more interested in a  modified
moduli spaces.
These are the original moduli spaces modulo a 5-dimensional
group $H$.
$H$ is a subgroup of the conformal group of $S^4$. It is generated by
translations and  dilations on $R^4$. Any conformal
transformation induces an action on  moduli spaces. So we can
define $\mathcal{M}^b_E = \mathcal{M}_E(S^4)/H$.
Alternatively
we can describe this moduli space more  explicitly:
let $\hbar$ be a constant less than $4\pi^2$. Define
$\tilde{\mathcal{M}}^b_K$ to be  the collection of ASD-connections $A$
with the properties
\begin{equation}
\int_{R^4} x |F(A)|^2 =0
\mbox{ and }
\int_{B(1)} |F(A)|^2 = \hbar.
\end{equation}
The first equation requires that the mass center of $A$ is 0
and the second one requires the energy in unit disk $B(1)$
to be equal to a small constant.
We call (2) the {\it balanced condition}.
An ASD-connection $A$ is {\it balanced}
if it satisfies (2).
Then $\mathcal{M}^b_E = \tilde{\mathcal{M}}^b_E(S^4)/\mathcal
{G}_E$, where $\tilde{\mathcal{M}}^b_E(S^4)$ is the set of
balanced ASD-connections.

Fix a point  $x_0\in X$. When $X=S^4$, we always
choose $x_0=\infty$.
 Suppose 
$\mathcal{G}_E^0$ is the subgroup of $\mathcal{G}_E$
that consists of gauge transformations with identity at point
$x_0$. Then
$$
u: \mathcal{A}^\ast_E/\mathcal{G}^0_E 
\to \mathcal{B}^\ast_E 
$$
is an $SO(3)$-bundle. We denote the bundle
by $\mathbb{P}^0_K$ when $G=SU(2)$, and by $\mathbb{P}^0_{K,c}$
when $G=SO(3)$. If no confusion is caused, we just write
it as $\mathbb{P}^0$. Alternatively, this bundle can be given by
$$
u: \mathcal{A}_E^\ast \times_{\mathcal{G}_E} P_{x_0}
\to \mathcal{B}^\ast_E.
$$
Let $\mathcal{M}^0_E(X)=u^{-1}(\mathcal{M}_E(X))$.
There is another important $SO(3)$-bundle over $\mathcal{B}^\ast
\times X$ which plays a crucial role in the  Donaldson theory. 
This bundle is defined by
$$
v: \mathcal{A}_E^\ast \times_{\mathcal{G}_E} P
\to \mathcal{B}^\ast_E\times X.
$$
We denote this bundle by $\mathbb{P}_K$ or $\mathbb{P}_{K,c}$
depending on $G$. In general, we  write it as $\mathbb{P}$.
Similarly, $\mathcal{M}^0_{E,X}=v^{-1}(\mathcal{M}_E(X)
\times X)$.
Using $p_1(\mathbb{P})$, one  defines
$$
\mu: H_\ast (X;\mathbb{Z}) \to H^{4-\ast}(\mathcal{B}^\ast_E;
\mathbb{Q})
$$
by taking  the slant product
$\mu([-])= -\frac{1}{4} p_1(\mathbb{P})/[-]$.
 For the 0-class $[x_0]$, $\mu([x_0])=
- \frac{1}{4}p_1(\mathbb{P}^0)$.

    As we mentioned in the introduction, the bubble tree
    compactification is our initial step towards our new compactification.
It
is based on the following compactness theorem.
\begin{theorem} 
     [Parker-Wolfson] Suppose that
$[A_n]_{n=1}^\infty$ is a sequence of instantons in 
$\mathcal{M}_E(X)$. Then there exists a subsequence of $[A_n]$
that converges to a bubble tree instanton $[\mathbf{A}]$
in $\overline{\mathcal{M}}_E(X)$.
\end{theorem}
Definitions of $\overline{\mathcal{M}}_E(X)$ and 
{\it bubble tree instantons} are given in \S 3.1.
Suppose $P,E, K$ are as above. 
Let $E_{-l}, K\geq l>0$, be the bundles such that
$c_1(E_{-l})=c_1(E)$ and $c_2(E_{-l})=c_2(E)-l$.
Roughly speaking, a bubble tree instanton
$[\mathbf{A}]$
is a sequence of finite instantons $([A_0],[A_1],\ldots,
[A_k])$
 on components:
$X$ and 4-spheres $S^4_i,1\leq i\leq k$.
$X$ is called the {\it principal component}
and $S^4_i$'s are called {\it bubbles}. 
We call $[A_0]$  the {\it background instanton} 
of $[\mathbf{A}]$. $[A_0]$ should be in $\mathcal{M}_{E_{-l}}(X)$
for some $l$. 
Then  we say that $[\mathbf{A}]$ is in the {\it $l$-level
stratum}.
Suppose $X\not= S^4$.
If $[A_0]$ is trivial, we say that 
$[\mathbf{A}]$ is trivial. If $[A_0]$
is reducible, we say $[\mathbf{A}]$ is reducible.
Let $\mathbf{F}_E$ 
be the collection of trivial bubble tree instantons. 
In this paper, we assume that $\mathbf{F}_E=\emptyset$.
For example, when $G=SO(3)$ and $w_2(P)\not= 0$, 
$\mathbf{F}_E=\emptyset$ (\cite{K}). 
Essentially, this is the case that we are interested in.
To avoid the confusion caused 
by the complexity of
notations, we restrict  to $G=SU(2)$ and make the assumption that 
$\mathbf{F}_E=\emptyset$. Also note that  $E$
 is determined by $K$ up to bundle isomorphisms. We replace 
the subscript ``$E$'' by $K$ in general. For example, 
$\mathcal{M}_K(X)$ is the same as $\mathcal{M}_E(X)$ from now on.
If we consider $G=SO(3)$, then the subscript should be $K,c$.

In the next section , we will define $\overline{\mathcal{M}}_K(X)$
and its topology, then
study its topological and smooth structure using the gluing
theory.

We now review the problem of the wall-crossing formula. 
In the rest of the chapter,
assume $b^+_2(X)=1$ and
$P\to X$ is an $SO(3)$ principal bundle with $w_2(P)\not= 0$.
Let $c$ be an integer lift of $w_2(P)$ and
$E$ be a $U(2)$-bundle over $X$
with $c_1(E)=c, p_1(E)=p_1(P)$. 
Set $K= c_2(E)-\bar{0}$ as before. 
Suppose $g_X$ is a generic metric. Let 
$\overline{\mathcal{M}}^u_K(X,g_X)$,$\overline{\mathcal{M}}_K(X,g_X)$
and
$\underline{\mathcal{M}}_K(X,g_X)$ 
be the compactified spaces
of $\mathcal{M}_K(X,g_X)$ in the sense of the Uhlenbeck 
compactification, the bubble tree compactification and the
smooth compactification which is given later in \S 4. The map 
$\mu: H_2(X, \mathbb{Z})\to H^2(\mathcal{M}_K(X,g_X), \mathbb{Q})$
defines 2-classes $\mu(\Sigma), \Sigma\in H_2(X)$. 
$\mu(\Sigma)$ can be extended to $\bar{\mu}(\Sigma)$ over
$\overline{\mathcal{M}}^u_K(X,g_X)$
(\cite{DK}, \cite{FM}). There are  two different approaches:
in \cite{DK}
a line bundle $\mathcal{L}_\Sigma$ over
$\overline{\mathcal{M}}^u_K
(X,g_X)$ is constructed
and $\bar{\mu}(\Sigma)$
is the $c_1$ of the bundle; in \cite{FM} 
the proof  is more algebraic-topological. Here
we follow the approach of  \cite{DK}. 
Note that there are natural projections 
from both $\overline{\mathcal{M}}_K(X,g_X)$
and
$\underline{\mathcal{M}}_K(X,g_X)$ to 
$\overline{\mathcal{M}}^u_K(X,g_X)$.  So $\bar{\mu}(\Sigma)$
is also well defined over these two spaces via pull-back maps. 
The Donaldson invariants are defined to be the pairing of
compactified moduli space with 
cohomology classes given by $\mu$. 
One can show that the Donaldson invariants are independent of compactified
spaces we
choose.
However, when $b^+_2(X)=1$, the pairing is no longer
metric independent:  let $g_{-1}, g_1$ be two generic metrics
and 
$\lambda=g_t$ be a path of metrics connecting them. Define
$$
\mathcal{M}_K(X, \lambda):=
\{ ([A],t)| [A]\in \mathcal{M}_K(X,g_t)\}.
$$
This forms a ``cobordism'' between $\mathcal{M}_K(X,g_i),i=-1,1$.
Unfortunately, it is not an actual cobordism since
there may be reducible connections in the space. 
We can still apply compactification theories to this space. 
Then $\overline{\mathcal{M}}^u_K(X,\lambda)$ 
forms a ``cobordism'' between $\overline{\mathcal{M}}
^u_K(X,g_i)$, and so is true for other two different 
compactifications. $\bar{\mu}(\Sigma)$ are not well defined
wherever
reducible connections appear. Even worse, reducible connections
can be at lower strata. Suppose $A$ is
a reducible connection on the top stratum.
Then it gives a reduction
of $P$ to an $S^1$-bundle $Q_{\alpha}$, $P= Q_\alpha\times_{S^1}SO(3)$.
Here $c_1(Q_\alpha)=\alpha, \alpha\in H^2(X,\mathbb{Z})$.
We say that $A$ is {\it reducible with respect to $\alpha$}. 
Let $L_\alpha$ be the line bundle associated to $Q_\alpha$.
Correspondingly, $E$ splits as 
$$
E= L_\alpha^{1/2}\otimes L_c^{1/2}
\oplus L_\alpha^{-1/2}\otimes L_c^{1/2}.
$$
If $A\in \mathcal{M}_K(X,g)$, $F^+_A=0$. Suppose that
$\omega(g)$ 
is the self-dual harmonic 2-form with respect to $g$. Then 
$\alpha\cdot \omega(g)=0$. In general, the reducible connection
$[A]$ with respect to $\alpha$ is in
$\overline{\mathcal{M}}^u_K(X, g)$ if and only if 
$$
\omega(g)\in W^\alpha:=\{
x\in H^2(X,\mathbb{R})| x^2>0, x\cdot \alpha=0\}.
$$
where $H^2(X,\mathbb{R})$
is identified with the space of harmonic forms. We call
$W^\alpha$ {\it the wall associated to $\alpha$}.
$W^\alpha$ is called a wall of $P$-type 
if $w_2(P)$ is the reduction of $\alpha$ modulo 2 and 
$0>\alpha^2\geq p_1(P)$.
The set of all walls $W^\alpha$ of $P$-type is called
$P$-walls. The set of chambers $\Delta_P$ 
of $X$ is the set of components of the complement, in the positive cone
of $H^2(X,\mathbb{R})$, of the $P$-walls.
The Donaldson invariants
associated to a bundle $P$  depend on the chamber where
the period point of the metric is located. The invariants vary when 
period points of metrics pass through $P$-walls. Let
$d= -p_1(P)-3$. 
$\Phi^X_{P,c}: \Delta_P\to Sym^d(H_2(X,\mathbb{Z}))$
is the map such that 
$\Phi^X_{P,c}(C)$ 
defines the Donaldson invariants for the chamber $C$.

We now state the Kotschick-Morgan conjecture in terms of
the Uhlenbeck compactification. 
Suppose $[A]\in \overline{\mathcal{M}}^u_K(X,\lambda)$
 is a reducible connection with respect to 
$\alpha$. Let $r=(\alpha^2-p_1(P))/4$. Then the family
of reducible connections is
$[A]\times Sym^r(X).$ Let $D(\alpha)$ denote the 
link 
of the family of reducible connections. 
This defines an element $\delta_P(\alpha)$
in $Sym^d(H_2(X,\mathbb{Z}))$ as 
follows: let
$z\in Sym^d(H_2(X,\mathbb{Z}))$,
$$
\delta_P(\alpha)(z)= \langle \bar{\mu}(z), D(\alpha)
\rangle.
$$
If $C_{-1}$ and $C_1$ are chambers, then  
$$
\Phi^X_{P,c}(C_1)-\Phi^X_{P,c}(C_{-1})
=\sum_\alpha \epsilon(c,\alpha)\delta_P(\alpha).
$$
where $\epsilon(c,\alpha)=(-1)^{(c-\alpha)^2/2}$.
\begin{conjecture}
[\cite{KM}] With notations  as above,
$$
\delta_{K}(\alpha)= \sum_{i=0}^r a_i(r,d,X)q^{r-i}
\alpha^{d-2r-2i}.
$$
where $q$ is the intersection form of $X$ and the coefficients
$a_i$ depend only on $r,d$ and homotopy type of $X$.
\end{conjecture}

\section{Gluing Theory, The Topological and Smooth Structures 
of Moduli Space $\overline{\mathcal{M}}_K(X)$}
    Before we resolve the singularities of
    $\overline{\mathcal{M}}_K(X)$, we first have to understand its
    smooth structure and the structure of singularities. 
The chapter is organized as follows: in \S 3.1, we define the
bubble tree compactified space $\overline{\mathcal{M}}_K(X)$;
The topological
structures  and the smooth structures
are discussed in \S 3.2 and \S 3.3 respectively
by  using the gluing theory.

\subsection{Bubble Trees and Strata }
Let $\mathcal{M}_K(X), \mathcal{M}^b_K$ be as before.
We denote their compactified spaces by
$\overline{\mathcal{M}}_K(X), \overline{\mathcal{M}}_K^b$.
They are stratified spaces. We begin
with the description of  their strata.

Each stratum is
associated with a  connected
tree with  additional requirements.
Let $T= (V, D)$ be a tree,
where $V$ is the set of vertices  and $D$ is the set of edges.
Choose a vertex $v_0\in V$ to be the {\it root} of $T$.
For each $v\in V$
there is a unique path connecting $v$ and $v_0$.
The length of the path (i.e,
the cardinality of edges in the path) is called $depth(v)$.
We say that vertex
$v_i$ is an {\it ancestor} of vertex $v_j$ if the unique path
 connecting $v_i$
and $v_j$ does not pass through $v_0$ and
$depth(v_i) < depth(v_j)$, and
we also call $v_j$  a {\it descendant}
 of $v_i$. Moreover,
if $depth(v_j)=depth(v_i)+1$,
$v_i$ is called a {\it parent} of $v_j$
and $v_j$ is a {\it child} of $v_i$.
In  other words, $v_i$
is an ancestor of $v_j$ and there is an edge $e=(v_i, v_j)$ in $D$.
From now on, when we write an edge $e=(v_i, v_j)$ we mean that $v_i$
is the parent of $v_j$.
For each vertex $v\in V$, let $V(v)$ be the set of vertices  consisting
of $v$ and all descendants of $v$. Define $t(v)= (V(v), D(v))$ to  be
 the subtree
of $T$ induced by $V(v)$, i.e, $D(v)$
consists of all edges in $D$ connecting  vertices
in $V(v)$. We now assign each
 vertex $v$ a nonnegative integer $w(v)$.
$w(v)$ is called the  {\it weight} or {\it charge} of $v$.
We call such a $T$ to be a {\it weighted} tree.
Given a vertex $v$, define $W(v)$ to be the sum of charges
of all vertices
in $t(v)$. We call $W(v)$ to be the {\it total charge of} $v$.
Let $child(v)$ be the set of all children of the vertex $v$.
\begin{definition}
A weighted tree $T$ with a prechosen root $v_0$
is a \emph {bubble tree} if for any vertex $v$
we have either
\begin{itemize}
\item
$ w(v)\not= 0, $ or
\item
$|child(v)| \geq 2$, and $W(v_i) >0$ for all $v_i\in child(v)$.
\end{itemize}
We write the bubble tree  as $(T,v_0)$ or $T$.
\end{definition}
Note that for each vertex $v$ in a bubble tree $T$ the associated
subtree $t(v)$ is  still a bubble tree. And $t(v_i)$ is called
a {\it bubble tree component} of $v$ if $v_i\in child(v)$.
Given an edge $e=(v_i,v_j)$, we say that the
tree $T'=(V',D')$ is the {\it contraction of $T$ at edge $e$}
if $V'=V\setminus\{v_j\}$ and $D'$ is the union of edges
  in $D$ that do not connect $v_j$ and 
new edges $e'=(v_i,v_l), v_l\in child(v_j)$. $w(v_i)$
is updated to $w(v_i)+w(v_j)$. If $T$ is a bubble tree, $T'$
is also a bubble tree with the same root. 
For a bubble tree  $(T, v_0)$,
define the  {\it total charge $W(T)$ of $T$} to be $W(v_0)$.
Given an integer $K>0$, let
$$
\mathcal{T}_K = \{(T,v_0)| (T,v_0) \mbox{ is a bubble tree and }
 W(T)=K \}.
$$
We introduce
a {\it partial order} on $\mathcal{T}_K$: $T_1 < T_2$ if
$T_2$ is obtained from $T_1$ by a sequence of contractions.
It is easy to prove
\begin{lemma}
$|\mathcal{T}_K| <\infty$.
\end{lemma}

Now we focus on $X=S^4$ for a moment.
Suppose $(T,v_0)$ is a bubble tree and $T=(V,D)$.
We assign  4-spheres $S^4_i$ to each vertex $v_i\in V$.
For each edge $e=(v_i, v_j)$ or $v_j\in child(v_i)$, we assign a point
$d_T(e)=p_{ij}\in S_i^4$. Moreover $d_T(e)\not= d_T(e')$ if
$e\not= e'$.
Abstractly, we define
a bubble tree space $Y$ of $T$
to be a pair
$
Y=(\coprod_{i} S^4_i, d_T(D)).
$
Geometrically, $Y$ is realized as a quotient
space $\coprod S^4_i/\sim$: for each $e=(v_i,v_j)$,
$d_T(e)\sim \infty_j$. We denote $Y$ by
$\coprod_i S^4_i/d_T$.
For any $v_i$ and edge $(v_i,v_j)$, $d_T(v_i,v_j)$
is called a {\it bubble point} on $S^4_i$.
For $\mathfrak{p}=(p_1,\ldots, p_n)\in 
(S^4\setminus\{\infty
\})^n$ we assign charges (or energies)
$\mathfrak{w}=(w_1,\ldots, w_n)\in (\mathbb{Z}^+)^n$
 to it, i.e, 
$w_i$ is the charge of $p_i$.
$S_n$, the
$n$-permutation group, acts on $\mathfrak{p}$
and $\mathfrak{w}$ in the standard
way.
Define $S_\mathfrak{w}<S_n$ to be the kernel of action on 
$\mathfrak{w}$ and
$$
N_\mathfrak{w}: = (S^4\setminus\{\infty\})^n/S_\mathfrak{w}.
$$
The class of $\mathfrak{p}$ in $N_\mathfrak{w}$
 is denoted by
$[\mathfrak{p}]=[p_1,\ldots, p_n]$.
Note that the group action $H$ on $S^4\setminus\{\infty\}$
is also well defined on $N_\mathfrak{w}$.
\begin{definition}
A \emph{ generalized instanton}
on $S^4$ is an instanton $[A]\in \mathcal{M}_{k_0}(S^4)$
with
$[\mathfrak{p}]=[p_1,\ldots, p_n]\in N_\mathfrak{w}$,
where $\mathfrak{w}= (k_1,\ldots, k_n)$.
Charges $k_i$ are called the $\delta$-\emph{ mass}
at $p_i$. The moduli space of generalized instanton
denoted by  $\mathcal{M}_{k_0,\mathfrak{w}}
(S^4)$ is $\mathcal{M}_{k_0}(S^4)\times N_\mathfrak{w}$.
Define
$$
\mathcal{M}^b_{k_0,\mathfrak{w}}=
\mathcal{M}_{k_0,\mathfrak{w}}(S^4)/H.
$$
The element $[[A],[\mathfrak{p}]]$ in
$\mathcal{M}^h_{k_0,\mathfrak{w}}$ is called a {\it balanced
 instanton}.
The representative
of $[[A],[\mathfrak{p}]]$ can be chosen to satisfy
the \emph{ balanced condition}:
\begin{enumerate}
\item
the (generalized) mass center  $m$ is the north pole, namely
$$
m: = x(m([A]))k_0 + x(p_1)k_1 +\cdots+ x(p_n)k_n =0
$$
\item
one of the following  cases holds:
\begin{enumerate}
\item
the charge of $[A]$ at south hemi-sphere equals the
constant $\hbar$ and $p_i$'s are located in the open north semi-sphere,
\item
$[A]$ is nontrivial and the
charge of $[A]$ at south hemi-sphere is less than $\hbar$.
Also $p_i$'s are  located in north hemi-sphere and
there exists at least one $p_i$ on the equator.
\item
$[A]$ is trivial; then
\begin{equation}
x^2(p_1) +x^2(p_2)+\cdots + x^2(p_n)=\sum_{i=1}^n k_i.
\end{equation}
\end{enumerate}
\end{enumerate}
The instanton is denoted by $[[A], [\mathfrak{p}]]_\mathfrak{w}$.
\end{definition}
\begin{example} 
Suppose $k_0=0, \mathfrak{p}=(p_1,p_2)$ and $\mathfrak{w}=(1,1)$.
Then 
$$
\mathcal{M}^b_{0,\mathfrak{w}}= S^3/\mathbb{Z}_2 = RP^3.
$$
\end{example}
Directly from the definition, we know that
{\lemma
$\mathcal{M}_{k_0,\mathfrak{w}}^b$ is a smooth manifold.
}
\begin{definition}
$(T,v_0)$ is a bubble tree and $T=(V,D)$.
Suppose $Y$ is a bubble space  of
$(T,v_0)$, where
$Y = (\coprod_{v_i\in V}S_i^4)/d_T$,  and instantons $[A_i]
\in \mathcal{M}_{w(v_i)}(S^4_i)$.
For each $v_i$, suppose $child(v_i)=\{v_{i1},\ldots
v_{il_i}\}$. Let
$$
\mathfrak{p}_i=(d_T(v_i, v_{i1}), \ldots, d_T(v_i, v_{il_i}))
$$
and
$$
\mathfrak{w}_i = (W(v_{i1}, \ldots, W(v_{il_i})).
$$
Data $(Y, [A_i])_{v_i\in V}$
generate generalized instantons $[[A_i], [\mathfrak{p}_i]]
_{\mathfrak{w}_i}$
on $S^4_i$. 
If each element $[[A_i],[\mathfrak{p}_i]]_{\mathfrak{w}_i}$
is balanced, then $([[A_i],[\mathfrak{p}_i]]
_{\mathfrak{w}_i})_{v_i\in V}$
is called a \emph{ bubble tree instanton}.
The \emph{ stratum}, denoted by $\mathcal{S}_{T}(S^4)$,
is the set of such bubble tree
instantons.
\end{definition}
Let $\tilde{M}_{T}$ be the set of   bubble spaces of some
bubble tree
instantons
in $\mathcal{S}_{T}(S^4)$. Suppose $Y\in
\tilde{M}_{T}$.
Note that $Y$ can be represented as $(\coprod S^4_i, d_T)$.
For each $v\in V$, suppose $child(v)=\{v_1,\ldots, v_n\}$
and $\mathfrak{w}_v= (W(v_1),\ldots, W(v_n))$. 
The action $S_{\mathfrak{w}_v}$
on $(d_T(v,v_1),\ldots, d_T(v,v_n))$ induce an action on $d_T$
and so on $\tilde{M}_{T}(S^4)$.
Define
$M_{T}(S^4)=\tilde{M}_{T}/S_{\mathfrak{w}_v}$. Elements of
$M_{T}(S^4)$
are called  {\it bubble
tree manifolds}.

we make a convention
on representations
of bubble trees and bubble tree manifolds. 
We denote by $R(T)$   the representation of $T$. We define
it inductively: Suppose $v_0$ is the root and 
$child(v)=\{v_1,\ldots,v_k\}$. Their weights are $w_i,0\leq i\leq k$.
Then 
$$
R(T)= [w_0[w_1\cdot R(t(v_1)),w_2\cdot R(t(v_2))
,\ldots, w_k\cdot R(t(v_k))]].
$$
In this equation, $\cdot$ is only used  to make the equation
easy to read.
Similarly, we can represent a bubble tree manifold in the same way.
Suppose $Y$ is a bubble tree manifold of $T$ defined by 
$d_T$. To represent $Y$ is essentially  to represent 
$d_T$  if  $X$ is given. We denote by $Y(d_T)$ the representation
of $d_T$. Suppose $p_i= d_T(v_0,v_i), 1\leq i\leq i$,
then
$$
Y(d_T)=[p_1\cdot Y(d_{t(v_1)})  , p_2\cdot
Y(d_{t(v_2)}), \ldots, p_k\cdot Y(d_{t(v_k)}) ].
$$
Follow the notations,
one can  see that
\begin{lemma}
$M_{T}(S^4), \mathcal{S}_{T}(S^4)$ are smooth
manifolds.
\end{lemma}
\noindent
{\bf Proof: }
We prove this by induction. If $T$ consists of 
one vertex, then the lemma follows from lemma 
3.5.
Else, suppose
$$child(v_0)=\{v_1,\ldots,v_n\}
\mbox{ and } \mathfrak{w} = (W(v_1), \ldots, W(v_n)).$$
$\mathcal{S}_{T}(S^4)$
is a fiber bundle over $\mathcal{M}^b_{w(v_0),\mathfrak{w}}$
and the fiber is a product of $\mathcal{S}_{t(v_i)}$.
Then $\mathcal{S}_T(S^4)$ is smooth  if
$\mathcal{S}_{t(v_i)}$ are smooth. This can be finished
inductively. The proof for the $M_{T}(S^4)$'s is the same.
q.e.d.
\vskip 0.1in
\noindent
Define
$$
\overline{\mathcal{M}}^b_K = \cup_{T\in \mathcal{T}_K}
 \mathcal{S}_T(S^4).
$$

Now we explain how these definitions can be modified and 
generalized to general $X$.
Given a manifold $X$ and a bubble tree $(T,v_0)$,
we define  bubble tree manifolds and bubble tree instantons
inductively.
First assign the root $v_0$ for the manifold $X$. Then for any
edge $e_i=(v_0, v_i)$ in $T$, define a bubble point
$p_i:=d_T(e_i)\in X$.
For all vertices in the subtree $t(v_i)$, the assigned spheres are
$\widehat {TM}_{p_i}$,
the one-point compactification of the tangent space
$TM_{p_i}$.
{\it Bubble tree instantons} are the modified bubble
tree manifolds with instantons on each  component and
 they satisfy the  same requirements given in definition 3.6
except the one on root $X$.
The instanton on $X$ is in $\mathcal{M}_{w(v_0)}(X)$.
Again strata are defined
to be the set of bubble tree instantons. 
We denote  the stratum by $\mathcal{S}_T(X)$. $M_T(X)$
is defined similarly.
It is not  obvious that the strata defined
this way are still smooth,
since the strata are parametrized by bubble points on $X$
and we define the space  pointwise. However this is still true,
indeed it is a fiber bundle over $X^n/S_\mathfrak{w}$
 for an appropriate $\mathfrak{w}$.
This can be seen easily from a better version of description
 given in \S 3.2.
Define
$$
\overline{\mathcal{M}}_K(X) = \cup_{T\in \mathcal{T}_K}
\mathcal{S}_T(X).
$$

We now introduce a set of terminologies concerning ghost bubbles.
Suppose $T=(V,D)$ is a bubble tree and $M$ is a bubble
tree manifold of $T$. By a {\it ghost vertex}
 $v$ we mean that
$w(v)=0$. $G_T$ is the set of ghost vertices of $T$.
The corresponding component in $M$ is called a {\it 
ghost bubble}.
If $T$ contains any ghost vertex, we say $T$ is a {\it ghost bubble
tree} and $\mathcal{S}_T(X)$ is a {\it ghost stratum}.
Let $\mathcal{G}_K\subset \mathcal{T}_K$
be the set of all ghost trees. Define
$$
S_K(X)=\cup_{T\in \mathcal{G}_K} \mathcal{S}_T(X).
$$
We call $S_K(X)$ the {\it singular set} of 
$\overline{\mathcal{M}}_K(X)$. Similarly, all definitions
apply to $\overline{\mathcal{M}}^b_K$.

\subsection{The Gluing Theory}
To patch all strata together and give $\overline{\mathcal{M}}$
a manifold  (or orbifold)  structure, the standard tool is the gluing theory,
which has been developed
intensively
(\cite{DK},\cite{T}).
The basic idea is that the gluing method constructs 
local charts for the compactified space.
In this section, we
go over the gluing theory. In addition, 
we consider how to ``patch'' gluing maps for different
strata. Therefore, we are enable to 
conclude the {\it smoothness} of the whole compactified space. 
This is explained in \S 3.3 and necessary estimates are given in 
this subsection.
This type issue was first discussed by Ruan \cite{R} in the quantum
 cohomology theory.

We are going to encounter a large collection of notations 
concerning fiber
bundles in this paper. Suppose $\pi: X\to Y$ be a fiber bundle
with fiber $Z$. To make the notations more suggestive,
we write $X= Y\dot{\times} Z$ if no confusion is caused.

Given a stratum $\mathcal{S}_T(X), T\in \mathcal{T}_K$,
there is a so-called
{\it gluing parameter} $Gl_T$.
Roughly speaking,
the gluing  theory  is to study a gluing map
$$
\Psi_T: \mathcal{S}_T(X)\dot{\times} Gl_T
\to {\mathcal{M}}_K(X).
$$
Also for any $T<T'$ there is an associated gluing parameter
$Gl_{(T,T')}\subset Gl_T$.
$
\Psi_T $ maps $\mathcal{S}_T(X)\dot{\times} Gl(T,T')$
 to $\mathcal{S}_{T'}(X)$.
In the gluing construction,  a very important step is to construct
approximating solutions via  ``splicing''.
In particular, this is the case when we expect
a global gluing. This step is to
construct a splicing  map
$$
\Psi_T' : \mathcal{S}_T(X)\dot{\times} Gl(T,T')
\to \mathcal{B}_{T'}
$$
for all $T<T'$.

We  study the base case. Suppose bubble tree $T$
 consists of only two vertices $v_0, v$
and one edge $e=(v_0,v)$. The weights are $k_1=w(v_0), k_2=
w(v)$ and let $K=k_1+k_2$.  $R(T)= [k_0[k_1]]$.
Assume $k_1k_2\not= 0$.
Bubble tree manifolds are given by $d_T(e)$.
Suppose $d_T(e)=p, p\in X$, then $Y(d_T)=[p]$.
So $M_T$
is parameterized by $X$.
To be consistent with notations below, let $X_1=X, X_2=S^4$
and $p_1=p, p_2=\infty$.

The stratum $\mathcal{S}_T(X)$
is described as follows:
Let $Fr(X)$ be the frame bundle of $X$,
$P_1\to X$ and $P_2\to S^4$ be $SU(2)$ principal bundles
with $c_2(P_i)=k_i$. $E_i$ are associated vector bundles
of $P_i$.  Define
$$
\mathcal{M}_K(S^4,X)= Fr(X)\times_{SO(4)} \mathcal{M}_K(S^4),
$$
and
$$
\mathcal{M}^b_K(S^4, X) = Fr(X)\times_{SO(4)}\mathcal{M}^b_K.
$$
Then
$$
\mathcal{S}_T(X) = \mathcal{M}_{k_1}(X)\times \mathcal{M}^b_{k_2}(S^4,X)
$$
and it is also a bundle over $X$. 
The {\it gluing data} of stratum $\mathcal{S}_T(X)$ is a bundle
$$
gl: \mathbf{GL}_T=
\mathcal{M}^0_{k_1,X}\times Fr(X)
\times_{SO(3)\times SO(4)}
\mathcal{M}_{k_2}^{b,0}(S^4)
\times R^+
\to \mathcal{S}_T(X),
$$
where the action $SO(3)\times SO(4)$ on
elements in $\mathcal{M}_{k_2}^{b,0}(S^4)$
comes from the action of $SO(4)$ rotating $S^4$
and the action of $SO(3)=SU(2)/\{\pm 1\}$ on the base frame.
Note that $M_T$ is a factor of $\mathcal{S}_T(X)$.
In this case $M_T=X$.
 Locally for a $p\in X$
the geometric meaning of
the fiber $gl^{-1}(*,p)$ is
$$
Gl(p):=Hom_{SO(3)}(\mathfrak{g}_{E_1}(p), \mathfrak{g}_{E_2}(\infty))
\times R^+.
$$
This is isomorphic to $SO(3)\times R^+= R^4\setminus\{0\}/Z_2:= Gl_T
\setminus\{0\}$, where
$Gl_T = R^4/Z_2$.
We call  the parameter $R^+$ in $Gl_T$
the {\it gluing radius}.
Use $\mathbf{GL}_T(r)$ for the set in $\mathbf{GL}_T$
 with radius $\leq r$.

We say a set $U$ in a space is {\it proper} 
if its closure in the space is
compact.
The gluing theorem is now stated as
\begin{theorem}[\cite{T}]
For any open proper set $U\in \mathcal{S}_T(X)$ there exists
a small constant $\epsilon_0$, depending on $U$,
and
a gluing map $\Psi_T$
$$
\Psi_T: \mathbf{GL}_T(\epsilon_0) \cap gl^{-1}(U) \to \mathcal{M}_K(X)
$$
such that
$\Psi_T$ is a diffeomorphism.
\end{theorem}
%The rest of this section is contributed to the proof of this theorem.
The proof consists of two parts:  part 1, the
local diffeomorphism of $\Psi_T$, and  part 2,
the injectivity of $\Psi_T$. 

\vskip 0.1in
\noindent
{\bf Part I,} {\it local diffeomorphism of $\Psi_T$}
\vskip 0.1in
The idea is that
$\Psi_T$ can be constructed  gauge equivariantly. 
So we can work on $\mathcal{A}^\ast$ instead of 
$\mathcal{B}^\ast$. 
Also the local diffeomorphism is a local issue. 
We set up the
local coordinates for gluing.
Suppose that $[A_i^0],i=1,2$ are in $\mathcal{M}_{k_1}(X),
\mathcal{M}^b_{k_2}$; $U_1,U_2^b$ are their neighborhoods;
$p\in X$
and $B_p(1)$ is a ball at $p$ with radius 1.
We identify
the ball with the unit ball in $R^4$. (Readers are
referred to \cite{T} for a careful explanation.) Then locally
$\Psi_T$ is a map from
$$
\Psi_T: \mathbf{U}(\epsilon):=
U_1\times U_2^b \times B_p(1)\times SO(3)\times (0,\epsilon)
\to \mathcal{M}_K(X).
$$

It is convenient to compare this gluing map with the one 
for a more classic model. This comparison also suggests 
a reasonible treatment of $B_p(1)\times R^+$ in
the domain of $\Psi_T$.

We modify the domain of $\Psi_T$ by:
fixing  a $r\in (0,\epsilon)$ and a point $p\in X$,
replacing
$\mathcal{M}^b_{k_2}$  by $\mathcal{M}_{k_2}(X_2)$.
For this set-up, $X_2$ can be arbitrary. 
This is  the case considered in
\cite{DK}. But eventually, we are only interested in $S^4$.
Moreover identify $S^4\setminus\{\infty\}= TX_p$.
Define gluing data to be the bundle
$$
\mathbf{GL}'(p,r)=\mathcal{M}^0_{k_1}(X_1)
\times_{SO(3)} \mathcal{M}^0_{k_2}(X_2)
$$
over $\mathcal{M}_{k_1,k_2}(X_1,X_2):=\mathcal{M}_{k_1}(X_1)
\times\mathcal{M}_{k_2}(X_2) $
and gluing map $\Phi(p,r): \mathbf{GL}'(p,r)\to \mathcal{M}_K(X)$ 
which is defined later.
\begin{prop}
Given  proper sets $U$ in $\mathbf{GL}'(p,r)$,
$\Phi(p,r)$ is a local diffeomorphism   for any  $r<\epsilon_0$.
$\epsilon_0$ depends only on $U$.
\end{prop}
Locally write
$$
\Phi(p,r): U_1\times U_2 \times Gl(p)\to
\mathcal{M}_{K}(X)
$$
where $U_i\subset \mathcal{M}_{k_i}(X_i)$. We prove that
this  map, constructed in the proof, is a local diffeomorphism.
The proof
consists of four steps.
\\
{\it Step 1, splicing maps and   approximate solutions:}
Identify
$B_{p_i}(1)\setminus\{p_i\}$ with the cylinder
$R^1\times S^3$ by mapping points the polar coordinates $(s,\theta)$
to $(\log s, \theta)$.
Let $X'_i = X_i\setminus B_{p_i}(N^{-1}r) $
for some large constant $N$. Identify the 
annulus $B_{p_i}(N^{-1}, Nr):= B_{p_i}(Nr)\setminus
B_{p_i}(N^{-1}r)$ with the tube $[\log r-T, \log r +T ]
\times S^3$, where $T=\log N$.
Denote the tube by $C_{p_i}[\log r-T, \log r+T]$, or
$C_{p_i}$ if no confusion is caused.
Define the connected sum $X_r = X\sharp S^4$ by gluing
$X_i$ at $C_{p_i}$ with identification
$(t,\theta)\sim (2\log r-t, \theta)$.
$S^4$ should be naturally treated as one-point compactification
of $(TX)_p$, so we can simply assume that
$X_r$ is conformal to $X$.
Suppose $g_i$ are standard metrics of $X_i$. We 
choose a metric $g_r$ on $X_r$ with $g_r=m_ig_i$
over $X'_i$ such that $1\leq m_i\leq 2$ on the gluing
area and equal to $1$ elsewhere. The Sobolev
norms $L^{1,p}, L^q$ on $X_r$ are defined with
respect to this metric. Also we fix $p,q$ such that
$$
2 < p <4 \mbox{ and } q= \frac{4p}{4-p}.
$$
There are several cut-off functions that  are frequently used:
let $i=1,2$,
\\
$\gamma_{i,p,r}$:
equals  $0$ for $t  < \log r -1$ on $C_{p_i}$ and equals
$1$
 for $t> \log r +1$.
Moreover

\hskip 0.1in
$
\gamma_{1,p,r} + \gamma_{2,p, r} = 1
$ on $X_r$
\\
$\beta_{i,p,r}$:  equals  $0$  for 
$t < \log r - T$ on $C_{p_i}$ and
equals  $1$ for
$t> \log r -1$.
\\
$\eta_{i,p,r}$:  equals  $0$ 
 for $t < \log r + T + 1$ on $C_{p_i}$ and
equals  $1$ for
$t> \log r +2T + 1$
\\
Note that $\beta_{i,p, r}=1$ over the support of $\gamma_{i,p, r}$.
We drop  subscripts $p,r$ from these notations unless we need
to address their relations to $p,r$.
\vskip 0.1in
\noindent
Convention: in this section, all constants independent of $r$
are denoted
by $C,  C_i$ (or $\epsilon, \epsilon_i$)
 if we require them  bounded from above
(below). They may depend on proper sets $U$.
\begin{lemma}
$\beta_i$ can be chosen such that
$ \|\nabla \beta_i\|_{L^4} \leq C(\log N)^{-3/4}. $
\end{lemma}
\noindent
{\bf Proof: }
For any function $f$, $\|\nabla f\|_{L^4}$
is conformally invariant. We can construct $\beta_i$
over cylinder or $R^4$. We do this on cylinder.
Define
$$
\tilde{\beta}(t, \theta) = \left\{
\begin{array}{ll}
0  & t\leq 0 \\
1 & t\geq \log N \\
 t (\log N)^{-1} & otherwise
\end{array}
\right.
$$
$\|\nabla \tilde{\beta}\|_{L^4}$ satisfies the requirement,
However $\tilde{\beta}$ is not smooth.
Let $h(x)$ be a smooth, compact supported  function over $R^1$ such that
$\int h=1$.
Define
$$
\beta(t,\theta)= \int_{R^1}\tilde{\beta}(s,\theta)h(t-s)ds.
$$
It is easy to show that $\beta$ satisfies the estimates.
$\beta_i$ can be obtained from $\beta$ by shifting
and  reflection on $t$ coordinate. So $\nabla\beta_i$
have same estimates as $\nabla\beta$.
q.e.d.
\vskip 0.1in
\noindent
Set $\delta= C(\log N)^{-3/4}$.
The following two lemmas study  derivatives of
our cut-off functions with respect to $r,p$.
\begin{lemma}
Let $f_{r} = \gamma_{i,r}, \beta_{i,r},\mbox{ or } \eta_{i,r}$.
$$
|\frac{\partial}{\partial r} f_r|<Cr^{-1},
$$
and for $f_r=\beta_{i,r}, \eta_{i,r}$
$$
\|\nabla(
\frac{\partial}{\partial r} f_r)\|_{L^4}
<Cr^{-1}.
$$
\end{lemma}
Here the derivatives of functions on $X_2$ is interpreted as
follows:
these functions are well defined over $X_r$
and hence they induce functions  over $X$. In fact, they 
are supported in
$B_p(Nr)$.
The derivatives make sense
when these functions are treated as functions over $X$.
\vskip 0.1in
\noindent
{\bf Proof: }
We use $(x,\theta)$ to denote the coordinate
of points on $X$ and $(s,\theta), (t,\theta)$
for points on cylinder.
Suppose $f_r=\beta_{1,r}$. By definition
$\beta_{1,r}(x,\theta) = \beta(\log x-\log r +T, \theta)$.
Here $\beta$ is defined in lemma 3.10.
A direct computation shows that
$$
\frac{\partial}{\partial r}
\beta_{1,r}(x, \theta)
=
r^{-1}\int \beta(\log x-\log r -t) h'(t) dt.
$$
So
$
|\frac{\partial}{\partial r} \beta_{1,r}|
< C r^{-1}.
$
Same argument clearly works for $\gamma_1$.\\
Since $|\nabla f|_{L^4}$ is  conformally invariant,
we can compute  it
in terms of cylinder coordinates. So
$$
\nabla
(\frac{\partial}{\partial r}\beta_{1,r}(s, \theta))
= r^{-1} \int \beta'(s - \log r -t) h'(t) dt.
$$
The second inequality in the 
  lemma for $\beta_{1,r}$
 follows from this expression. Of course,  the proof
here works for $\eta_1$. \\
For $\gamma_2, \beta_2,\eta_2$ note that the derivatives
are only supported in annulus $B(N^{-1}r,Nr)$.
The proof essentially has no difference.
q.e.d.
\begin{lemma}
Let $f_{p,r} = \gamma_{i, p,r}, \beta_{i,p,r},
\mbox{ or } \eta_{i,p,r}$.
$$
|\frac{\partial}{\partial p} f_{p,r}|<Cr^{-1},
$$
and for $f_r=\beta_{i,p,r}, \eta_{i,p,r}$
$$
\|\nabla(
\frac{\partial}{\partial p} f_r)\|_{L^4}
<Cr^{-1},
$$
\end{lemma}
\noindent
{\bf Proof: }
Note that
$$
f_{q,r}(x) = f_{p,r}(|x-q|).
$$
We explain how this works for the estimates by describing 
one case.
Fix a point $p$, suppose $f_{p,r}= \beta_{1,p,r}$.
$\beta_{1,p,r}(x) = \beta(\log |x| -\log r +T)$.
Then
$
\beta_{1,q,r}=\beta(\log |x-q| -\log r +T).
$
Note that the derivative is only supported in
annulus $B_p(N^{-1}, Nr)$. So
$$
|\frac{\partial}{\partial q}|_{q=p}
\beta_{1,q,r}|
= |\beta'(\log |x-q|-\log r +T)\frac{\nabla_q|x-q|}{|x-q|}||_{q=p}
\leq C r^{-1}.
$$
Other cases are similar to the computations in lemma 3.11.
q.e.d.
\vskip 0.1in
\noindent

For any $[A_i]\in U_i$,
choosing representatives of the class is based on
\begin{lemma}
Suppose $E$ is a trivial bundle over unit
ball of $R^4$. Then there exists a representive of a connection
with $A_r = 0$, unique up to $A\to uAu^{-1}$ for a constant gauge
 transformation $u$. Moreover
$$
|A(x)|\leq |x|sup|F_A|.
$$
\end{lemma}
This is a well known result. We skip the proof.
If $A$ is chosen in this way, we say that  $A$
is a {\it r-gauge connection}
at 0.

%Suppose we fix a frame of $E_i$ and use it to trivialize
%bundles $E_i$ over $B_{p_i}(1)$. The representatives
%of $[A]$ are chosed to be r-gauge. The set of
%choices is still infinite
%dimensional. This can be  improved.
%\begin{prop}
%The $SO(3)$-principle bundle
% $\mathcal{M}^0_{k,X}$
%over $\mathcal{M}_k(X)\times X$ can be constructed such that
%connections on each fiber of $([A],p)$
%are $r$-gauge connections at $p$ in class $[A]$.
%\end{prop}
%\noindent
%{\bf Proof: }
%Indeed, $\mathcal{M}^0_{k,X}$
%is already given in \S 2.  We need to show that
%given a pair $(A,h)\in \tilde{\mathcal{M}}\times
%P|_p, p\in X$, the $\mathcal{G}_P$ orbit $[A,h]$
%of $(A,h)$ defines a unique $r$-gauge connection.
%Suppose $(A,h)$ is given.
%$A$ determines a section $\phi(A,h)$
%of $P$ over $B_p(1)$ by the parallel transport
%of $h$ along radial geodesics through $p$.
%The trivialization of $A$ with respect to
%this section of frame is $r$-gauge and
%independent of
%gauge transformations. So the proof is concluded.
%q.e.d.
%\vskip 0.1in
%\noindent
Back to 
 the construction of  approximating solutions.
Suppose that $A_i$ are r-gauge.
Define $A'_i = \eta_i A_i$.
Note that
$A'_{i}$ are trivial
on the gluing area. Then
\begin{equation}
|A_i - A'_{i } | \leq C r,  |F(A'_{i})| \leq C,
\end{equation}
and
\begin{equation}
\|F^+(A'_{i})\|_{L^2}, \|A_1- A'_{i}\|_{L^4} \leq C r^2.
\end{equation}
For any $\rho \in Gl$, one can glue bundles
$E_i$ to get a bundle $E_{\rho}$ over $X_r$ using
$\rho$.
$A'_{i}$ automatically give a connection $A'_{\rho}$.
In general, we write $\mathbf{A}^0 =
(A^0_1, A^0_2, \rho_0)$ and let $\mathbf{a}=(w_1, w_2, \nu)$
be a tangent vector  at $\mathbf{A}^0$.
\begin{prop}
$A'_{\rho}$ are approximating
 ASD-connections,
$$
\|F^+(A'_{\rho}) \|_{L^p} \leq C  r^{1+ 4/q}.
$$
\end{prop}
\noindent
{\bf  Proof. } Note that $F^+(A'_{\rho}) = F^+(A'_{1 }) 
+F^+(A'_{2 })$.
$$
\|F^+(A'_{1})\| _{L^p}\leq  \|\nabla \eta_1 A_1\|_{L^p} +
\|(\eta - \eta^2) A_1^2\|_{L^p}=: I_1 +I_2.
$$
Note that the integral area is in $B_{p_1}(Cr)$ and
$$
I_1\leq \|\nabla\eta_1\|_{L^4}\|A_1\|_{L^q}\leq C\|
\nabla\eta_1\|_{L^4}\|A_1\|_{L^p_1}.
$$
We can prove that
$I_1 < C  r^{1+ 4/q} $
 by using lemma 3.10 and (4). The estimate of $I_2$ 
also
follows from  (4). So
$$
\|F^+(A'_{1}) \|_{L^p} \leq C  r^{1+ 4/q}.
$$
Similarly, we have estimates for
$\|F^+(A'_{2})\| _{L^p}$.  q.e.d.

\begin{prop}
 Let $\mathbf {a} = 
(w_1,w_2, e)$ be  as above, then
for $z = w_i$
$$
\|\frac{\partial}{\partial z} F^+\|_{L^{p}} \leq C
r^{1+4/q}\|z\|_{L^{1,p}}.
$$
\end{prop}
\noindent
{\bf Proof. }
By direct computation,
$$
\|F^+_z\| \leq C(\|\nabla \beta z\|_{L^p} + \|(\beta-\beta^2) z\wedge
A\|_{L^p}.
$$
Note that $w_i$ are tangent vectors of a finite dimensional moduli space,
  they are bounded and $C^\infty$.
Moreover
they are  chosen so that   $w_i(p_i)=0$. So we can assume
$|w_i(x)| \leq C|x|,$ for $x\in B_1(p)$.
Now the rest of the proof
is same as the previous proposition.
q.e.d.
\vskip 0.1in
\noindent
The reason for using $r$-gauge connections is to 
keep the gluing process gauge equivariant.
\\
{\it Step 2:  the right inverse of $d^+_{A'_{\rho}}$: }
We assume that $H^2_{A_i}=0$. This is  the case that we are concerned
with
in this paper: when $X=S^4$, $H^2_A=0$ for all ASD-connections including
trivial connections; when $X\not= S^4 $, 
it follows our  assumption that the metric $g$
is regular and $\mathbf{F}_E=\emptyset$ (rf. \S 2.1).
 Therefore there are right inverses $P_i$
$$
P_i:
 \Omega^+_{X_i} (\mathfrak{g}_{E_i})
 \to \Omega^1_{X_i} (\mathfrak{g}_{E_i}),
 $$
to the  operators $d^+_{A_i}$.
$P_i$ is uniquely determined by  its image. In particular,
we can choose $P_i$ by requiring that the
image of $P_i$ is perpendicular to the
 kernel of $d_{A_i}^+$ in $L^2(\cap L^p)$.
So
$$
\|P_i \xi\|_{L^q} \leq C \|\xi\|_{L^p}.
$$
Since the moduli
 spaces $U_i\subset\mathcal{M}_i$
 are smooth, the finite
dimensional proper open sets
 and operators $d_{A_i}^+, P_i$ are smoothly parameterized by
 the manifold
$\mathcal{M}_i$,
so we have for $z= w_i$
\begin{equation}
\|\frac{\partial}{\partial z} P_i \|
\leq C \|z\|_{L^{1,p}(X_i)}
\mbox{ and }
\|\frac{\partial}{\partial z} d_{A_i}^+\|
\leq C \|z\|_{L^{1,p}(X_i)}.
\end{equation}
Define $Q_i = \beta_i P_i \gamma_i$  
and let $Q=Q_1+ Q_2$.
We have
\begin{equation}
\|d^+_{A'_{\rho}} Q \xi - \xi \|_{L^p(X)} \leq 
 C\delta \|\xi \|_{L^p(X)}.
\end{equation}
The proof can be found in \cite{DK}.
\begin{prop}
For small $r$ there exists a right
 inverse $P$ to $d^+_{A_\rho}$ over $X_r$ satisfying
\begin{equation}
\|P\xi\|_{L^q}\leq C \|\xi\|_{L^p}.
\end{equation}
 Moreover for $w_i\in
(T\mathcal{M}_i)_{A_i}$
\begin{equation}
\|\frac{\partial}{\partial w_i} P\xi\|_{L^{1,p}}
\leq C \|w_i\|_{L^{1,p}(X_i)} \|\xi\|_{L^p}.
\end{equation}
\end{prop}
\noindent
{\bf Proof:} By (7) , 
we know that $d^+_{A'_\rho} Q$ is invertible if $\delta$
is small enough.
Letting $r$ small, we 
can choose a constant $N$ to let $\delta$ small. (see lemma
3.10). 
  Set
$P= Q(d^+_{A'_\rho} Q)^{-1}$,  (8) is obvious then.
 By definition of  $Q_i$ and (6)  
$$
\|\frac{\partial}{\partial z} Q_i \xi\|_{L^{1,p}(X)}
\leq C \|z\|_{L^{1,p}(X_i)}\|\xi\|_{L^p(X)}.
$$
Same estimate holds for $Q$.
(9) follows from the following expression of
 the derivative of $P$ and existing estimates for $Q$: 
$$
P' = (Q(d^+_{A'_\rho} Q)^{-1})' = Q'(d^+_{A'_\rho}Q)^{-1} +
P(d^+_{A'_\rho}Q)'(d^+_{A'_\rho}Q)
^{-1}.
$$
This finishes the proof. q.e.d.

\vskip 0.1in
\noindent
{\it Step 3, constructing solutions: } Here we use Taubes' argument.
It is said that all solutions are of the form $A_\rho =
A'_\rho + a$, where $a = P(\xi)$ for some $\xi$.
Then equation $F^+(A_\rho)=0$  changes to
\begin{equation}
d^+_{A'_\rho} a  + a\wedge a = - F^+(A'_\rho).
\end{equation}
Apply $P$ to the equation,
\begin{equation}
a + P(a\wedge a) = P(-F^+(A'_\rho)).
\end{equation}
\begin{prop}
Let $\delta_0 >0$ be a small constant
depending on $U_1\times U_2$.
For small  $r\leq \delta_0$ equation (11) has unique solution with
$$
\|a\|_{L^{1,p}} \leq C r^{1+3/q}.
$$
Let $w_i$ as before
\begin{equation}
\|\frac{\partial}{\partial w_i} a\|_{L^{1,p}} \leq C r^{1+3/q}
\|w_i\|_{L^{1,p}},
\end{equation}
\end{prop}
\noindent
{\bf Proof:}
recall that  $\| F^+(A'_\rho)\|_{L^p}\leq C r^{1+4/q}$.
By (8), the same estimate holds for $P(-F^+(A'_\rho))$.
Let
$$
B =\{ f\in L^{1,p}(\Omega^1_X(\mathfrak{g}_{E_\rho}))| 
\|f\| \leq C
r^{1+3/q}\}.
$$
Define a map $H: L^{1,p}\to L^{1,p}$ by
$$
H(a) =  -P(a\wedge a) - P(F^+(A'_\rho)).
$$
It is routine to show that
 $H$ maps $B$ to $B$ and it satisfies  the contraction
mapping
principle(\cite{DK}).  This proves the first statement.

Taking derivative on equation (11),  we have
$$
a' =- 2P(a\wedge a')- P'(a\wedge a) + P' 
(-F^+(A'_\rho)) +P((-F^+(A'_\rho)').
$$
Using the estimates for $P', (F^+(A'_\rho))'$,
we have bound for
$L^{1,p}$ norms
of  the last three terms.  For the first term
\begin{eqnarray*}
\|P(a\wedge a')\|_{L^{1,p}} &\leq& C\|a\wedge a'\|_{L^p}
\leq C \|a\|_{L^4} \|a'\|_{L^q} \\
&\leq& C \|a\|_{L^{1,p}}\|a'\|_{L^{1,p}}
\leq C r^{1+3/q} \|a'\|_{L^{1,p}}\leq \frac{1}{2}\|a'\|_{L^{1,p}}.
\end{eqnarray*}
Absorb this term on the left side, we have (12).
q.e.d.

\vskip 0.1in
The map $\Phi$ now can be defined as
$$
\Phi(A_1,A_2,\rho)= A'_\rho + a,
$$
where $a$ is solved in proposition 3.17
 and it can be treated as a map defined on
$U_1\times U_2\times Gl$. Define $\Phi'(A_1, A_2,\rho) = A'_\rho$
to be the splicing map.
So
$\Phi = \Phi' +a$.
\begin{corollary}
For any $z = w_i\in T\mathcal{M}_{k_i}(A_i)$
or  $z = \nu \in T_1(Gl)$
\begin{equation}
\|da (z)\|_{L^{1,p}} \leq C r^{1+3/q}
\|d\Phi' (z)\|_{L^{1,p}}
\end{equation}
Namely $\|da(z)\|\lll \|d\Phi'(z)\|$.
\end{corollary}
\noindent
{\bf Proof: } By the definition of  the metric $g_r$ of $X_r$,
it is obvious that
$C'\|w_i\|\leq \|d\Phi'(w_i)\|\leq C\|w_i\|$. 
So  (13) is a consequence of (12).
This result is still true for the tangent vector $\nu$.
A simple way to see this is to replace $\mathcal{M}_{k_2}(X_2)$ by
$\mathcal{M}^0_{k_2}(X_2)$ ($\cong Gl\times \mathcal{M}_{k_2}$ locally)
 and
repeat the  gluing theory.
Now $Gl$ is embedded in $\mathcal{M}^0_{k_2}$
as a subspace. $\nu$ can be treated as  a vector
 in $\mathcal{M}^0_{k_2}$
and so by the same argument we have (12) for $z= \nu$ and so (13).
q.e.d.
\vskip 0.1in
\noindent

Since we are concerned about the property of derivatives which is 
a local computation,
locally we can  put $Gl\times U_2$ together and replace it by $U^0_2\subset
\mathcal{M}^{0}_{k_2}$ from now on. This 
is just what we did in the proof. 
Note that $d\Phi' \sim Id$, so the corollary implies that
$\|da\|\leq C r^{1+3/q}$.

To study the diffeomorphism issue, it is useful to give a local 
coordinate chart for $\mathcal{M}_E$. The next few pages (up to 
remark 3.24) is contributed to this issue.
 
We begin with the study of  the local structure of $\mathcal{A}$ around 
$Im(\Phi')$.
Locally, let $W= \Phi'(U_1\times U_2^0)$. Let $v_0\in W$
be a fixed connection. $\mathcal{A}_E$ is indentified as
a vector space with $v_0=0$.
Define a map
$$
\mathcal{Y}: \Omega^0(\mathfrak{g}_E)\times
W \times \Omega^2_+\to \mathcal{A}_E
$$
given by
$$
\mathcal{Y}(\xi, v, \eta)=
\exp(\xi)(v+ P_v \eta), 
$$
where $P_v$ is the right inverse to $d^+_v$.
In the definition, $\Omega^0(\mathfrak{g}_E)$ is assigned
$L^{2,p}(X_r)$-norm;
$\Omega^2_+$ is assigned $L^p(X_r)$-norm;
$V$ and $\mathcal{A}_E$ are assigned $L^{1,p}$-norms.
We would like to conclude that $\mathcal{Y}$
is a diffeomorphism in the neighborhood of $V$. 
\begin{prop}
Let $V_0= (0,v,0)$. Let $T= D\mathcal{Y}_{V_0}$ be the tangent map.
$$
T(\xi, z, \eta)= 
d_{v} \xi + z+ P_{v}\eta.
$$ 
Then $T$ is an isomorphism and
$$ 
\|T\|\leq C, \|T^{-1}\|\leq C.
$$
\end{prop}
\noindent
{\bf Proof: }
$\|T\|\leq C$ follows from the definition of norms.
Conversely, we first show that
$$
\|(\xi,z,\eta)\|\leq C\|\alpha\|,
$$
where $\alpha= T(\xi,z,\eta)$.

We know
$$
d_v\xi + z + P_v\eta=\alpha,
$$
apply $d_v^+$ to the equation,
$$
d_v^+d_v\xi + d_v^+z +\eta= d_v^+ \alpha.
$$
It is easy to have
\begin{eqnarray*}
\|d_v^+\alpha\| &\leq&  C\|\alpha\| \\
\|d_v^+z\|&\leq& C r \| z\| \\
\|d_v^+d_v \xi\| &\leq & Cr \|d_v\xi\|.
\end{eqnarray*}
On the other hand, similar to the argument in \cite{DK},
by the fact that 
$H^0(A)=H^0(B)=0$ (where $v=\Phi'(A,B)$), we have 
$$
\|d_v\xi +z\|\geq C(\|d_v\xi\|+\|z\|).
$$
Combine these inequalities together
\begin{eqnarray*}
\|\eta\|
&=& \|d^+_v\alpha - d_v^+d_v\xi -d_v^+ z\|\\
&\leq & 
C\|\alpha\| + C r(\|d_v\xi\|+ \|z\|) \\
&\leq &
C\|\alpha\| + C' r(\|d_v\xi+ z\|) \\
&=& 
C\|\alpha\| + C' r(\|\alpha- \eta\|).
\end{eqnarray*}
Reorganize the terms of $\|\eta\|$, we have
$$
\|\eta\|\leq C\|\alpha\|.
$$
For $\xi,z$, the estimates follow from
$$
\|d_v\xi + z\| =\|\alpha - P_v\eta\|\leq C \|\alpha\|.
$$ 
To see that $T^{-1}$ exists, one can use the argument 
of using index theory (\cite{DK}).
q.e.d.
\vskip 0.1in
In general,
\begin{prop}
For any $V= (\xi_0,v,\eta_0)$ in the $\epsilon$-ball of
$V_0$, $D\mathcal{Y}_V$ is isomorphic. Moreover,
$$
\|D\mathcal{Y}_V\|\leq C, 
\|D\mathcal{Y}_V^{-1}\|\leq C,
$$
where $\epsilon, C$ are independent of $r$.
\end{prop}
The idea is to compare the operators
$D\mathcal{Y}_V$ and $D\mathcal{Y}_{V_0}$
by a direct computation. A similar computation is already given in \cite{DK}.
One can find that the difference between two operators are  controlled
by $\epsilon$. Hence proposition 3.19 implies this propostion.

As a consequence,
\begin{corollary}
$\mathcal{Y}$ is a local diffeomorphism around $W$. 
$\mathcal{Y}(W\times \Omega^2_+) $ gives a slice for $\mathcal{A}_E$.
\end{corollary}
By the proposition 3.17, 
we know that the moduli space  $\mathcal{M}_E$ near
$W$ is isomorphic to $W$. In fact, the isomorphic map is 
given by
$$
\hat{a}(w) = w+ a\circ (\Phi')^{-1} (w), w\in W.
$$
So 
\begin{corollary}
$(W,\hat{a})$ is a local coordinate chart for $\mathcal{M}_E$.
\end{corollary}
\vskip 0.1in
\noindent
{\bf Proof of proposition 3.9: }
Clearly, $\Phi'$ is a local diffeomorphic map.
 Using the coordinate chart explained 
above, $\Phi$ itself is $\Phi'$. So $\Phi$ is 
a local diffeomorphism. q.e.d.

\vskip 0.1in
Now
suppose we 
have  a smooth map
 $\tilde{f}: W\to \tilde{\mathcal{M}}_E\subset \mathcal{A}_E$.
Treating that $W$ as a chart of $\mathcal{M}_E$,
$\tilde{f}$ induces a map from $f= \tilde{f} \hat{a}^{-1}$
 from $\mathcal{M}_E$ to itself. 
\begin{prop}
Let 
$
\hat{f}(x)= \tilde{f}(x)- x.
$ Suppose
$$
\|\hat{f}\|\leq \epsilon, 
\|d\hat{f}\| \leq \epsilon_0,
$$
for some small constant $\epsilon_0>0$.
Then $f$ is a diffeomorphism. 
Here the norm used for $\mathcal{A}_E$ is $L^{1,p}$.
Furthermore, if $\tilde{f}_i,i=1,2$ are two such maps,
then $f_2\circ f_1^{-1}$ is locally diffeomorphic. 
\end{prop}
\noindent
{\bf Proof: }
Let $v\in W$.
We may assume that 
$W$ is a linear vector space. Otherwise, we can identify 
$W$ with the tangent space $W_{v}$. Assuming $v=0$ in $W$.
By $\mathcal{Y}$, we can assume that 
$\mathcal{A}_E$ is $\Omega(\mathfrak{g}_E)\times W\times
\Omega^2_+$. Under this identification, suppose
$$
\tilde{f}(z)= (\xi(z), v(z), \eta(z)).
$$
$f$ can be interpretted as $f(z)= v(z)$  when
$W$ is identified with $\mathcal{M}_E$. By proposition 3.19, we see that
$$
\|dv- I\|\leq C\|d\hat{f}\|\leq C\epsilon_0.
$$
So when $\epsilon_0$ is chosen small such  that 
$C\epsilon_0< 0.5$, $dv$ is an isomorphism. Therefore $f$
is locally diffeomorphic. q.e.d.
\begin{remark}
The proposition 3.23 is useful for local diffeomorphism issue. 
It says that if we can control 
difference of two maps $\tilde{f}_1-\tilde{f}_2$
up to its derivatives in terms of the norm of Banach space(!),
we can conclude the diffeomorphism on the moduli spaces.
One can see that this argument works because we are taking 
particular norms on the Banach spaces. Roughly speaking, one of the key
is that
our norms are with respect to the
metric of $X_r$ instead of $X$.
\end{remark}

We now discuss our gluing map $\Psi_T$.
By comparing to $\Phi$ and so following the remark 3.24,
 we are able to show that  it is a 
local diffeomorphism.
The main difference between $\Psi_T$ and $\Phi$ 
is that for $\Psi_T$
we allow the  gluing
radius $r$ and the bubble point $p\in X$ to vary.
In order to make $\Psi_T$ diffeomorphic, we know just by dimension
counting that we  have to
replace $\mathcal{M}_{k_2}(S^4)$ by $\mathcal{M}^b_{k_2}$.
Locally, the map $\Psi_T$ is defined to be
$$
\Psi_T(A_1,A_2, \rho, p, r) = \Phi_{p,r} (A_1,A_2, \rho).
$$
This is well defined globally (\cite{T}).
Also,
 by our convention,
locally we combine $(A_2,\rho)$ together to be an $A_2\in U_2^{b,0}$. Let
$$
\Psi_T'(*,p,r) = \Phi'_{p,r}(*).
$$
Fix a point $\mathbf{B^0}= (A_1^0,A_2^0, r_0, p)$.
Let $V_1(\epsilon_1)$ and
 $V_2^{b,0}(\epsilon_1)$
be $\epsilon_1$ neighborhood of
$A_i^0$ in $\mathcal{M}_{k_1}(X_1)$
and $ \mathcal{M}^{b,0}_{k_2}$.
We construct a natural  smooth map
$
D_{p,r_0}: V_1(\epsilon_1)\times V_2^{b,0}(\epsilon_1)
 \to \mathcal{M}_{k_1}(X)
\times \mathcal{M}_{k_2}^0\times R^+\times B_p(1)
$:
Let
$\mathfrak{h}:\mathcal{M}^{0}_{k_2}(X_2)\to \mathcal{M}
^{b,0}_{k_2}$ be the projection map and
$ V_2^0= \mathfrak{h}^{-1}(V_2^{b,0})$.
For any balanced ASD-connection
$A_2\in V_2^{b,0}(\epsilon_1)$, the fiber
 $\mathfrak{h}^{-1} A$ is in the form
of $A_2(t(\cdot - y))$, where $t\in [1-\epsilon,1+\epsilon],
y\in B(\epsilon)$.
Define
$$
D_{p,r_0}(A_1,A_2(t(\cdot- y))= (A_1,
A_2, r_0t^{-1}, p+ r_0^2 y).
$$
By the geometric meanings of $t, y, r, q$, it is clear that
$\Psi_T\circ D_{p,r_0}\sim \Phi_{p,r_0}$ over $V_1\times V_2^0
(\epsilon)$.
\begin{prop}
Let $\tilde{\Psi}_T=\Psi_T\circ D_{(p,r_0)}$.
$\tilde{\Psi}_T$ is diffeomorphic over
$V_1(\epsilon)\times V_2^0(\epsilon)$ when $\epsilon$ is small.
\end{prop}
\noindent
{\bf Proof: }
Let $\Psi_\delta = \tilde{\Psi}_T - \Phi_{p,r_0}$.
We know that
$\Phi_{p,r_0}$ is diffeomorphic by proposition ?. 
By proposition ?. It is sufficient
to show lemma ?.
q.e.d.
\begin{lemma}
$\|d\Psi_\delta\| \leq C r_0^{3/q}\|d\Phi'_{p,r_0}\|.$
\end{lemma}
\noindent
{\bf Proof: }
Let $\mathbf{A^0} = D_{p,r_0}^{-1}\mathbf{B^0}$. Suppose
that $z$ is a tangent vector  at $\mathbf{A^0}$ and
$z' = dD_{p,r_0}( z)$.
Let $\tilde{\Psi}'_T = \Psi_T'\circ D_{p,r_0} $.
The difference $\Psi_\delta$
is generated by two terms: 1)
$\tilde{\Psi}'_T-\Phi'_{p,r_0}$  and
2) $\tilde{\Psi}_T-\Phi'_{p,r_0}$.
When  tangent vector $z$ is in
$(TV_1)_ {A_1^0}\times (TV_2^{b,0})_{A_2^0}$, $\Psi_\delta=0$.
So only $r$-direction vectors and
$p$-direction vectors are nontrivial.
To simplify the notations, we
consider them separately. The case for combinations of
different directions can be proved by the same arguments.
\vskip 0.1in
\noindent
Case 1: $z' = r_0\frac{\partial}{\partial r}$.\\
Here with factor $r_0$, $\frac{\partial}{\partial r}$ 
is normarlized such that $\|z\|\sim 1$.
The proof consists of two lemmas.
\begin{lemma}
$\|d\tilde{\Psi}_T'(z)-
d\Phi_{p,r_0}'(z)
\|\leq Cr^{4/q}_0.
$
\end{lemma}
\noindent
{\bf Proof: }
Suppose $A_r$ is the path  representing
$z$ at  $A_{r_0} = A_2$. And
$$
A_r = A_2(\frac{r_0}{r}\cdot).
$$
By definition $\tilde{\Psi}'_T(A_1,A_r)= \Phi_{p,r}(A_1,A_2)$.
On $X$
$$
\tilde{\Psi}_T'(A_1, A_r)
= \eta_{1,r}(A_1) + \eta_{2,r}A_r.
$$
So
$$
\frac{\partial}{\partial r}|_{r=r_0}
\tilde{\Psi}_T'(A_1,A_r)
=
\frac{\partial}{\partial r}|_{r=r_0} \eta_{1,r} A_1
+ \frac{\partial}{\partial r}|_{r=r_0} \eta_{2,r}
A_2 + \eta_{2,r_0}z.
$$
Hence
$$
\|d\tilde{\Psi}'_T(z)-d\Phi_{p,r_0}(z)\|_{L^{1,p}}
\leq \|\frac{\partial}{\partial r}|_{r=r_0} \eta_{1,r} A_1\|_{L^{1,p}}
+
\|\frac{\partial}{\partial r}|_{r=r_0} \eta_{2,r} A_2\|_{L^{1,p}}
$$
Using lemma 3.11, the estimate  is routine.
q.e.d.
\vskip 0.1in
\noindent
Similarly we  also have
\begin{lemma}
\begin{eqnarray*}
\|\frac{\partial}{\partial z} F^+(\tilde{\Psi}'_T(A_1,
A_r))\|_{L^p} \leq C r^{4/q} \\
\|\frac{\partial}{\partial z} Q(\tilde{\Psi}'_T(A_1,
A_r))\|_{L^{1,p}} \leq C r^{-1} \\
\|\frac{\partial}{\partial z} P(\tilde{\Psi}'_T(A_1,
A_r))\|_{L^{1,p}} \leq C r^{-1}
\end{eqnarray*}
and
\begin{equation}
\|\frac{\partial}{\partial z} \Psi_\delta(A_1,
A_r)\|_{L^{1,p}} \leq C r^{3/q}
\end{equation}
\end{lemma}
\noindent
{\bf Proof: }
Most of the proofs are similar to the previous one.
We explain the proof for $\frac{\partial}{\partial z} Q$.
$Q=Q_1+Q_2$. We prove for $Q_2$. $Q_1$ is similar but
slightly simpler. By definition
$Q_2(r) = \beta_{2,r}P_{A_r}\gamma_{2,r}$. For any $\xi$
\begin{eqnarray*}
\|\frac{\partial}{\partial z}
Q_2(r)\xi\|_{L^{1,p}}
&\leq&
\|\frac{\partial}{\partial r}\beta_{2,r}P_{A_2}
\gamma_{2,r_0}\xi\|_{L^{1,p}}
+
\|\beta_{2,r_0}\frac{\partial}{\partial z}P_{A_r}
\gamma_{2,r_0}\xi\|_{L^{1,p}}  \\
&+&
\|\beta_{2,r_0}P_{A_2}\frac{\partial}{\partial z}\gamma
_{2,r}\xi\|_{L^{1,p}}
=:  I_1+I_2+I_3.
\end{eqnarray*}
Estimate $I_1$:
\begin{eqnarray*}
I_1
&\leq &
 \|\nabla(\frac{\partial}{\partial r} \beta_{2,r})
\|_{L^4}\|P_{A_2} \gamma_{2,r_0}\xi\|_{L^{1,p}}
+
|\frac{\partial}{\partial r}\beta_{2,r}|
\|P_{A_2}\gamma_{2,r_0}\xi\|_{L^{1,p}} \\
& \leq &
(Cr^{-1})\|\xi\|_{L^p} \leq Cr^{-1} \|\xi\|_{L^p}.
\end{eqnarray*}
Here we use lemma 3.11. \\
Estimate $I_2$:
$$
I_2
\leq
\|\nabla \beta_{2,r_0}\|_{L^4}
\|\frac{\partial}{\partial z} P_{A_{r}}
(\gamma_{2,r_0} \xi)\|_{L^{1,p}}
\leq
C\|\xi\|_{L^{p}}.
$$
The estimate of $I_3$ is similar.
\\
For $P$ we use the expression in proposition
3.16  and the estimates for $Q$.
\\
(14) is  a consequence of  the expression of $a$ (see (11))
 and
all estimates we just have.
q.e.d.
\vskip 0.1in
\noindent
Combining these two lemmas, we prove the case 1.
\vskip 0.1in
\noindent
Case 2: $z'= r_0^2 \frac{\partial}{\partial x}$.\\
Here $\partial/\partial x$
is a unit vector in $TM_p$. The factor $r^2_0$
normalizes the vector $z$ such that $\|z\|$
is some constant $\sim 1$. This is
due to the definition
of $D_{p,r_0}$.
\begin{lemma}
$\|d\tilde{\Psi}'_T(z)-
d\Phi'_{p,r_0}(z)
\|\leq Cr^{4/q}_0.$
\end{lemma}
\noindent
{\bf Proof: }
Let $p(t)=p+tx$ be a path in $M$.
Besides the changes of cut-off functions
in the construction of gluing maps, there is another subtle
change. Recall  the map $\Phi'_{p,r_0}(A_1,A_2)$.
We require that $A_1$ is r-gauge with respect to
$p(t)$ before  applying $\Phi'$ to it. Denote
connections as $A_{1,p(t)}$. Set
$$
\tilde{z}= \frac{\partial}{\partial t}A_{1,p(t)}.
$$
This vector only depends on $A_1$ and $p$.
So $|\tilde{z}|_{L^{1,p}}\leq C$.
Let $(A_1, A_{2,t})= D_{p,r_0}^{-1}(A_1, A_2,p(t))$.
By definition,
$$
\tilde{\Psi}'_T(A_1,A_{2,r^2_0 t})
=
\eta_{1, p(r^2_0 t), r_0}A_{1,p(r^2_0 t)} +
\eta_{2, p(r^2_0 t), r_0}A_2(\cdot - tx).
$$
Take derivative with respect to $t$:
$$
\frac{\partial}{\partial t}
\tilde{\Psi}'_T(A_1,A_{2,r^2_0 t})
=
r^2_0 (\frac{\partial}{\partial p}\eta_{1,p,r_0} A_1
+ \eta_{1,p, r_0} \tilde{z}
+\frac{\partial}{\partial p}\eta_{2,p,r_0} A_2
)
+\eta_{2,p,r_0} z
$$
So
$$
|(d\tilde{\Psi}_T' - d\Phi_{p,r_0})z|
\leq r_0^2
|\frac{\partial}{\partial p}\eta_{1,p,r_0} A_1
+ \eta_{1,p, r_0} \tilde{z}
+\frac{\partial}{\partial p}\eta_{2,p,r_0} A_2
|
$$
Because of
the factor $r_0^2$, we  indeed get better
estimates than the one stated in the lemma.
q.e.d.
\vskip 0.1in
\noindent
Similarly, we also have
\begin{lemma}
$\|\frac{\partial}{\partial z} \Psi_\delta(A_{1,p(t)},
A_2)\|_{L^{1,p}} \leq C r^{3/q}$.
\end{lemma}
Combining these two cases, we prove lemma 3.26. q.e.d.
\vskip 0.1in
\noindent
{\bf Part II: } The gluing map is injective.
\vskip 0.1in
\noindent
Here
For simplicity of the notation, we prove the injectivity for 
$\Phi$.   We state the result as
\begin{prop}
$\Phi$ is injective. 
\end{prop}
{\bf Proof: }
Suppose there exists two points in $\mathcal{M}_E$
$$
B= \Phi (A_1,A_2,\gamma); B'= \Phi(A_1',A_2',\gamma')
$$
such that $[B]=[B']$, i.e, there exists
$g\in \mathcal{G}_E$ 
$$
B=gB'.
$$
Assume $K=\Phi'(A_1,A_2,\gamma)$ and $K'=\Phi'(A_1',A_2',\gamma')$. 
Without the loss of generality, we may 
assume that $L^{1,p}$-norms of $A_i,A_i'$'s are bounded. 
By proposition 3.17, 
$$
\|K-gK'\|\leq \epsilon.
$$
Now, we want to find $(A''_1,A''_2,\gamma'')$  such that for 
$K''=\Phi'(A_1'',A_2'',\gamma'')$
$[K']=[K'']$ and 
$$
\|K-K''\|\leq C\epsilon.
$$
This can be done as following: 
since $K,K'$ 's norm are small at the gluing area, so is $\|dg\|$.
Now fix a point $x_0$ in the gluing area, for example, 
$|x_0-p|=r_0$. 
Let $g_0=g(x_0)$. 
we can choose $g_1=g$ over $X\setminus B_p(r^{1/2}_0)$
and be constant $g_0$
in $B_p(2r_0)$. Similarly, (on the $S^4$ side), let
$g_2=g$ over $B_p(r^2_0)$ and be constant $g_0$ 
outside  $B_p(r_0/2)$. One can show that the choice 
$$
A''_1= g_1A_1', A''_2= g_2A_2',
\gamma'' = g_0^{-1} \gamma' g_0
$$
satisfies the requirement. Let $B''=\Phi(A_1'',A_2'',\gamma'')$.
We know that $[B]=[B'']$ and 
$$
\|B-B''\|\leq C\epsilon.
$$
However, by proposition 3.21, they are located in the local diffeomorphic
region. So it provides the injectivity. q.e.d.
\vskip 0.1in
\noindent
Modulo some trivial consideration, the main idea of the proof for $\Psi$
is identical. We skip the proof here. So we complete the proof of 
proposition 3.8. q.e.d.

Now suppose $k_2=0$. Then $\mathcal{M}^b_0=\{0\}$,
%$\mathcal{M}^{b,0}_0=SO(3)$.
$\mathcal{S}_T(X)=\mathcal{M}_{k_1}(X)$.
$\mathbf{GL}_T$ has a  natural $SO(3)$ action on itself.
In fact, $\mathbf{GL}_T = \mathcal{M}^0_{k_1}(X)$.
The gluing theory does not make sense for $\Psi_T$
for this case. The reason is that $T$ is not a bubble tree.
But a similar result to proposition 3.9 can be obtained:
\vskip 0.1in
\noindent
{\bf Proposition 3.9': }{\it
For any open proper set $U\subset \mathcal{M}_{k_1}(X)$
there exists a small constant $\epsilon_0$, depending
on $U$, and a gluing map
$$
\Phi: (\mathbf{GL}'(p)\cap gl^{-1}(U))/SO(3) \to
\mathcal{M}_{k_1}(X).
$$
such that
$\Phi$ is diffeomorphic.
}
\vskip 0.1in
\noindent
$SO(3)$ is the isotropy group of the  trivial connection $A_2$. 
The treatment for nontrivial isotropy group is
already given as in \cite{DK}. 
This is the obstruction of constructing a smooth compactification
via the bubble tree methods.

We now discuss the general cases.
Suppose a bubble tree $(T, v_0)$ is given.
We describe $\mathcal{S}^b_T$, $\mathcal{S}_T(X)$ and define
 $\mathcal{S}^{b,0}_T$.
This is done inductively.
The base case is that $T$ has only one vertex $v_0$.
Then $\mathcal{S}^b_T = \mathcal{M}^b_{w(v_0)}$,
$\mathcal{S}_T(X)= \mathcal{M}_{w(v_0)}$ and
$\mathcal{S}^{b,0}_T:=\mathcal{M}^{b,0}_{w(v_0)}$.
%%%herehere
We first assume that $T$ is not a ghost tree.
\vskip 0.05in
\noindent
{\it Case 1: }$\mathcal{S}^b_T$.\\
Let $v$ be a vertex of $T$.
Let  $child(v)
=\{v_1,\ldots, v_n\}$, and
Suppose that $w_i=w(v_i)$. Let
$\mathfrak{w}_v=(w_1,\ldots, w_n)$.
Suppose that
for all subtrees $(t(v_i), v_i),$
$\mathcal{S}^b_{t(v_i)}$ has been constructed.
Let $S^4_v$ be the sphere assigned to $v$. Then
$$
\mathcal{S}_{t(v)}= \mathcal{M}^b_{w(v),\mathfrak{w}(v)}
\times \prod_{i=1}^n \mathcal{S}_{t(v_i)}.
$$
Let $P$ be an $SU(2)$-principal bundle over $S^4_v$
with $c_2(P)=w(v)$.
Define bundle
$$
\mathbb{P}^0_v=
[\tilde{\mathcal{M}}^b_{w(v)}\times_{\mathcal{G}_P}
 (P|_\infty\times\prod_{i\in child(v)}
P)\to \mathcal{M}^b_{w(v)}\times (S^4_v\setminus\{\infty\})^n]
/S_{\mathfrak{w}_v}.
$$
Note that $\mathcal{M}^b_{w(v),W(v)}$
is a subspace
of $\mathcal{M}^b_{w(v)}\times (S^4_v\setminus\{\infty\})^n/
S_{\mathfrak{w}_v}$.
$\mathbb{P}^0_v$
defines a bundle over $\mathcal{M}^b_{w(v),W(v)}$.
$$
\mathcal{S}^{b,0}_{t(v)}
= \mathbb{P}^0_v \times_{\prod_{v_i\in child(v)} SO(3)}
(\prod_{v_i\in child(v)}\mathcal{S}^{b,0}_{t(v_i)}\times R^+).
$$
Drop $P|_\infty$ from $\mathbb{P}^0_v$ and denote
 the  bundle by  $\mathbb{P}_v$.
Then we define 
$$
\mathbf{GL}_{T}
= \mathbb{P}_{v_0}
\times_{\prod_{v_i\in child(v_0)}SO(3)}
(\prod_{v_i\in child(v_0)}(\mathcal{S}^{b,0}_{t(v_i)}\times R^+)).
$$
{\it Case 2: } $\mathcal{S}_T(X)$.\\
Suppose $child(v_0)=\{v_1,\ldots, v_n\}$ and
$\mathfrak{w}_{v_0}=(W(v_1), \ldots, W(v_n))$.
For each $v_i$, we already have
$\mathcal{S}^b_{t(v_i)}, \mathcal{S}^{b,0}_{t(v_i)}$.
$X$ is assigned to $v_0$ and each edge $(v_0,v_i)$ is assigned
a bubble point $p_i$.  Tuples
$(p_1,\ldots,p_n)$ are parameterized by
$B_{v_0}:=(X^n\setminus\Delta)/S_{\mathfrak{w}_{v_0}}$, 
where $\Delta$
is the big diagonal.
Let $P\to X$
be the $SU(2)$ principle bundle
over $X$ with $c_2(P)=w(v_0)$.
Define
$$
P_i
= Fr(X)\times_{SO(4)} \mathcal{S}^b_{t(v_i)}
\to
X
$$
and set
$$
\mathcal{S}_T(X)
=[\mathcal{M}_{w(v_0)}(X)
\times (\prod_{v_i\in child(v_0)} P_i)
\to
(\prod_{v_i\in child(v_0)} X)]/S_{\mathfrak{w}_{v_0}}.
$$
Also, define
$$
\mathbb{P}_{v_0}
=\tilde{\mathcal{M}}_{w(v_0)}
\times_{\mathcal{G}_P}
\prod_{v_i\in child(v_0)} (P \times Fr(X))
\to
\mathcal{M}_{w(v_0)}
\times \prod_{v_i\in child(v_0)}X.
$$
The bundle of gluing data is 
$$
\mathbf{GL}_T
=[ \mathbb{P}_{v_0}
\times_{\prod_{v_i\in child(v_0)}(SO(3)\times SO(4))}
(\prod_{v_i\in child(v_0)}
 (\mathcal{S}^{b,0}_{t(v_i)}\times R^+))]/S_{\mathfrak{w}_{v_0}}.
$$
Define $Gl_T$ to be the fiber of $\mathbf{GL}_T$.
By definition of $\mathbf{GL}_T$,
$$
Gl_T \cong \prod_{e\in D} (R^+\times SO(3)),
$$
where $D$ is the edge set of $T$. Namely, each
edge $e$ corresponds to a gluing parameter
$R^+\times SO(3)$ and we write this $SO(3)$ to be
$SO(3)_e$.  Treat $SO(3)$ as $SU(2)/\mathbb{Z}_2$, we get
$R^+\times SO(3)= R^4\setminus\{0\}/\mathbb{Z}_2$.
Compactify $\mathbf{GL}_T$ by adding $0$-section
to the bundle, it is then a $R^4/\mathbb{Z}_2$
bundle.  We denote
$\overline{\mathbf{GL}}_T$ for these compactified bundles.
However, for simplicity, we still use ${\mathbf{GL}}_T$
instead of $\overline{\mathbf{GL}}_T$. Including 0-section or not can
be distincted by the contexts.  
Correspondingly, $R^4/\mathbb{Z}_2$
generated by the edge $e$ is written
as $R^4_e/\mathbb{Z}_2$.
For simplicity,
at this moment, we do not want to be bothered
by the  stabilizer $S_\mathfrak{w}$'s. 
For the time being, we assume that all $S_\mathfrak{w}=1$.
Given a vertex $v$ in $T$, suppose $e_1,\ldots, e_n$
 are edges connecting $v$.
Define  $\mathbf{GL}(v)$ to be a sub-bundle
of $\mathbf{GL}_T$ with fiber isomorphic to
$\prod_{i} (R^4_{e_i}/\mathbb{Z}_2)$.

Suppose $T$ is a ghost bubble tree. One can still define
a bundle $\mathbf{GL}_T$ as above. But there are more structures
on this bundle. Suppose $v$ is a ghost vertex. There is
an $SO(3)$ action on the fiber of $\mathbf{GL}(v)$. The action is
along fiber and it
 acts on
each $SO(3)_e$ as group multiplications on suitable side.
%The action is similar to
%the base case that we explained. 
Hence, there is an $SO(3)$-action
on $\mathbf{GL}_T$ given by
each ghost vertex. Suppose $v_1,\ldots, v_l$  are ghost vertices
of $T$. Each $v_i$ corresponds to an $SO(3)$ and denoted by
$SO(3)_{v_i}$. Define the isotropy group of $T$
$$
\Gamma_T= \prod_{i=1}^l SO(3)_{v_i}.
$$
The gluing data now is 
$\mathbf{GL}_T/\Gamma_T$. For non-ghost bubble tree,
define $\Gamma_T = 1$.

Suppose edge set of $T$ is $D=\{e_1,\ldots,e_n\}$.
$Gl_T= \prod_i R^4_i/\mathbb{Z}_2$.
We now define a map $\psi: Gl_T\to \mathcal{T}_K$.
Suppose $\underline{x}=(x_1,\ldots, x_n)\in Gl_T$.
We contract $T$ at edges $e_i$ for those $x_i\not= 0$
and set $\psi(\underline{x})$ to be the resultant
tree. For any index set $I\subset \{1,,\ldots,n\}$
let $\underline{x}=(x_1,\ldots,x_n)$ be any point
such that $x_i\not= 0$ iff $i\in I$.
Let $T'=\psi(\underline{x})$. Define
$$
Gl_{T,T'}= (\prod_{i\in I} (R^4_i\setminus\{0\})
/\mathbb{Z}_2)\times\{0\} \subset Gl_T
$$
and $\mathbf{GL}_{T,T'}$ to be the sub-bundle of $GL_T$
with fiber $Gl_{T,T'}$.
If nontrivial  $S_\mathfrak{w}$ is also considered,
$\Gamma_T$ and  $Gl_{T,T'}$
are still well defined.
Then the gluing theory for the general case is
\begin{theorem}[\cite{T}]
Let $(T,v_0)\in \mathcal{T}_K$. For any proper
set $U\subset \mathcal{S}_T(X)$
there is a small constant $\epsilon_0$ and gluing map
$$
\Psi_T: (\mathbf{GL}_T(\epsilon_0)|_U)/\Gamma_T
\to \overline{\mathcal{M}}_{K}(X),
$$
such that
$\Psi_T$ maps $\mathbf{GL}_{T,T'}/\Gamma_T$
to $\mathcal{S}_{T'}$
diffeomorphically for any $T'>T$.
\end{theorem}
Let
$$
\mathcal{D}(X,K)
=\{
(\Psi_T^{-1}((\mathbf{GL}_T(\epsilon_0)|_U)/\Gamma_T),
 \Psi_T^{-1})|
T\in \mathcal{T}_K
\}
$$
where $U,\epsilon_0$ are chosen as  in the theorem. This
defines an {\it atlas}
for $\overline{\mathcal{M}}_K(X)$ and makes it 
into a  topological space. Moreover
\begin{corollary}
$\overline{\mathcal{M}}_K(X)\setminus S_K(X)$
is an orbifold.
\end{corollary}
\noindent
{\bf Proof: } The only thing that needs to be  checked is that
the transition maps between different charts are continuous.
This is proved in
\S 7.2.8 \cite{DK} or  referred to \S 3.3.
q.e.d.
\vskip 0.1in
\noindent
All the constructions and results works for
$\overline{\mathcal{M}}^b_K$
parallelly. 

We now have topology defined on $\overline{\mathcal{M}}_K(X)$
via gluing. On the other hand, we also have 
the Parker-Wolfson bubble tree compactification
theorem. It has been  known that 
the convergence defined in  bubble tree compactification
is compatiable to the given  topology  on $\overline{\mathcal{M}}_K(X)$.
We now  apply the Parker-Wolfson bubble tree compactification
theorem to show the compactness of $\overline{\mathcal{M}}_K(X)$
and $\overline{\mathcal{M}}^b_K$.
\begin{theorem}[\cite{PW}]
1) For any sequence $\{[A_n]\}_{i=1}^\infty
\subset \mathcal{M}_K(X)$
there exists a subsequence such that it
converges to a bubble tree instanton
in $\overline{\mathcal{M}}_K(X)$
via bubble tree compactification process.
2) The bubble tree compactification 
is consistent with the topology defined by $\mathcal{D}(X,K)$.
3) For $X=S^4$, $\overline{\mathcal{M}}^b_K$
is compact. Otherwise, if $\mathbf{F}_K=\emptyset$,
$\overline{\mathcal{M}}_K(X)$ is compact.
\end{theorem}
\noindent
{\bf Proof: } The statement (1) is just the Parker-Wolfson theorem.
The second statement can be proved by the
gluing theorem.
 More precisely, suppose the limit $\mathbf{A}\in \mathcal{S}_T$.
Apply $\Psi_T$ to  $U\times Gl_T(\epsilon)$, where  $U$
is a small neighborhood of $\mathbf{A}$ in $\mathcal{S}_T$.
Since $\Psi_T$ is a diffeomorphism map, one can  show
that $[A_n]$ is in the image of $\Psi_T$ when $n>N$ for some $N$.
We skip the details.  One can find the proof, for example,
in \cite{T}. Now we show (3).
Suppose that $\{\mathbf{A}_n\}_{i=1}^\infty$ is a sequence
in $\overline{M}_K(X)$.
Take a positive sequence $\delta_n\to 0$.
If $\mathbf{A_n}$ is not in $\mathcal{M}_K(X)$,
it is in some chart
$(\Psi_T^{-1}(\mathbf{GL}_T(\epsilon_0)|U), \Psi_T^{-1})$.
Suppose $\mathbf{A}_n= \Psi_T(A'_n, b_n)$, where
$A'_n\in U, b_n\in Gl_T$. Choose any $b'_n \in Gl_T$ such that
$|b_n-b'_n|< \delta_n$ and
  $\tilde{A}_n= \Psi_T(A'_n, b'_n)\in \mathcal{M}_K$.
Now we have a sequence  $\tilde{A}_n$ in
$\mathcal{M}_K(X)$.
Apply the bubble tree compactification to this
new sequence  and get a limit $\mathbf{A}
\in \overline{\mathcal{M}}_K$.  Clearly, this is also the
limit of $\mathbf{A}_n$.
q.e.d.

\subsection{Smoothness of $\overline{\mathcal{M}}_K(X)$}
From 
the previous section, we know that $\overline{\mathcal{M}}_K(X)\setminus S_K$
is a manifold. It is natural to ask if it is smooth.
The problem is to study the gluing maps given by different strata. 
 A straightforward way is to show that the transition
maps between two gluing maps have certain degrees of smoothness.
In \cite{R}, a similar situation is treated and
the author proved that the transition maps are
$C^1$. 
By remark ?, a result as ? is the key.  
In particular, it already provides $C^\infty$-smoothness.
In this subsection, 
we change the strategy
slightly here:
we ``patch'' the gluing maps by perturbing one of the
gluing maps. Hence the transition maps
are relatively simple
 and they are smooth automatically.

We illustrate the difficulty of the problem
for the first  nontrivial case. Let
$T$ be the bubble tree corresponding to bubble
tree manifolds 
$$X \coprod S^4_1\coprod S^4_2/\{p\sim \infty_1, q\sim \infty_2\}$$
where $p\in X, q\in S^4_1$.
$Gl_T = R^4_1/Z_2\times R^4_2/Z_2$
and the coordinate on it is denoted by  $(x_1,x_2)$.
For a proper open set $U\subset \mathcal{S}_T$,
let
$(\Psi_T^{-1}(U\dot{\times} Gl_T(4\epsilon_0)), \Psi_T^{-1})
$ be a chart.
At the mean while, there are two other charts intersecting
this chart. Define the bubble trees $T_i = \psi_T(R^4_i\setminus
\{0\}\times \{0\})$. Set $V_i=B_i(\epsilon_0,
4\epsilon_0)\times \{0\}\subset R^4_i\times \{0\}$
and $U_i=\Psi_{T,T_1}(U\times V_i)$.
$U_i\subset \mathcal{S}_{T_i}$ are proper in their strata.
These  also  provide coordinate charts
$(\Psi_{T_i}^{-1}(U_i
\dot{\times} Gl_{T_i}(\epsilon_i)), \Psi_{T_i}^{-1}).
$
The problem is whether the transition maps on the overlaps
of these three charts are smooth.
For example, consider charts
$$
(\Psi_T^{-1}(U\dot{\times} Gl_T(4\epsilon_0)), \Psi_T^{-1})
\mbox{ and }
(\Psi_{T_1}^{-1}(U_1\dot{\times} Gl_{T_1}(\epsilon_1)),
\Psi_{T_1}^{-1}).
$$
We  prove
\begin{prop}
There is a perturbed gluing map
$\overline{\Psi}_{T_1}$ of
$\Psi_{T_1}$ such that the transition map
between charts
$(\Psi_T^{-1}(U\dot{\times} Gl_T(2\epsilon_0)), \Psi_T^{-1})$
and $(\overline{\Psi}_{T_1}^{-1}(U_1\dot{\times} Gl_T(\epsilon_1)),
\overline{\Psi}_{T_1}^{-1})$ is smooth for small 
$\epsilon_1$, where $\epsilon_1$ depends
on $\epsilon_0$.
\end{prop}
\noindent
{\bf Proof: }
We split the perturbation of $\Psi_{T_1}$ into two steps.
\\
Step 1:
We know that the map  $\Psi_{T,T_1}: U\dot{\times} V_1\to \mathcal{S}_{T_1}$
is diffeomorphic by theorem 3.32.
For any proper $U$ that we consider in the gluing theory,
when $\epsilon_0$,  depending on $U$, is small,
we may assume that $\mathbf{GL}_T$ is a product of two bundles
$\mathbf{GL}_T^i$ whose fiber is $R^4_i/\mathbb{Z}_2$.
Pull back $\mathbf{GL}_T^2$ to $U\dot{\times} V_1$ and denote it by
 $\mathbf{GL}_{T,1}^2$.
$(\Psi_{T,T_1})_* \mathbf{GL}_{T,1}^2$ gives a trivialization  $B_{T,T_1}$ 
for $\mathbf{GL}_{T_1}$ over $\Psi_{T,T_1}(U\dot{\times}V_1)$
by a proper bundle map
$$
B_{T,T_1}: (\Psi_{T,T_1})_* \mathbf{GL}_{T,1}^2
\to \mathbf{GL}_{T_1}.
$$
So, without the loss of  generality, we can  use coordinates of
$\mathbf{GL}_T$ for both $\mathbf{GL}_T$ and $\mathbf{GL}_{T_1}$.
Define
$$\Psi_{T}^1: U\dot{\times} Gl_T(4\epsilon_0)
\to \mathcal{S}_{T_1}\dot{\times} Gl_{T_1} ,
\Psi_T^1(\mathbf{A}, x_1,x_2)
=(\Psi_T(\mathbf{A}, x_1,0),x_2).
$$
We first construct a map
$
\Psi_1
$
that is in the intermediate stage between
$\Psi_T$ and $\Psi_{T_1}\circ\Psi_T^1$.
For simplicity, we assume $U$ consists of only one point
$\mathbf{A}=(A_1,A_2,A_3)$.
Recall that $\Psi_T(\mathbf{A},x_1,x_2)$ is constructed
by:  splice $A_1,A_2,A_3$ together
using gluing parameter $x_1, x_2$ first, then solving for 
instanton. Formally write these two steps as
$$
\Psi_T'(\mathbf{A},x_1,x_2)
= A_1\sharp_{x_1}A_2\sharp_{x_2}A_3
= (A_1\sharp_{x_1}A_2)\sharp_{x_2}A_3
$$
and
$$
\Psi_T(\mathbf{A},x_1,x_2)
=(A_1\sharp_{x_1}A_2)\sharp_{x_2}A_3
+ a( (A_1\sharp_{x_1}A_2)\sharp_{x_2}A_3),
$$
where $A_1\sharp_x A_2$ denotes the splicing of 
$A_1$ and $A_2$ with respect to gluing parameter
$x$. 
Now $\Psi_1$ is defined to be
\begin{eqnarray*}
\Psi_1(\mathbf{A},x_1,x_2)
&=&
(A_1\sharp_{x_1}A_2+a(A_1\sharp_{x_1}A_2))\sharp_{x_2}A_3\\
&+& 
a((A_1\sharp_{x_1}A_2+a(A_1\sharp_{x_1}A_2))\sharp_{x_2}A_3).
\end{eqnarray*}
The real meaning of this
expression is that we
glue $A_1,A_2$ to an instanton
by $x_1$ first and then  glue the resultant
instanton with $A_3$ by $x_2$.
Our first perturbation is
\vskip 0.05in
\noindent
{\it Claim 1: }{$\Psi_1$ can be perturbed
to a new gluing map
$
\overline{\Psi}_1$ defined on $U\dot{\times}Gl_T(4\epsilon_0)
$ such that
$\Psi_T=\overline{\Psi}_1$
over $U\dot{\times} Gl_T(2\epsilon_0)$.
}
\\
{\it Proof: }
Let $\gamma$ be a cut-off function
such that $\gamma(t)=1, t>3\epsilon_0$
and is supported in $t>2\epsilon_0$. Define
\begin{eqnarray*}
\overline{\Psi}_1(\mathbf{A},x_1,x_2)
&=&
[A_1\sharp_{x_1}A_2+
\gamma(|x_1|)a(A_1\sharp_{x_1}A_2)]\sharp_{x_2}A_3\\
&+&
a[(A_1\sharp_{x_1}A_2+\gamma(|x_1|)
a(A_1\sharp_{x_1}A_2))\sharp_{x_2}A_3].
\end{eqnarray*}
Obviously, $\overline{\Psi}_1$ satisfies our goal.
We need to show that it is also a gluing map.
As we did before,  we compare the map with $\Psi'_T$
and study the derivative
of $\overline{\Psi}_1-\Psi'_T$.
As usual, we can show that
$$
\|d\Psi'_T-d\overline{\Psi}_1\|
\leq C(|x_1|+|x_2|)^{3/q}.
$$
and
$$
\|\overline{\Psi}_1(\mathbf{A},x_1,x_2)
-\Psi'_T(\mathbf{A},x_1,x_2)\|_{L^{1,q}}
\leq C(|x_1|+|x_2|)^{1+3/q}.
$$
This is sufficient to imply that $\overline{\Psi}$
is a  gluing map. 
q.e.d.
\vskip 0.05in
\noindent
Let $\Psi'_{T_1}=
\overline{\Psi}_1(\Psi^1_T)^{-1}$.
We have shown that the transition map
between 
charts $(\Psi_T(U\dot{\times} Gl_T(2\epsilon_0)),
 \Psi_T^{-1})$
and
$(\overline{\Psi}_1(U\dot{\times} V_1\times
B_2(4\epsilon_0))),
(\Psi'_{T_1})^{-1})$ is $\Psi_T^1$.
\\
{\it Step 2: }
One might think that
$\Psi'_{T_1}$ (or $\overline{\Psi}_1$)
is $\Psi_{T_1}$ intuitively, however they are different
because  the metrics used
for two gluing maps  are different. This problem is dealt with
in
\vskip 0.1in
\noindent
{\it Claim 2:
there exists a small constant $\epsilon_1>0$ and
a gluing map $\overline{\Psi}_{T_1}$
such that
$\overline{\Psi}_{T_1}=\Psi'_{T_1}$ over
$\Psi_T^1(U\dot{\times} 
(B_1(2\epsilon_0,3\epsilon_0)
\times B_2(\epsilon_1)
 )$ and
$\overline{\Psi}_{T_1}=\Psi_{T_1}$ over
$\Psi_T^1(U\dot{\times} (B_1(\geq 4\epsilon_0)
\times B_2(\epsilon_1)
)$.
}
\\
{\it Proof: }
The difference between these two maps
over $\Psi_T^1(U\dot{\times} B_1(3\epsilon_0,4\epsilon_0)
\times B_2(\epsilon_1)
 )$ can be seen by what follows.
We make a convention that the metric we refer to when we do
comparing is the one induced by
$X\sharp_{x_1}S^4$. 
Let $\mathbf{B}=\Psi_{T,T_1}(A_1,A_2, x_1,0)$. Then
\begin{eqnarray*}
\Psi'_{T_1}(\Psi_T^1(A_1,A_2,A_3,x_1,x_2))
&=& \mathbf{B}\sharp_{x_2}A_3 + a'(\mathbf{B}\sharp_{x_2}A_3) \\
\Psi_{T_1}(\Psi_T^1(A_1,A_2,A_3,x_1,x_2)
&=& \mathbf{B}\sharp_{x_2}A_3 + a(\mathbf{B}\sharp_{x_2}A_3).
\end{eqnarray*}
Here, 
$a'$ and $a$, given
in proposition 3.17 by using different metrics on $X$, are different.
$\Psi_{T_1}$ uses the standard metric while
$\Psi'_{T_1}$ uses the metric over $X\sharp_{x_1}S^4$.
We know that
$$
\|da'\|\leq Cr_2^{3/q}\|d\Psi'_{T_1}\|.
$$
We also have the  same estimate for $\Psi_{T_1}$
using the standard metric on $X$. Translate these estimates
to the metric we are using now, then  
$$
\|da\|\leq C(\epsilon_0)r_2^{3/q}\|d\Psi'_{T_1}\|.
$$
To glue two maps together, one can do either of the following ways:
\\
1, deform the metric on $X$ with respect to gluing parameter $x$
such that the metric is standard one when $|x|>3\epsilon_0$;\\
2, deform the right inverse $P$'s. To be precise, 
let $P,P'$ be the right inverse used to define  $a,a'$. 
Note that we have a
simple fact:$\lambda P+ (1-\lambda)P'$ is still a right inverse.
With this, we can easily deform $P$'s by cut-off functions.\\
If we choose $\epsilon_1$ small enough, then $r_2$ is small.
By either way, 
 $\Psi'_{T_1}$ can be deformed so that it matchs
 $\Psi_{T_1}$ when $|x_1|> 3\epsilon_0$. This glues 
$\Psi'_{T_1}$ and $\Psi_{T_1}$ together.
q.e.d.
\vskip 0.1in
\noindent
We have  finished the proof of proposition 
3.35.
q.e.d.
\begin{remark}
We make a remark on a compatibility property
of gluing maps $\Psi_T$ and $\overline{\Psi}_{T_1}$.
$$
\begin{array}{cccccc}
(\Psi_{T,T_1})_\ast\mathbf{GL}_T^2
&
\xrightarrow{B_{T,T_1}}
&         
\mathbf{GL}_{T_1}
& 
(\Psi_{T,T_1})_\ast\mathbf{GL}_T^2
&
\xrightarrow{B_{T,T_1}}
&         
\mathbf{GL}_{T_1}
\\
\Big\downarrow\vcenter{\rlap{$\Psi_T$}} 
& 
&
\Big\downarrow\vcenter{\rlap{$\Psi_{T_1}$}}
& 
\Big\downarrow\vcenter{\rlap{$\Psi_T$}} 
& 
&
\Big\downarrow\vcenter{\rlap{$\overline{\Psi}_{T_1}$}}
\\
V_T
&
\xrightarrow{?}
&
V_{T_1}
& 
V_T
&
\xrightarrow{\overline{\Psi}_{T_1}B_{T,T_1}\Psi_T^{-1}}
&
V_{T_1}
\end{array}
$$
In the first diagram, the map ``?'' does not
agree with  the transition map provided by the gluing theory,
while in the second one, we know the transition map is
the composition of three smooth maps. 
Therefore, the transition map is smooth. 
So, roughly speaking, the patching is  
the patching of gluing bundles by $B_{T,T_1}$
up to 
$\Psi_{T,T_1}$.
\end{remark}
We  can prove the main theorem in this section
\begin{theorem}
Let $K>0$ be any integer. There is an atlas
$\mathcal{D}(X,K)$
for $\overline{\mathcal{M}}_K(X)$
such that the transition maps
are $C^\infty$. So $\overline{\mathcal{M}}_K(X)
\setminus \mathbf{S}_K(X)$
is a smooth orbifold.
\end{theorem}
\noindent
{\bf Proof: }
The construction is based on the previous proposition.
Recall that we defined a partial order on $\mathcal{T}_K$.
We introduce the charts for strata with respect to this order
inductively and require
that the atlas satisfies  several requirements ({\bf R1, R2, R3})
\\
{\bf R1: } To each  stratum $\mathcal{S}_T$
only one chart $(\overline{\Psi}_T(U_T
\dot{\times} Gl_T(\epsilon_T)), \overline{\Psi}_T^{-1})$ is assigned,
where $U_T$ is a proper open subset of $\mathcal{S}_T$.
To simplify notations, we denote  charts
as  $chart(V_T,\overline{\Psi}_T)$, where
$V_T= \overline{\Psi}_T(U_T
\dot{\times} Gl_T(\epsilon_T))$.
\\
{\bf R2: } For any two charts $chart(V_T, \overline{\Psi}_T),
chart(V_{T'}, \overline{\Psi}_{T'})$, $V_T\cap V_T'\not=
\emptyset$  iff either $T<T'$ or $T'<T$.
\\
We start with lowest strata which are compact.
Choose a $\epsilon_0>0$ such that
for each $T$, $(\Psi_T(\mathcal{S}_T\dot{\times}
Gl_T(\epsilon_0), \Psi_T^{-1})$ gives a coordinate chart
for the neighborhood of $\mathcal{S}_T$. Set
$\overline{\Psi}_T=\Psi_T$.
Moreover, for small $\epsilon_0$, {\bf R2} can be satisfied.
Let $\mathcal{U}\subset\overline{\mathcal{M}}_K(X)$
be the union of all open sets in given charts.
This set is updated as charts are created.
For any stratum $T'$ whose chart has  not been constructed,
gluing maps of given charts induce a gluing map
over $\mathbf{GL}_{T'}|_{(U\cap \mathcal{S}_{T'})}$
(we simply say ``over $U\cap\mathcal{S}_{T'}$'' later).
The statement is
\\
{\bf R3: } There is a gluing map  $\Psi^0_T$
over $\mathcal{U}\cap \mathcal{S}_T$ induced
by defined gluing maps of lower strata.
\\
This is defined as following:
for any point $(\mathbf{A}, x)$, there exists a
lower strata $T<T'$ such that $\mathbf{A}
=\overline{\Psi}_{T,T'}(\mathbf{B}, x')$. Then define
the induced gluing map
$$
\Psi^0_{T'}(\mathbf{A},x) = \overline{\Psi}_{T}(\mathbf{B},x',x).
$$
Note that
 may not be well defined since
there may be more than one choice for $T$.
To solve this, it is equivalent to require
\\
{\bf R3': }
For any two charts $chart(V_{T_1},\overline{\Psi}_{T_1})$ and
$chart(V_{T_2},\overline{\Psi}_{T_2})$,
$\overline{\Psi}_{T_1}$ agrees with $\overline{\Psi}_{T_2}$
on the overlap of their charts.
\\
There is no such a  problem
for the base case.
Now suppose we are ready
to give a chart for  some $T\in \mathcal{T}_K$.
This means that
all strata lower than $\mathcal{S}_T$ already have
charts. So
$\mathcal{S}_T\setminus \mathcal{U}$ is a proper set
in $\mathcal{S}_T$. By induction hypothesis {\bf R3},
we already have a gluing map $\Psi_T^0$
over $M:=\mathcal{U}\cap
\mathcal{S}_T$
induced from given charts. Now  shrink $M$
to $M_1$ and then to $M_2$.
Using $\Psi_T$ over $\mathcal{S}_T\setminus
M_2$ and arguments for proposition 3.35,
 we glue $\Psi_T$ and $\Psi_T^0$ together
and define a new gluing map
$\overline{\Psi}_T$ over  $\mathbf{GL}_T(\epsilon_T)
\setminus M_2$ such that it agrees with $\Psi_T^0$
over $M_1$.
Here $\epsilon_T$ is a small constant
depending lower charts. Now modify the lower charts
such that $\mathcal{U}\cap \mathcal{S}_T$
is $M_1$
other than $M$. (This $\mathcal{U}$
does not count current chart.)
This will guarantee {\bf R3'}
and so  {\bf R3}.
Of course, we can make
$\epsilon_T$ small enough to guarantee {\bf R2}.
{\bf R1} is automatically true.
The transition maps are smooth as explained in proposition 3.35.
The new atlas is denoted by $\mathcal{D}(X,K)$. 
q.e.d.
\vskip 0.1in
$\mathbf{R3}$ is the generalization of remark 3.36.
 We call it the {\it compatibility property}
of the atlas $\mathcal{D}(X,K)$. 
Because of this property,
we can identify open set in every chart as
a gluing bundle for corresponding stratum and they are patched together
by bundle maps $B_{T,T'}$. 
This  is very useful 
to  our flip resolutions in \S 4.3. In fact, this was one of the main
motivations of this section.

\section{Flip Resolution and Smooth Compactification}
In the  previous chapter, we constructed the compactified space 
$\overline{\mathcal{M}}$ via bubble tree compactification.
Although certain smoothness and orbifold structures 
are achieved, there are still singularities wherever 
ghost strata occur. Suppose $T\in G_K$ is a ghost tree with $g(T)$
ghost vertices. The non-finite stabilizer is the product of $g(T)$
copies of $SO(3)$.  Due to these obstructions, 
the advantage we can get from the  bubble tree compactification,
as compared to the Uhlenbeck compactification,
is limited.  Our task in this paper is to
construct an orbifold compactification. 
The ingredient to solve this problem is the following key observation:
$M_T$ is obtained from some manifold by
$g(T)$ ``blow-up''s.  This structure is very similar to the
Fulton-MacPherson (briefly indicated by FM)  compactified space.
Using these ``blowing-up''s to absorb the stabilizer
$SO(3)$'s, we can apply "blowing-down"s
and hence resolve singularities at $\mathcal{S}_T$. We call
this procedure the {\it flip resolution}.

\subsection{Flip Resolution and Smoothing
$\overline{\mathcal{M}}_{K}(X)$ for $K=2$}

We start with the simplest case  $K=2$. 
For $K=2$ we only need one flip resolution, so
it is easy to see how the flip resolution
works. For general $K$, we need to apply flip resolutions
sequentially. 
We will describe this later.

Suppose $K=2$.
The moduli space $\overline{\mathcal{M}}_K(X)$
contains only one ghost stratum $\mathcal{S}_T$ where
$T$ is described as follows:
the vertex set $V$ of $T$ is $V= \{v_0,v_1,v_2,v_3\}$.
$child(v_0)= \{v_1\}, child(v_1)=\{v_2,v_3\}$ 
and $v_1$ is a ghost vertex. $w(v_2)=w(v_3)=1$.
$R(T)=[K-2[0[1,1]]]$. 

Let $STX$ be the sphere bundle of $TX$ with the standard
$\mathbb{Z}_2$ action fiber-wisely and 
$$
\mathcal{Z}_T= \mathcal{M}_{K-2}(X)\times X.
$$ 
Then
\begin{lemma}
$\mathcal{S}_T=\mathcal{M}_{K-2}(X)\times STX/\mathbb{Z}_2$
is a bundle over $\mathcal{Z}_T$.
The gluing bundle  $\mathbf{GL}_T$ is the bundle with fiber
$$Gl_T = R^4_0/\mathbb{Z}_2\times R^4_1/\mathbb{Z}_2\times R^4_2/\mathbb
{Z}_2.$$ 
$\Gamma_T= SO(3)$. The gluing data is $\mathbf{GL}_T/\Gamma_T$.
\end{lemma}
The lemma follows from example 3.4 and
the description at the end of \S 3.2. We skip the 
proof. By the construction of $\mathbf{GL}_T$, 
we know that $\mathbf{GL}_T$ is a pull-back bundle, say
$\mathbf{GL}'_T$, 
over $\mathcal{Z}_T$. In  other words, 
$$
\mathbf{GL}_T = \mathcal{S}_T\times_{\mathcal{Z}_T} \mathbf{GL}'_T.
$$
The neighborhood of the stratum is $\mathbf{GL}_T/SO(3)$.
In particular, $SO(3)$ only acts on $\mathbf{GL}'_T$.
It is worth to take a look at   other strata around 
$\mathcal{S}_T$ and  see how they change after flip
resolutions:
let $(x_0,x_1,x_2)$ denote
 the point in $Gl_T$, set 
$T_0= \Psi_T(1,0,0)$ and $T_1=\Psi_T(0,1,1)$. $\mathcal{S}_{T_0}$
is the stratum that contains two bubble points on $X$.
$\mathcal{S}_{T_1}$  is the stratum that contains only one
bubble point on $X$ with energy $2$.
For simplicity, we assume that $\mathcal{M}_{K-2}(X)$ is a point.
Then the compactification of $\mathcal{S}_{T_0}$
is the real blow-up of $X^2$ along diagonal modulo $Z_2$ action.
The compactification
of $\mathcal{S}_{T_1}$ is still singular at $\mathcal{S}_T$. 
They intersect at $\mathcal{S}_T$. 
\begin{definition}
Suppose $N$ is a compact space
and $U\to N$
is a disk bundle. Let $S\to N$
be the sphere bundle of $U$. The space
$\hat U =[0,1]\times S$
is called the \emph{semi-blow-up} of $U$
along $N$.
We say that
$\{0\}\times S$ is the
semi-blow-up boundary of
$U$.
Conversely, if  $\hat U = [0,1]\times S$
where $S$ is  a sphere bundle of some disk
bundle $U$ over $N$,
$U$ is call the semi-blow-down of $\hat U$.
\end{definition}
This is a rather trivial definition, however
it is the case corresponding to our general model.
Now let us see how we apply semi-blow-up and semi-blow-down
to $\mathbf{GL}_T/SO(3)$. To simplify the notations,
we ignore the orbifold structures and replace $R^4/\mathbb{Z}_2$
\& 
$SO(3)$ by $R^4$ \& $SU(2)$. In general, 
it is easy to see that the 
construction  still works for orbifold cases.
The resolution consists of two steps: 
\\
Step 1: semi-blow up $\mathbf{GL}_T$ along $\mathcal{S}_T$.
We call the resultant space $\mathbf{GL}_T^\sharp$. 
This space is described as follows: 
semi-blow up $\mathbf{GL}'_T$ along 
$\mathcal{Z}_T$. The resultant bundle is 
$$
(\mathbf{GL}'_T)^\sharp
=\mathcal{Z}_T\dot{\times} ([0,1)\times S^7)
=:[0,1)\times \mathbf{SGL}'_T.
$$
So
$$
\mathbf{GL}_T^\sharp/SU(2)
= [0,1)\times\mathcal{S}_T\times_{\mathcal{Z}_T}(\mathbf{SGL}_T'/SU(2)).
$$
Step 2: semi-blow down $[0,1)\times \mathcal{S}_T$ along 
$\mathcal{Z}_T$. Since $\mathcal{S}_T$ is a sphere bundle
over $\mathcal{Z}_T$ of some vector bundle, say $\mathcal{S}_{T,d}$,
we can apply the semi-blow down 
and replace $[0,1)\times \mathcal{S}_T$
by $\mathcal{S}_{T,d}$. Eventually, $\mathbf{GL}_T$
is replaced by 
$$
\tilde{\mathbf{GL}}_T = \mathcal{S}_{T,d}\times _{\mathcal{Z}_T}
 (\mathbf{SGL}_T/SU(2)).
$$
Clearly, none of points has  non-finite stablizer
in $\tilde{\mathbf{GL}}_T$. The procedure here is similar to the
flip in algebraic geometry. So we call such an  operation {\it flip
resolution}. The following is the formal definition.
\begin{definition}
Let $G$ be a group acting on $S^{n-1}$
freely and it induces an action on $R^n$.
Let $SN, SV$ be  sphere bundles of  the vector bundles $N,V$
over $X$, where $V$ is $n$-bundle.
  There is a natural diffeomorphism
$$
f: SN  \times_X (V\setminus X)/G
\to  N\setminus X \times_X SV/G
$$
by
$$
f(n, [v]) = (|v|n, [\frac{v}{|v|}]).
$$
Replacing $SN\times_X V/G$ by
$N\times_X SV/G$ via $f$ is called
\emph{ flip resolution.} Treat $X$ as 0-section
of $N$, we call
 $X\times_X SV/G = SV/G$  the \emph{ exceptional
divisor} of $N\times_X SV/G$.
\end{definition}
We denote the new space by $\underline{\mathcal{M}}_K(X)$.
So we have
\begin{lemma}
$\underline{\mathcal{M}}_K(X)$ is a smooth orbifold.
\end{lemma}
To have a better understanding of this new space, we explain what  the
completions of $
\mathcal{S}_
{T_0}$ and $\mathcal{S}_{T_1}$  are in the new space. The exceptional
divisor takes the place of
the  ghost stratum
$\mathcal{S}_T$. It is 
 $\mathcal{M}_{K-2}(X)\times X\dot{\times} (S^{11}/SU(2))$.
One can check that the completion 
of $\mathcal{S}_{T_0}$ intersects the exceptional
divisor at $\mathcal{M}_{K-2}(X)\times X\times \{[1,0,0]\}=
\mathcal{M}_{K-2}(X)\times X$. 
So
\begin{lemma}
$\overline
{\mathcal{S}}_{T_0}$ in $\underline{\mathcal{M}}_K(X)$
is $\mathcal{M}_{K-2}(X)\times Sym^2(X)$.
\end{lemma}
One can also prove that the completion of $\mathcal{S}_{T_1}$
is a smooth sub-orbifold in $\underline{\mathcal{M}}_K(X)$
and it is 
$$
\mathcal{M}_{K-2}(X)\times Fr(X)_{SO(4)}
\times \underline{\mathcal{M}}^b_2.
$$
where
$\underline{\mathcal{M}}^b_2$ is obtained from
$\overline{\mathcal{M}}^b_2$ by flip resolution in a similar way.
Note that original $\overline{\mathcal{M}}^b_2$
is also singular.

\subsection{FM-Compactification and Ghost Strata}
We just explained that how to apply
the flip resolution to the moduli space when $K=2$.
As $K$  gets bigger, the ghost strata are very complicated
and interwoven. Fortunately, they  basically 
follow the rules of the FM-compactification of
configuration space ${F(X,n)}$.
\vskip 0.1in
\noindent
{\it  Review of  the FM-compactification.}
\vskip 0.05in
\noindent
Suppose $X$ is a smooth, compact complex manifold.
Let $X^{[n]}:= X^n/S_n$, where $S_n$
is the permutation group. This space is also denoted by
$Sym^n(X)$.
(In 
\cite{FuM}, the construction is given for 
$X^n$.) 
Define
$$
F(X,n) = X^{[n]}\setminus \Delta,
$$
where $\Delta$ is the big diagonal in $X^{[n]}$.
 In general, for a tuple
$I= (i_1,\ldots, i_k), 1\leq i_1<\cdots<i_k=n$
define
$$
\Delta_I = \{ (x_1,\ldots,x_n)\in X^n|
x_{i_1}=\cdots = x_{i_k}\}.
$$
Let $s: X^n\to X^{[n]}$ be the projection. $\Delta_I$
also stands for $s(\Delta_I)$, for simplicity. So
$\Delta = \cup \Delta_{(a,b)}$. In \cite{FuM} a compactified
space $\overline{F(X,n)}$ is given.
The construction can be described as follows:
\\
1) (complex)-blow up $X^n$ along diagonal
$\Delta_{(1,2,\ldots, n)}$ and denote the new
manifold by $X^n_n$. Accordingly,  all other $\Delta_I$
are changed after blow-up. 
We still use  $\Delta_I$ for the updated diagonals.
Note that diagonals $\Delta_I, |I|=n-1,$
are disjoint in $X^n_n$.
\\
2) Inductively, suppose $X^n_k$ is created.
If $k=2$, we are done. Else, one can also
prove that diagonals $\Delta_I, |I|=k-1,$
are disjoint in $X^n_k$ by induction.
Blow up
$X^n_k$ along all these diagonals and create
$X^n_{k-1}$. As a result, $\Delta_I, |I|=k-2,$
are disjoint.
\\
3) The action of $S_n$ on $X^n$ has a natural extension
to $X^n_k$. Let
$$
\emph{}\overline{F(X,n)} = X^n_2/S_n.
$$
There is a projection map
$FM: \overline{F(X,n)}\to X^{[n]}$.
It is clear that reversing the blow-up process and blowing-down
$\overline{F(X,n)}$ we retrieve $X^{[n]}$.
Intuitively,
we
introduce blow-ups  wherever
two or more points run into one another.

When $X$ is a real manifold,  one can  also
construct $\overline{F(X,n)}$ using  real blow-up.
Here, we use ``semi-blow-up'' to
form a ``stratified'' space $\overline{F(X,n)}$.
Each stratum
is isomorphic to some $M_T$. 
We  explain 
the construction of $\overline{F(X,n)}$
for
$n\leq 3$. It 
 can be generized to all $n$.
\begin{example}[$n=2$]
The space $X^2_2$ is the semi-blow-up of $X^2$
along diagonal. And
$\overline{F(X,n)} = X^2_2/S_2$. Let $(X^2)_\Delta$
 be the real blow-up of $X^2$
along $\Delta$. Then
$$
\overline{F(X,n)}= (X^2)_\Delta/S_2.
$$
So when $n=2$, this coincides with FM-compactification
via real blow-ups.
$X^2_2$ has two strata $\Omega_1$ and
$\Omega_2$: $\Omega_1$ is the boundary $B$ of $\overline{F(X,2)}$ and
$\Omega_2=F(X,2)$.
\end{example}
\begin{example}[$n=3$]
$X^3_3$ is the semi-blow-up of $X$ along the diagonal 
$\Delta_{(1,2,3)}$. It is 
a manifold with boundary $B_{123}$.
$b_{123}:B_{123}\to X=\Delta_{(1,2,3)}$
is a sphere bundle over $X$.
We can describe points on $b_{123}^{-1}(x)$ as
follows: suppose $z\in b^{-1}_{123}(x)$, $z$ can be
written as $z=(x,[x_1,x_2,x_3]), x_i\in TX_x$ such that
\begin{equation}
x_1+x_2+x_3 =0, x_1^2+x_2^2+x_3^2 = 3.
\end{equation}
In fact, $z$ is treated as the limit
of triples $(x+tx_1, x+tx_2,x+tx_3)$ in $X^3\setminus
\Delta$.  $[x_1,x_2,x_3]$ is the equivalent class
$(x_1,x_2,x_3)/\sim$, where $(x_1,x_2,x_3)\sim
(x'_1, x_2',x_3')$ when $x_i=\lambda x_i, \lambda>0$ ([FM]).
The coordinates of neighborhood of $B_1$
can be identified as
$
(x, [x_1,x_2,x_3], t).
$
When $t\not= 0$, $(x,[x_1,x_2,x_3],t) = (x+tx_1,x+tx_2,x+tx_3)\in
X^3$. Usually we drop $t$ if it is equal to 0.
Now we blow up $X^3_3$ along
$\Delta_I, I=(1,2), (1,3), (2,3)$ and get $X^3_2$.
Consider $\Delta_{(1,2)}$. There is a natural
normal bundle for $\Delta_{(1,2)}\setminus B_{123}$
in $X^3_3\setminus B_{123}$: The fiber over
point $(x,x,y)$ is identified with normal directions given by
$(x+tx',x-tx',y), 
 x\in TX_x$.
This bundle can be naturally extended over 
$\Delta_{(1,2)}\cap B_{123}$.
We describe this
in the neighborhood of $B_{123}$. Given a point
$(x,[x_1,x_1,x_3],t)=(x, x+tx_1,x+tx_1,x+tx_3)$ with $t>0$,
a point $y\in TX_{x+tx_1}$  in the normal direction  is 
of the form
$$
\frac{d}{ds}|_{s=0} (x+tx_1+tsy, x+tx_2-tsy, tx_3).
$$
The path can be given in terms of our new coordinate
system. One can show that
it can be extended over $t=0$. This shows the bundle converges
to
the normal vector of $\Delta_{(1,2)}\cap B_{123}$ in
$B_{123}$.
Identify the neighborhood of $\Delta_{(1,2)}$  as
this normal bundle and we can apply the semi-blow-up
along $\Delta_{(1,2)}$.
The new boundary is denoted by $B_{12}$. This is a stratified space.
We denote the coordinates for the smallest 
stratum $B_{12}\cap B_{123}$ in
 the following way:
$$
(x[x_1 [y, -y], x_3]).
$$
Note that $x$ is followed by a bracket $[x_1[...], x_3]$ and 
$x_1$ is followed by $[y, -y]$
without comma. We explain the format of the notations: 
$x_1,x_3\in T_xX$, $y, -y \in T_{x_1}(T_xX)$.
The geometric meaning is that there are two points
(since the bracket after $x$ contains two points) that  run into 
$x$ and cause 
 a blow-up at $x$ and then 
there are two points run into $x_1$
and cause the second blow-up at $x_1$.
Moreover, the coordinates  satisfy
$$
2x_1 +x_3 =0; 2x_1^2 +x_3^2=1; y^2+(-y)^2=1.
$$
Correspondingly, blow up along $\Delta_{(2,3)}
$, $\Delta_{(1,3)}$ and  create the other
two boundaries
and denote them by $B_{23}, B_{12}$.
The resultant space is $X^3_2$ and $\overline{F(X,3)}=X^3_2/S_3$.
Totally, there are
four strata $\Omega_i, 1\leq i\leq 4$ in the space.
We describe the strata and formulate
the representations of points in the strata as follows.
\begin{eqnarray*}
\Omega_1 &=& \{(B_{12}\cup B_{23}\cup B_{13})\cap B_{123}\}/S_3,
[x, [y [z,-z], w]]; \\
\Omega_2 &=& \{B_{123}\setminus
(B_{12}\cup B_{23}\cup B_{13})\}/S_3, [x [y,z,w]] \\
\Omega_3 &=& \{(B_{12}\cup B_{23}\cup B_{13})\setminus B_{123}\}
/S_3,
[x [y, -y], z]; \\
\Omega_4 &=& F(X,3), [x,y,z]
\end{eqnarray*}
We explain the notations of points:
for $[x,y,z]$,  $x,y,z\in X^3\setminus
\Delta$ and $[x,y,z]=(x,y,z)/S_3$;
for $[x[y,-y],z]$,  $x,z\in X^2\setminus \Delta$.
$y,-y\in TX_x$, this corresponds to a blow-up
at $x$ and described by $x[y,-y]$;
for $[x[y,z,w]]$ this means there is a blow-up
at $x$ and $y,z,w\in (TX_x)^3\setminus \Delta$;
Similarly $[x[y[z ,-z] , w]]$ says there is a blow-up
at $x$ and  $y,w\in (TX_x)^2\setminus \Delta$,
and then followed by a blow-up
at $y$ and $z,-z\in T(TX_x)_y$.
\end{example}
In general,
a point in $\overline{F(X,n)}$
is denoted by
$$
p= [x_0[\ldots[\ldots] ], x_1[\ldots],\ldots].
$$
Any point $y$ included in the bracket following
a point $x$ is called {\it descendant} of $x$.
Two points in the same level, like $x_0,x_1$ in the example,
are called  {\it siblings}. By the construction, siblings
should  differ from each other.
It is clear that the representation of
points in same stratum
has same format. So we also say it is the format of
the stratum.
There is a natural way to construct bubble trees
in terms of the formats of strata.
\begin{example}
 We continue example 4.7. We intend
to build bubble trees $T_i$ such that
$M_{T_i}= \Omega_i$. Suppose $T_i\in \mathcal{T}_K$.
In the construction, we do not
care what the weight
 of the root $v_0$ is. $v_0$ is always assigned with $X$.
 Weights for other nontrivial
vertices
are all equal to 1 in the construction. Namely
\begin{eqnarray*}
R(T_1)&=& [K-3[0[1,0[1,1]]]] \\
R(T_2)&=& [K-3[0[1,1,1]]] \\
R(T_3)&=& [K-3[1,0[1,1]]] \\
R(T_4)&=& [K-3[1,1,1]]
\end{eqnarray*}
and representations of bubble tree mainfolds of $T_i$
are
\begin{eqnarray*}
Y(d_{T_1})&=&[x_1[x_2,x_3[x_4,x_5]]] \\
Y(d_{T_2})&=&[x_1[x_2,x_3,x_4]] \\
Y(d_{T_3})&=&[x_1,x_2[x_3,x_4]] \\
Y(d_{T_4})&=&[x_1,x_2,x_3]
\end{eqnarray*}
\end{example}
Following these examples, one can generalize the correspondences
for $n>3$.
A general condition similar to (18) is  the
{\it balanced condition} (3) 
 on ghost bubble: suppose
$\{x_1,\ldots, x_k\}$ are points on ghost bubble $S^4$,
then
$$
\sum_{i=1}^k x_i=0, \mbox{ and } \sum_{i=1}^k |x_i|^2=k.
$$

We now generalize the construction to
the weighted cases.
Let $\mathfrak{w}= (w_1,\ldots, w_n)\in Z^n, w_i>0$ and
$S_\mathfrak{w}$ be the associated permutation group.
Define $X^n_\mathfrak{w}$ to be the {\it weighted space
} of $X^n$. It is the collection of
points $(x_1,\ldots, x_n)\in X^n$  with
assigned weights $w_i$ to $x_i$.
Define $F_\mathfrak{w}(X,n) = X^n_\mathfrak{w}
\setminus \Delta$.
Similarly, we can define
$$
\overline{F_\mathfrak{w}(X,n)} = X^n_{\mathfrak{w},2}/S_\mathfrak{w}.
$$
\begin{example}
$\overline{F_\mathfrak{w}(X,3)}$.\\
$X^3_{\mathfrak{w},3}$: the key point is to
choose a normal bundle and its sphere bundle.
Define a sub-bundle $N_\mathfrak{w}\subset TX\times_X TX\times_X TX\to
\Delta_{(1,2,3)}=X$. For any $x\in X$ the fiber consists of
$$
\{(x_1,x_2,x_3)\in (TX_x)^3|
\sum_{i=1}^{3} w_ix_i=0\}
$$
Take sphere bundle $SN_\mathfrak{w}$ of $N_\mathfrak{w}$ by imposing
$$
w_1|x_1|^2+ w_2|x_2|^2 +w_3|x_3|^2=w_1+w_2+w_3
$$
with $\sum w_ix_i=0$. This is the { balanced condition}.
In general 
 if there are $\{x_1,\ldots,x_k\}$ points
with assigned weights $w=(w_1,\ldots, w_k)$, then
$$
\sum_{i=1}^k w_ix_i=0 \mbox{ and } \sum_{i=1}^k w_i|x_i|^2=\sum w_i
$$
We call this  the {\it balanced condition of weight $w$}.
These definitions are made to be consistent with (3). 
The neighborhood of $\Delta_{(1,2,3)}$ in $X^3\setminus \Delta
_{(1,2,3)}$
can be identified with $N_\mathfrak{w}$ and so with
$
(0,1)\times SN_\mathfrak{w}.
$ The weighted blow-up
is to replace the neighborhood by $[0,1)\times S_\mathfrak{w}$.
One can mimic  example 4.7 and do (weighted)-blow-up
of $X^3_3$ along $\Delta_{(1,2)}$ and so on.
The choice of bundle and its sphere bundle is made to
be consistent with balanced conditions.
\end{example}
\vskip 0.1in
Weighted spaces are natural objects for our model
and they are the generalization of FM-models.
Given $\mathfrak{w}= (w_1,\ldots,w_n)$,
a tree $T$ is {\it $\mathfrak{w}$-type}
if it has only $n$ nontrivial
vertices, say $v_1,\ldots, v_n$, with weights
$w(v_i)=w_i$ and all vertices are ghosts except
leaves.
Let $G_\mathfrak{w}$ be the collection of all 
$\mathfrak{w}$-type bubble
trees. One can show that
\begin{lemma}
$\overline{F_\mathfrak{w}(X,n)}$ is a stratified space. Moreover
$$
\overline{F_\mathfrak{w}(X,n)}= \cup_{T\in G_\mathfrak{w}} M_T.
$$
\end{lemma}

Now we know that $\overline{F_\mathfrak{w}(X,n)}$
is obtained by a sequence of blow-ups. And  example 4.7
also shows how to identify strata of 
$\overline{F_\mathfrak{w}(X,n)}$
with $M_T$. If we are just looking for an orbifold
compactification of $F(X,n)$, obviously our
FM-space
is a not a wise choice. The easiest one,  of course, is $X^{[n]}$ itself.
Note that to go from $\overline{F_\mathfrak{w}(X,n)}$ to
$X^{[n]}$, we can  reverse the blow-up process and do blow-downs.
$\overline{F_\mathfrak{w}(X,n)}$ is the right model to describe
$\overline{\mathcal{M}}_{K}(X)$. The redundancy of 
$\overline{F_\mathfrak{w}(X,n)}$ comparing to $X^{[n]}$ 
 suggests that we might 
simplify $\overline{\mathcal{M}}_{K}(X)$ and get a better 
compactification.

\vskip 0.1in
\noindent
{\it Moduli spaces with marked points and Ghost Strata}
\vskip 0.1in
\noindent
The concept of marked points is very common  in
the theory of $J$-holomorphic curves. It is natural to
introduce the same concepts to our case without  additional difficulties.
This is not essential in our compactification, however, they make the 
explanation simpler.
\begin{definition}
Let $\mathfrak{w}=(w_1,\ldots, w_n)$.
A marked point $\mathfrak{w}$-instanton of $X$
is an instanton $A$ with an unordered
 tuple $(x_1,\ldots, x_n)$,
$x_i\in X$ and
$x_i\not= x_j, i\not= j$.
Moreover $x_i$ are associated with 
charges $w_i$.
Denote such an element by $[A, (x_1,\ldots,x_n)]$.
\end{definition}
For $\mathfrak{w}=(w_1,\ldots,w_n)$
define
\begin{eqnarray*}
\mathcal{M}_{K, \mathfrak{w}}&=&
\{[A, (x_1,\ldots,x_k)] \mbox{ is a marked point }\\
&& \mathfrak{w}-\mbox{instanton}|  
 A\in \mathcal{M}_{K_0}(X)\}
\end{eqnarray*}
where $K_0= K-\sum w_i$.
Basically, the bubble tree compactification and gluing theory can 
also be generalized
to this case. We denote the compactified space
by $\overline{\mathcal{M}}_{K,\mathfrak{w}}(X)$.
We outline how these generalizations can be done:
\\
{\it Bubble tree compactification of $\mathcal{M}_{K,\mathfrak{w}}$.}
Suppose there is a sequence of marked point
instantons $\{[A_n,( x_{n1},\ldots, x_{nk}]\}_{n=1}^\infty$.
Either bubble points of $\{A_n\}$ or points $x$
satisfying
$$
\lim x_{ni}=\lim x_{nj}=x, i\not= j.
$$
are called the bubble points.
We consider the bubble point of second type.
Suppose $p$ is a bubble point and $\lim x_{n1}=\lim x_{n2}=p$.
Take unit disk $B_p(1)$ but treat
$(A_n, x_{n1}, x_{n2})$ as generalized connections. Using the
construction of bubble tree compactification \cite{PW}
 such as dilating, transition  etc.,
we get a  marked point instanton on sphere and so on. It is easy to
see that if we fix $A_n=A$ the compactification
coincides with FM compactification of weighted cases.

We can mimic \S 3.1 to define the compactified space which
is consistent with the bubble tree compactification. Some definitions
should be changed slightly:
\begin{enumerate}
\item
Marked point bubble trees:
Besides the general notations we mentioned in \S 3.1,
$w_i$ are assigned to marked points.
$W(v)$ is defined to be sum of $w_i$'s assigned to $v$
and $W(v_i)$ as defined in \S 3.1.  Then definition 3..1
can be used to define {\it marked point bubble trees}.
\item
Marked point bubble space of a marked point
bubble tree: suppose $X$ is the bubble space
assigned to the bubble tree.
For each vertex $v$ suppose the
assigned manifold is $X_v$. For all $w_i$
assigned to $v$ they define distinct marked points $q_i\in X_v$
which are different from the bubble points on
$X_v$, and $w_i$ are weights assigned to these points $q_i$.
\item
Marked point bubble tree instanton:
treat $q_i$ as part of generalized instanton on $S^4_v$
and use the definition 3.3.
\end{enumerate}
We can similarly define $\overline{\mathcal{M}}_{K,\mathfrak{w}}(X)$
and etc.
The gluing theory can be extended to these spaces similarly.

We now introduce some terminologies that are related
 to ghost strata.
Given a ghost tree $T$, we suppose $v_1,\ldots,v_l$ are
all ghost vertices of $T$. If $v_i$
satisfies  $W(v_i)= \min_{1\leq j\leq l}W(v_j)$,
we call $v_i$ to be an {\it ends} of $T$ and $W(v_i)$
to be the energy of the end.
Suppose $T=(V,D)$ is ghost tree and $v$ is an end.
Let $d(v)$ be the set of vertices that are descendants of $v$ and
$V'= V\setminus d(v)$. $T'=(V',D')$ is the subtree of $T$ induced
by $V'$, moreover assign $W(v)$ of $T$ to $v$ to make $T'$
 a marked point bubble tree. We say $T'$ is obtained by
{\it cutting $T$
at $v$}.

We now consider
a  typical  ghost tree  $(T,v_0)\in \mathcal{T}_K$:
$child(v_0)= \{v_1,\ldots, v_k\}$ and $v_i, 1\leq
i\leq k$ are ghost; For each $v_i$
$child(v_i)= \{v_{i1}, \ldots, v_{il_i}\}$
and  $v_{ij}$ are leaves. We now describe
the stratum $\mathcal{S}_T$.
Let
$$
\mathfrak{w} := (w_0, w_1,\ldots, w_k):=
(w(v_0), \sum_{j=1}^{l_1} w(v_{1j}),
\ldots, \sum_{j=1}^{l_k}w(v_{kj}))
$$
 and $K= \sum_{i=1}^k w_i$. It is clear that $\mathcal{S}_T$ is a fiber
bundle over $\mathcal{M}_{K, \mathfrak{w}}(X)$. The fiber is contributed
by two components. One component is bubble points on ghost bubbles. This
can be described as follows: let $\pi: \mathcal{M}_{K,\mathfrak{w}}(X) \to
F_\mathfrak{w}(X,k)$ be the projection map on marked points. Define $$
S_T= \pi^* \prod_{i=1}^{k}SN_{\mathfrak{w}'_i}, $$
$\mathfrak{w}'_i=(w(v_{i1}),\ldots,w(v_{il_i}))$. The other component
comes from the moduli spaces on $v_{ij}$. Let
 $$ Z_T = \prod_{i=
1}^{i=k}\prod_{j=1}^{l_i}\mathcal{M}^b_{w(v_{ij})}. $$ 
and $\mathcal{Z}_T
= \mathcal{M}_{K,w}(X)\dot{\times} Z_T$. Then $$ \mathcal{S}_T =
S_T\times_{\mathcal{M}_{K,\mathfrak{w}} (X)}\mathcal{Z}_T. $$ Since
$SN_{\mathfrak{w}'_i}$ are (almost) sphere bundles, $\mathcal{S}_T$ are
product of $k$ sphere bundles over certain base. This structure is the key
to our flip resolutions. We explain how this can be done assuming that
$SN_{\mathfrak{w}'_i}$ are indeed sphere bundles:  first, note that The
gluing bundle $\mathbf{GL}_T$ whose fiber is $Gl_T$ is a pull back bundle
$\mathbf{GL}'_T$ over $\mathcal{Z}_T$. In other words, it is trivial
over fiber of $S_T$. The neighborhood of $\mathcal{S}_T$ therefore is $$
U(\mathcal{S}_T) \cong \mathcal{S}_T\times_{\mathcal{Z}_T}
\mathbf{GL}'_T/\Gamma_T. $$ For simplicity, we assume $k=1$. Then
$\Gamma_T=SU(2)$ up to some finite group. The flip resolution on this
model is exactly the same as what we did in \S 3.1.

\subsection{Flip Resolutions on $\overline{\mathcal{M}}
_K(X)$}
We already see how a single flip resolution
can be applied to resolve singularities. 
But when $K>2$, ghost strata are complicated.
We need a systematic process to apply flips.
The arrangement is given in the proof of the
following theorem. $K=3$ is a representative
 case to see
how this works. So along each step in the proof, 
the case of $K=3$
is associated as an example. We introduce some 
notations for $\overline{\mathcal{M}}_3(X)$ first:
there are seven ghost strata. The trees are 
$$
\begin{array}{ll}
T_1 = [\bar{0}[0[0[1,1],1]]]; &
T_2 = [\bar{0}[0[1,1[1]]]]; \\
T_3 =  [\bar{0}[0[1,2]]]; &
T_4 = [\bar{0}[0[1,1,1]]];\\
T_5 = [\bar{0}[1[0[1,1]]]]; &
T_6 = [\bar{0}[0[1,1],1]];\\
T_7 =[1[0[1,1]]]. &
\end{array}
$$
Other non-ghost trees are
$$
\begin{array}{llll}
T_{01}=[3]; & T_{02}= [2[1]]; & 
T_{03}=[1[1,1]]; &\\
T_{04}=[1[2]]; & T_{05}=[1[1[1]]]; & &\\
T_{1i}=[\bar{0}[T_{0i}]], & 1\leq i\leq 5; & &\\
T_{21} =[\bar{0}[1,2]]; &
T_{22}=[\bar{0}[1,1[1]]]; &
T_{23}=[\bar{0}[1,1,1]]. &
\end{array}
$$
\begin{theorem}
For any integer $K>0$ by finite steps of "flips",
the singularities of $\overline{\mathcal{M}}_K(X)$
can be resolved.
Denote the new spaces by $\underline{\mathcal{M}}_K(X)$.
\end{theorem}
\noindent
{\bf Proof: }
We prove this by induction for both
$\underline{\mathcal{M}}_K(X), \underline{\mathcal{M}}^b_K$.
The first nontrivial case is
$K=2$. This is proved in \S 4.1.

Now consider general $K$.
The resolutions are done inductively on
the end energies  $m\leq K$ of ghost trees.
\\
1, $m=2$.\\
Pick up  bubble trees in $\mathcal{T}_K$
which has at least one end with energy $m=2$.
Suppose $T$ is such a tree with $k$ such ends
and $v_1,\ldots, v_k$ are these ghost vertices.
Set $V=(v_1,\ldots,v_k)$.
Let $T_1$ be the $k-$marked bubble
tree obtained by cutting off these ends from $T$.
Let $\mathfrak{w}=(K-2k, 2,2,\ldots,2)$.
Then
$\mathcal{S}_{T_1,\mathfrak{w}}
\subset \overline{\mathcal{M}}_{K,\mathfrak{w}}$.
As before 
$$
\mathcal{S}_T
= S_T\times_{\mathcal{S}_{T_1,w}}
\mathcal{Z}_T.
$$
Here $S_T$ is the product of sphere bundles
$SN_i,1\leq i\leq k,$
over $\mathcal{Z}_T$, where the sphere structure
of $SN_i$ is obtained from the ghost vertex
$v_i$.
Actually, this is true only for a proper subset
$U\subset \mathcal{S}_{T_1,w}$. For simplicity we take it
as $\mathcal{S}_{T_1,w}$ right now. We write
$$
Gl_T/\Gamma_T = (Gl_{T,v_1}\times\cdots\times
 Gl_{T,v_k})\times Gl_{T-V}=: Gl_{T,V}\times Gl_{T-V},
$$
where
$$GL_{T,v_i}= R^4_0/\mathbb{Z}_2\times \prod_{j=1}
^{k_i}(R^4_j/\mathbb{Z}_2),$$
and $Gl_{T-V}$ is defined in terms of
(4.7). In this case $k_i=2$.
Correspondingly, there are bundles
$\mathbf{GL}_{T,V}$ and $\mathbf{GL}_{T-V}$ over
$\mathcal{Z}_T$ such that
$$
\mathbf{GL}_T
=\{\mathcal{S}_T\times_{\mathcal{Z}_{T}}\mathbf{GL}_{T,V}
\times_{\mathcal{Z}_{T}}\mathbf{GL}_{T-V}\}/
(\prod_i SO(3)_{v_i}\times \Gamma_{T-V})
$$
where $\prod_i SO(3)_{v_i}\times \Gamma_{T-V}=\Gamma_T$.
We can apply $k$ flip resolutions 
simultaneously to k-pairs $(SN_i, Gl_{v_i})$.
So the stabilizer for the new model is $\Gamma_{T-V}$.
As we pointed  out, the flip resolutions only
apply to the proper set of each stratum. 
Namely, we apply the flip resolutions to each chart given in 
$\mathcal{D}(X,K)$. We need that  resolutions
are compatible on different charts which intersect.
This compatibility is promised by the compatibility
condition {\bf R3'} in \S 3.3.

After doing flip resolutions, we
update $\mathcal{T}_K$:
\\
$[1]$, get rid of trees whose strata
appear in $\overline{\mathcal{S}}_T$.
\\
$[2]$, distribute  the exceptional divisors
to those corresponding trees.
We now explain how to do this.
Suppose $T$  has only one end to be
flipped. When the number of
flip resolutions is more than one,
the method is same. Suppose $v$
is the corresponding ghost vertex.
Let $v_{-1}$ be the parent of $v$
and $child(v)=\{v_1,\ldots,v_k\}$.   Then
$$
Gl_{T,v}=R^4_0/\mathbb{Z}_2
\times (R^4_1/\mathbb{Z}_2\times \ldots \times R^4_k/\mathbb{Z}_2).
$$
Here $R^4_0$ is the basic gluing parameter assigned to
the edge $(v_{-1},v)$ and $R^4_i$ are assigned to $(v, v_i)$.
The assignment of trees to points in  the exceptional divisor
is determined by its coordinates in $SGL_{T,v}$, the sphere
in $Gl_{T,v}$, in a natural way: suppose 
the coordinate of the point is $[x_0,\ldots, x_k]$.
Here $[\ast]$ is the equivalence class with respects to
 $\mathbb{Z}_2$'s
and isotropy group $SO(3)$ actions.
If $x_i\not= 0, i\geq 0$, we contract $T$  at its corresponding
edge. The resultant tree is assigned to the point.
% Since the strata for these new $T$ are different now.
%We would like to make marks. These are three cases:
%\begin{itemize}
%\item Only $x_0\not= 0$. We mark  vertex $v_{-1}$ red and
%call it red  vertex. Note that this corresponds to
%the blow-down in the remark .
%\item $x_0=0$. Now $v$ is no longer trivial vertex. We mark
%$v$ blue  and call it blue vertex.
%\item else. we mark $v_{-1}$ blue.
%\end{itemize}
%\vskip 0.1in
%\noindent
%{\bf Lemma *: }{\it
%Let $T$ be a basic bubble tree and $T'$ is the same type
%of $T$. Now Suppose $T'$ is marked.
%For any vertex $v\in T'$ if non of ascendant of $v$
%is marked
%black, then the  gluing parameters of $v$
%of original $T'$ can be extended over the marked $T'$.
%If $v$ is a ghost vertex,
%$\mathcal{S}_{T'}$ is a sphere bundle over
%$\mathcal{Z}_{T',v}$.
%}
%\vskip 0.1in
%\noindent
%We prove this lemma later.
We now illustrate this step for $K=3$. The ghost strata
to flip when $m=2$ are those of trees $T_1,T_5, T_6, T_7$.
Among them, only $\mathcal{S}_{T_1}$ is compact.
So for the other three strata, their charts only cover
some proper sets. The resolutions for these three strata
only resolve the singularities on these proper sets.
The resolution
of the rest  is done by the resolution
on $\mathbf{GL}_{T_1}$. They would patching together well
because of the compatibility condition. $T_5,T_6,T_7$
are replaced by certain exceptional divisors and they are no longer
singular. The  interesting one left is $T_1$.
 There are two components
in the exceptional divisor which are still singular. 
One is added to  $\mathcal{S}_{T_4}$ and the other component
is added to $\mathcal{S}_{T_3}$ (or $\mathcal{S}_{T_2}$).

Define the resulting space as $\overline{\mathcal{M}}_{K,2}(X)$.
\\
2, $m\leq K$.\\
Pick up bubble trees in $\mathcal{T}_K$
which has at least one end with energy $m$.
Suppose $T$ is such a tree.
For simplicity, we assume
$T$ contains only one such end. For multiple ends
the argument is same
as what we do for $m=2$.
Suppose $v$ is the ghost
vertex. Let
$child(v)= \{v_1,\ldots, v_k\}, w_i= w(t(v_i)),$ and $\mathfrak{w}'=
(w_1,\ldots,w_k)$. We make a further assumption
that $v_i$ are leaves.
The reason we can make such an assumption
is that by induction, $\mathcal{M}^b_{w_i}$, corresponding
to $t(v_i)$,
are completely resolved.
 So we can assume that
$\underline{\mathcal{M}}^b_{w_i}$ is assigned to vertex
$v_i$. For instance, in our example  $\mathcal{S}_{T_2}$
can be treated as a submanifold of $\mathcal{S}_{T_3}$
in this sense. 
By assumption $\sum w_i=m$.
Set $\mathfrak{w}= (K-m,m)$. Let $T_1$ be the $1$-marked  point
bubble tree obtained by cutting $T$ at $v$.
We still have sphere bundles $S_T$ and $\mathcal{Z}_T$,
which it is better to denote by $\underline{\mathcal{Z}}_T$,
over
$\mathcal{S}_{T_1,w}$.
One important fact is that $S_T$ is(!) a sphere bundle.
Note that this is not true until we finish previous
$(m-1)$-step filp resolutions which provide all necessary
blow-downs. In fact, if not,  there must
exist some ghost stratum in lower level which has an end with
energy less
than $m$. But this is impossible.
So we can continue the flip
resolution to the neighborhood of $\mathcal{S}_T$ as before.
By finite steps, all ghost trees should disappear from
$\mathcal{T}_K$.

As an example, let us consider $K=3$ again. 
The ghost strata to resolve for $m=3$ are 
$T_2,T_3, T_4$. As we said, $T_2$ and $T_3$ can be treated 
together to be one stratum. We still denote it by 
$\mathcal{S}_{T_3}$. 
For $T_4$, originally $S_{T_4}$ is almost a sphere bundle 
but it is not a sphere bundle! To see this, 
it is $B_{123}$ in example 4.7
 diagonals are replaced
by the exceptional divisors of semi-blow-ups. Applying
the flip resolutions to the step
for $m=2$, the model now corresponds
to $X^3_3$.
these blow-ups are killed by blowing-down. 
So $S_{T_4}$ is now a sphere bundle.  It actually corresponds
to $B_{123}$. Then we again apply flip resolutions 
to these strata. This completes the flip resolutions for $K=3$.
q.e.d.
\vskip 0.1in
We now derive some consequences from our smooth compactified
moduli spaces.
First consider 
the $SO(3)$-bundle $\mathcal{M}^{b,0}_K$ over $\mathcal{M}^{b}_K$.
After compactification, 
$\overline{\mathcal{M}}^{b,0}_K$
is a fibration over $\overline{\mathcal{M}}^b_K$, but the
bundle structure fails exactly at where the ghost bubbles occur
at the principal component.
By exactly same  constructions, 
the flip resolutions can be applied to both spaces. 
For greater convenience later on,
we consider a bundle $\mathcal{RM}^{b,0}_K$, which is a semi-blow-down
of $\mathcal{M}
^{b,0}_K\times R^+$. Basically, we just add a $R^+$ factor to 
the fiber $SO(3)$ and make it to be $R^4/\mathbb{Z}_2$. Note that
this bundle plays a role as part of gluing parameter
in $\overline{\mathcal{M}}(X)$. So by the same construction
given in theorem 4.12,
one can construct smooth orbifold bundle $\underline{\mathcal{RM}}^
{b,0}_K$ over $\underline{\mathcal{M}}^b_K$. The fiber is 
$R^4/\mathbb{Z}_2$. Take its sphere bundle and call it
$\underline{\mathcal{M}}^{b,0}_K$. 
It might be more natural to resolve $\overline{\mathcal{M}}
^{b,0}_K$ directly. In fact, this  gives the same bundle. 
As we said, the $SO(4)$ action on $R^4$
induces  an action on $\mathcal{M}_K^b$ and $\mathcal{M}_K^{b,0}$.
It is not difficult to check that 
the actions can be extended to all compactification spaces. So we have
\begin{corollary}
 $\underline{\mathcal{M}}^{b,0}_K$
is an $SO(3)$-bundle over $\underline{\mathcal{M}}^b_K$.
Moreover, the bundle is $SO(4)$-equivariant.
\end{corollary}
We now consider some useful
sub-orbifolds in $\underline{\mathcal{M}}_K(X)$ 
and their normal bundles.
Let $0<l\leq K$. We put all $l$-level strata 
together and  denote the set by 
$\mathcal{S}^l$.
Each component of $\mathcal{S}^l$, given
by a partition $\mathfrak{t}$
of $l$, is a sub-orbifold and  denoted 
by $\mathcal{S}^l_\mathfrak{t}$.
From the construction of $\underline{\mathcal{M}}_K(X)$,
we have
\begin{corollary}
Suppose $\mathfrak{t}=(l_1,\ldots,l_k)$
is a partition of $l$. In $\underline{\mathcal{M}}_K(X)$ 
$$
\mathcal{S}^l_\mathfrak{t}
=[\mathcal{M}_{K-l}(X)
\times (Fr^k(X)\times_{SO^k(4)}
\prod_{i=1}^k \underline{\mathcal{M}}_{l_i}^b)]/S_\mathfrak{t}.
$$
Its normal bundle, denoted by $\mathcal{NS}^l_\mathfrak{t}$, is 
$$
[(\tilde{\mathcal{M}}_{K-l}(X)\times_{\mathcal{G}_{P_{-l}}}
\times P^k_{-l})\times_{X^k} Fr^k(X))\times_{SO^k(4)\times
SO^k(3)} \prod_{i=1}^k 
\underline{\mathcal{RM}}_{l_i}^{b,0})]/S_\mathfrak{t}.
$$
In particular, when $\mathcal{M}_{K-l}$ is a point. 
\begin{equation}
\mathcal{S}^l_{\mathfrak{t}}
=
[Fr^k(X)\times_{SO^k(4)}
\prod_{i=1}^k \underline{\mathcal{M}}_{l_i}^b]/S_\mathfrak{t}.
\end{equation}
and 
\begin{equation}
\mathcal{NS}^l_\mathfrak{t}
= [\prod_{i=1}^k (P_{-l}\times_X Fr(X))\times_{SO(4)\times SO(3)}
\underline{\mathcal{RM}}^{b,0}_{l_i}]/S_\mathfrak{t}.
\end{equation}
\end{corollary}
Later, for the model of the wall-crossing formula,
$\mathcal{M}_{K-l}$ is not a point. But we 
only  take a neighborhood of a point. 
So the models we consider are  products of $\mathbb{C}^N$
with (19) and (20). Note that these two models 
are very nice. Essentially, they have product structures
and so usually the computations are reduced to  
the case that $\mathfrak{t}$ contains only one element.

\section{Kotschick-Morgan Conjecture on Donaldson 
Wall-crossing Formula }
In the second part of the paper
we  demonstrate
how the new compactifications  can be applied.
As an example,  we prove 
the KM-conjecture in this chapter.
It turns out that the localization theory of equivariant
cohomologies is very helpful to the wall-crossing formula
when we working on  smooth spaces.
The similar situation to the Seiberg-Witten theory 
was studied by Cao-Zhou (\cite{CZ}). 
Recall that the wall-crossing formula is expressed in terms of
the link. In \S 5.1, we explain how links are related to
the equvariant integration. Then in \S 5.2 \& 5.3, we apply 
 the localization formula
 to our new compactified moduli spaces and prove the 
conjecture.

\subsection{Review on Equivariant Theory and Localization Formula}
In this section, we simply review  the equivariant theory and 
concentrate on the localization formula. Readers are referred to
\cite{AB},\cite{BGV} and \cite{K} for details.

Suppose $G$ is a compact Lie group and $M$ is a topological space.
Set $M_G=EG\times_G M$, where $EG\to BG$ is the universal $G$-bundle.
The {\it equivariant cohomology} of $M$ is defined by
$$
H^\ast_G(M):=  H^\ast(M_G).
$$
Since $EG\times_G
M$ is a bundle over $BG$ with fiber $M$, $H^\ast_G(M)$
is a module over $H^\ast(BG)$. In particular, 
$H^\ast(BS^1)=\mathbb{C}[u], H^\ast(B\mathrm{Spin}(4))=
\mathbb{C}[c_R,c_L]$ (see lemma 5.10). 
Suppose $\pi: E\to M$ is a $G$-equivariant vector bundle.
The {\it equivariant Thom class} is defined to be
the Thom class of $\pi_G: E_G\to M_G$.

If $M$ is smooth, the equivariant de Rham  theory is also available.
Assume $dim M=2m$.
From now on, we only consider $G=S^1$.
Let $A^\ast(M)$ be the graded algebra of differential forms of $M$.
Counting $u$ as an element of
degree 2, we grade polynomial ring $A^\ast(M)[u]$
as well. 
So a degree $2k$ polynomial is in the form
\begin{equation}
\gamma = a_0 u^k + a_1u^{k-1} +\cdots + a_k.
\end{equation}
We call $a_0 u^k$ is the {\it leading term} of $\gamma$. Define
$\gamma_{[n]}= a_n$, the $2n$-form coefficient.
Suppose $V$ is the vector field on $M$ generated by the $S^1$ action.
Define an operator $d_{S^1}: A^\ast(M)[u]
\to A^\ast(M)[u]$ by
$$
d_{S^1}\gamma = (d- i_V)\gamma, \mbox{  } \gamma\in A^\ast(M)[u].
$$
We say that $\gamma
$ is an {\it $S^1$-equivariant form} if $\gamma$ is
$d_{S^1}$-closed. If 
 $\gamma$ has degree 
$2k$ given in the form as in (21), it is equivariant iff
$$
i_V a_n = d a_{n-1}, 1\leq n\leq k.
$$
These definitions can be easily generalized to $A^\ast(M)[u,u^{-1}]$
where $u^{-1}$ has degree $-2$.

We now explain the localization formula.
From now on, we assume that all equivariant classes are represented
by  de Rham forms.
Suppose $\mathbf{F}=\cup_i F_i$
is the union of all components of stationary sets of the action.
Let $N_i$ be   the $S^1$-equivariant
normal bundle of $F_i$ in $M$ and
$\Theta_i$ be its equivariant Thom classes. 
Then the equivariant Euler class of the  bundle $N_i$
is $E_i=i^\ast_{F_i} \Theta_i$.
We endow $M$ with an $S^1$-invariant 
metric. Define a 1-form $\theta$
on $M-\mathbf{F}$ such that 
$\theta(V)=1, \theta|_{V^\perp}=0$.
Here $V^\perp$    is the orthogonal complement
of $V$ in the tangent space.
Then $\theta$ is a connection on the principal
bundle $M-\mathbf{F}\to (M-\mathbf{F})/S^1$.
Equivariant forms have some properties over 
$M-\mathbf{F}$. We list two useful  lemma
here.
By a direct computation, one has
\begin{lemma}
Suppose 
$\gamma = \mu + fu$ is an equivariant
two form, then $r(\alpha)=\mu -\theta\wedge df+fd\theta$
is an $S^1$-invariant form over $M-\mathbf{F}$.
\end{lemma}
A simple result for equivariant forms with general degree
is also true (\cite{CZ}). We only need degree $2$ forms 
 in our application.
\begin{lemma}
If $\gamma$ is an equivariant form. Then 
$\gamma_{[m]}$ is exact outside $\mathbf{F}$.
\end{lemma}
\noindent
{\bf Proof: }
Following the definition, one can check that
$$
\gamma_{[m]}
=d(\theta\wedge\gamma_{[m-1]}
+\theta\wedge d\theta\wedge \gamma_{[m-2]}
+ \theta\wedge (d\theta)^2\wedge \gamma_{[m-3]}+\cdots ).
$$
This proves the lemma.q.e.d.
\begin{theorem}
[Atiyah-Bott] 
If $\gamma\in H^\ast_{S^1}(M)$,
\begin{equation}
\int_M \gamma = \sum_i\int_{F_i}\frac{\gamma}{E_i}.
\end{equation}
\end{theorem}
Note that the degree of $E_i$ is equal to the codimension of 
$F_i$ and its leading term does not vanish. The right hand 
side of (22) makes sense (in $A^\ast(F_i)[u,u^{-1}]$).
The same situation also holds  when we  restrict  $\Theta_i$
over a neighborhood of $F_i$.
\\
{\bf Proof: }
Let $U_i,U_i',U_i^{''}$ be $S^1$-invariant 
tubular neighborhoods of $F_i$ with
relations $U_i\supset U_i'\supset U_i^{''}$. Suppose that   $\Theta_i,
\Theta_i'$ are the Thom classes of the normal bundles and they are
supported in $U_i, U_i'$ respectively. Moreover, we require that
\[
\Theta_i|_{U_i^{''}}=\Theta_i'|_{U_i^{''}}.
\]
Let $\gamma_i=i^*_{F_i}\gamma$.
Compare
$$
\sum_{i}\int_{F_i}\frac{\gamma_i}{E_i}=\sum_{i}\int_{U_i}\frac{\gamma_i}{\Theta
_i}\wedge
\Theta_i',
$$
 to $\int_M \gamma$, we have
\begin{equation}
\sum_{i}\int_{U_i}\frac{\gamma_i}{\Theta_i}\wedge \Theta_i' -\int_M\gamma=
\sum_{i}\int_{U_i}\gamma(\frac{\Theta_i'}{\Theta_i}-1)-\int_{M\setminus
\cup
U_i} \gamma,
\end{equation}
note that $\frac{\Theta_i'}{\Theta_i}-1$
are supported  in
$M\setminus U_i^{''}$ and are equal to $-1$ when in $M\setminus U_i'$,
 thus $\gamma':=\sum_{i}\frac{\gamma_i}{\Theta_i}\wedge
\Theta_i'
-\gamma$ is an well defined   equivariant form which is supported away
from all
stationary sets.  By lemma 5.2, $\int_M \gamma'=0$. This implies that
the right hand side of (23) is 0. So we prove (22). 
q.e.d.
\vskip 0.1in
\noindent
The Localization theorem also has a version for
manifolds with boundary (\cite{Ka}, \cite{CZ}).
Suppose $M$ is a manifold with boundary and
$\partial M\cap \mathbf{F}= \emptyset$. 
We state the theorem for a special case
that we need later. 
\begin{theorem}
Suppose
$\gamma=\mu
+ fu$ is an equivariant 2-form 
of $M$ and $r(\gamma)$ is the one defined in lemma 5.1. Then 
\begin{equation}
\int_{\partial M/S^1} r^{m-1}(\gamma)
=\sum_i\int_{F_i} \frac{\gamma^{m-1}u}{E_i}.
\end{equation}
where $2m=\dim M, 2n_i=dim F_i$.
\end{theorem}
\noindent
{\bf Proof: } Let $\tilde{\gamma}= \gamma^{m-1}u$. Consider 
$\int_M \tilde{\gamma}$. By the same argument used in the proof
of theorem 5.3, one also has a form $\tilde{\gamma}'$ supported
outside $\mathbf{F}$. Moreover, $\tilde{\gamma}'=\tilde{\gamma}$
near $\partial M$.
It is sufficient to show
$$
\int_{M} \tilde{\gamma}'= \int_{\partial M/S^1} r^{m-1}(\gamma).
$$
By lemma 5.2, we know that $\tilde{\gamma}'
=d \eta$ for some $\eta$ over $M-\mathbf{F}$. Note
that $\eta$ is trivial over a neighborhood of $\mathbf{F}$. 
And 
in a neighborhood of  $\partial M$, one can show that
$$
\tilde{\gamma}'=\gamma
= d(\theta \wedge r^{m-1}(\gamma)).
$$
Therefore by Stokes theorem
$$
\int_M\tilde{\gamma}' = \int_{\partial M} \theta \wedge r^{m-1}(\gamma).
$$
Since $r^{m-1}(\gamma)$ is horizontal, applying 
integration along fiber to the right side, we show (24).
q.e.d.
\vskip 0.1in
\noindent
In (?), the left hand side corresponds to the ``link'' 
in our application. 
Finally, we make a remark on $\mathbf{F}$. If the stationary
sets are orbifolds, these results are still true up to 
 proper constants depending upon the orbifold structures.

\subsection{Reducible instantons}
Here we explain the local models of reducible connections
in $\underline{\mathcal{M}}_K(X, \lambda)$. 
All the constructions and techniques 
are well described
in \cite{DK}. Here, we state most of results  
without giving proofs. 
Readers are referred to corresponding
theories in \cite{DK}.

First assume the reducible connection $[A]$ with 
respect to $\alpha$ is in the top stratum
$\mathcal{M}_K(X,\lambda)$. It is known that
\begin{lemma}
The neighborhood of $[A]$ 
in $\mathcal{M}_K(X,\lambda)$ is
diffeomorphic to a neighborhood $U_A/\Gamma_A$ of $0$ in
$\mathbb{C}^N/\Gamma_A$,
where $N= -\alpha^2-2$ and $S^1$ action is complex multiplication. 
\end{lemma}
\noindent
{\bf Proof: } In general, the local model at $[A]$ is given by the
$\Psi^{-1}(0)$ for  a  $\Gamma_A$-equivariant map
$$
\Psi: H^1_{[A]} \to H^2_{[A]}.
$$
When $b^+_2=1$, $H^2_{[A]}$ is 1-dimensional. By choosing a
generic path $\lambda$, 1-dimension parameter $t$ can be mapped
onto $H^2_{[A]}$ (\cite{KrM}). So the lemma follows. q.e.d.
\vskip 0.1in
\noindent
The universal bundle $\mathbb{P}$ can  not be defined over 
$U/S^1\times X$. However, there is a bundle
$\tilde{\mathbb{P}}\to U_A\times X$ with compatible $\Gamma_A$
action such that $\mathbb{P}= \tilde{\mathbb{P}}/\Gamma_A$
over $U-\{0\}/S^1$. Our main strategy is to apply
the equivariant theory over $U_A$.
\begin{lemma}
For any $\Sigma\in H_2(X,\mathbb{Z})$,
the restriction of $\tilde{\mu}(\Sigma)$
to $U_A$ can be extended as an equivariant
class $\tilde{\mu}(\Sigma) + fu$ such that
$$
f|_{\{0\}}= -\frac{1}{2}\langle \alpha, \Sigma\rangle,
$$
\end{lemma}
The proof is exactly the same as the one in \S 5.1.4 (\cite{DK}).
Alternatively, one can construct $\Gamma_A$-equivariant 
line bundle representing 
$2\tilde{\mu}(\Sigma)$ over $U$. 
$\Gamma_A$ acts weightedly on the line bundle (\S 5.2 in \cite{DK}).
Then the  equivariant Chern class of the line bundle is 
$2\tilde{\mu}(\Sigma) -\langle \alpha, \Sigma\rangle u$.

Now consider 
the reducible connection $[A]$ with respect to 
$\alpha$ located in lower strata. Suppose
$r= (\alpha^2- p_1(P))/4$. The family of reducible connections
in $\overline{\mathcal{M}}^u_K(X,\lambda)$
is $[A]\times Sym^r(X)$. Let $N=-\alpha^2-2$ and $U$
is a neighborhood of $\{0\}$ in $\mathbb{C}^N$.
When $\alpha^2=-1$, one has to consider the obstruction
bundle. 
The localization formula can still be applied to the 
so called "virtual neighborhoods" \cite{CLRT}.  
The main result, proposition 5.9, is still true for this case.
Now assume $\alpha^2\not= -1$.
The local model of $[A]\times Sym^r(X)$
can be written as 
$(U_A\times \mathcal{GL}_{\alpha,r})/\Gamma_A$.
In general, 
$\mathcal{GL}_{\alpha,r}$ is rather complicated (\cite{FM},\cite{L}).
Given a $\Sigma\in H_2(X,\mathbb{Z})$, there is a line bundle
$L_\Sigma\to X$ whose $c_1$ is $\omega$, the Poincare dual 
of $\Sigma$ in $H^2(X,\mathbb{Z})$. For a line bundle 
$L\to X$, $L^n/S_n$ is a line bundle over $Sym^n(X)$. We denote
the bundle by $Sym^n(L)$. $Sym^n(\omega)$ is the $c_1$
of $Sym^n(L_\Sigma)$. 
Mimicking the construction in \S 7 of \cite{DK}, we have
\begin{prop}
There is a $\Gamma_A$-equivariant line bundle
$\mathcal{L}_{\alpha,r}$ 
over $U\times \mathcal{GL}_{\alpha,r}$ which is 
isomorphic to $\mathcal{L}_{\alpha,0}\otimes
Sym^r(L_{\Sigma}^2)$. 
$\mathcal{L}_{\alpha,r}/\Gamma_A$
restricted to $(U\times \mathcal{GL}_{\alpha,r}/\Gamma_A)
\setminus [A]\times Sym^r(X)$ represents $2\bar{\mu}(\Sigma)$.
\end{prop}

\subsection{The Proof of Kotschick-Morgan Conjecture}
In this section, we  apply the localization formula 
given in theorem 5.4 to prove the KM-conjecture.

As we mentioned, all properties in $\overline{\mathcal{M}}
^u_K(X,\lambda)$ can be lifted to both
$\overline{\mathcal{M}}_K(X,\lambda)$ and
$\underline{\mathcal{M}}_K(X,\lambda)$ by $\pi_i^\ast,
i=1,2$. Let $\mathcal{U}_{\alpha, r} = \pi^{-1}_2 (U_A\times
\mathcal{GL}_{\alpha,r})$. To apply theorem 5.4 to $\mathcal{U}_{\alpha,r}$
we first study the fixed loci of $\Gamma_A$ 
and how $\bar{\mu}(\Sigma)$'s are equivariantly extended.
\begin{prop}
 The fixed loci of 
$\Gamma_A$ in $\mathcal{U}_{\alpha,r}$ are disjoint smooth
components $F_\mathfrak{t}:=
[A]\times \underline{\mathcal{S}}^r_\mathfrak{t}(X)$,
where $\mathfrak{t}=(t_1,\ldots, t_l)$ runs over
all  partitions
of $r$.
The pull-back bundle $\pi^\ast_2\mathcal{L}_{\alpha,r}$
is a $\Gamma_A$ equivariant line bundle. Its equivariant
chern class restricted to the  fixed loci 
$F_{\mathfrak{t}}$ is 
 $\langle \alpha, \Sigma \rangle u +2 
Sym^r(\omega)$.
\end{prop}
\noindent
{\bf Proof: } 
We explain the first statement for the case $r=2$. For 
$r>2$, the phenomenon are same.

We begin with the model in $\overline{\mathcal{M}}_{K}(X,
\lambda)$. 
We explain how $\Gamma_A$ acts in this model and show
how the situation changes after applying flip resolutions.
Let $\mathcal{U}_1:= \pi^{-1}_1 (U\times
\mathcal{GL}_{\alpha,r})$. The partitions of $r=2$
are $\mathfrak{t}_1=(1,1)$ and $\mathfrak{t}_2=(2)$.
The reducible connections 
are parameterized by two components
 $F_1: =[A]\times \overline{Sym^2(X)\setminus\Delta}$ and
$F_2:= [A]\times \overline{X\dot{\times} \mathcal{M}_2^h}
$. Here $F_i$ corresponds to the partition $\mathfrak{t}_i$.
 Note that $F_1$ intersects with $F_2$
at ghost bubbles. 
 We now clarify the situation 
around ghost bubbles. 
The ghost stratum can be identified as 
$(U\times STX)/\mathbb{Z}_2$.
The gluing parameter 
of this stratum is 
$$
Gl= \prod_{i=0}^2 R^4_i/\mathbb{Z}_2.
$$
The gluing bundle
$\mathbf{GL}$ is the bundle over 
$U\times STX$ with fiber isomorphic to $Gl$.
In fact, it is constructed by a pull-back bundle from $X$
$$ 
\pi^\ast (\mathcal{M}^0_{K-2}(X)\times 
\mathcal{M}^{b,0}_1(S^4,X))^2\times R^+
\times R^+\times R^+.
$$
where $\pi: U\times STX\to X$ is the projection map.
Then the local model of the ghost stratum  is
$
\mathcal{U}_1 \subset \mathbf{GL}/(\mathbb{Z}_2\times SO(3)).
$
$\Gamma_A$  acts only on the first factor of $Gl$.
In general, it acts only on the gluing parameters   
associated to bubble points on $X$. 
Suppose $(x_0,x_1,x_2)\in Gl$. The fixed loci of 
$\Gamma_A$ in $\mathcal{U}_1$ are
those corresponding to $(0,x_1,x_2)$
and $(x_1,0,0)$. The first case follows from the 
action directly, the second case comes from the $SO(3)$
quotient.
They are indeed the 
neighborhoods of ghost stratum in $F_2, F_1$.
 This picture says that 
two different  fixed loci $F_1,F_2$ 
intersect at the ghost stratum.  
Clearly, this awkward situation is caused by the singularity
at the ghost stratum. Now we explain how the situation
is modified after flip resolutions.
After the resolution, the local model is
$$
\mathcal{U}= [U\times (TX/\mathbb{Z}_2)\dot{\times}
S^{11}/SU(2) ]/\Gamma,
$$ 
here $\Gamma$ is isomorphic to $\Gamma_A$
up to some finite group.  
By the definition of flip resolutions,
points in $S^{11}/SU(2)$ is
$
[x_0,x_1,x_2]
$ where $(x_0,x_1,x_2)\in S^{11}$ 
and the
equivalent relation is given by
$(x_0, x_1,x_2)
\sim (c^{-1} x_0,c x_1,c x_2)$ for some $c\in SU(2)$.
Now the fix loci of $\Gamma_A$ are separated. They correspond
to $[0,x_1,x_2]$ and $[x_0,0,0]$ respectively. On the exceptional
divisor after flip, $X\times [x_0,0,0]$ is just the natural closure
of $Sym^2(X)\setminus \Delta$. These two 
components are in $\underline{\mathcal{S}}^2_{\mathfrak{t}_i},i=1,2$.
This proves the first statement. For the second 
part, we use the line bundle given in proposition 5.3
and pull it back to $\mathcal{U}_{\alpha,r}$.
One can show that the weight of $\Gamma_A$-action
on factor $L_\Sigma$ is trivial. Then the assertion
about forms at fixed loci follows from the
case of top stratum.
q.e.d.
\vskip 0.1in
\noindent
So far, we have finished the study of equivariant forms
over $\mathcal{U}_{\alpha, r}$.
We take the neighborhood of each fixed  locus 
$F_{\mathfrak{t}}$  and denote it by
$\mathcal{U}_{\alpha,r,\mathfrak{t}}$.
We will apply theorem 5.4 to this model.
Here, $\gamma$ will be the equivariant 
extension  of $\bar{\mu}(\Sigma)$.  
Usually, 
computing the 
equivariant Euler classes of normal bundles at fixed loci
is highly nontrivial. We first recall some facts mentioned
at the end of
\S 4.3. These can be used to reduce the computations.
Recall that if  $\mathfrak{t}=(t_1,\ldots,t_k)$ is a partition of 
$r$, 
the corresponding fix locus is 
$F_\mathfrak{t} \cong \underline{\mathcal{S}}^r_{\mathfrak{t}}(X)$.
Then its normal bundle
is 
$$
N_{\mathfrak{t}}= U_A\times
\underline{\mathcal{NS}}^r_{\mathfrak{t}}(X).
$$
Since these bundles have 
 product structures when $k>1$, we consider
$k=1$ first.
Now the situation is: let $\mathfrak{t}=(r)$,
there is a bundle
$N_\mathfrak{t}$ over 
$[A]\times \underline{\mathcal{S}}_\mathfrak{t}(X)
=[A]\times X\dot{\times}\underline{\mathcal{M}}^b_r$.
The fiber is $\mathbb{C}^N\times c(SO(3))$, where
$c(SO(3))$ denotes the cone of  $SO(3)$.
We now focus on the nontrivial part of the bundle:
the subbundle $N_{\mathfrak{t}}'$ with fiber $c(SO(3))$.
Explicitly, 
$$
N_{\mathfrak{t}}'= [P_\alpha \times 
\underline{\mathcal{M}}^{b,0}_r(S^4,X)]/SO(3).
$$
Also
$$
\underline{\mathcal{M}}^{b,0}_r(S^4,X)
=Fr(X)\times_{SO(4)} \underline{\mathcal{M}}^{b,0}(S^4).
$$
By theorem 5.4, we have
\begin{prop}
Suppose $\theta$ is an equivariant form over 
$\mathcal{U}_{\alpha,r}$.
Take an equivariant neighborhood
of $F_\mathfrak{t}$ in $N_{\mathfrak{t}}$, say
 $U_\mathfrak{t}$. 
Then
\begin{equation}
\int_{\partial U_\mathfrak{t}/S^1} r(\theta)
= 2 \int_{F_{\mathfrak{t}}}\frac{1}{u^N}
\frac{\theta u}{-(u-\alpha)^2 +p_1(r)},
\end{equation}
where $p_1(r)$ is the $p_1$ of $SO(3)$-bundle 
$\underline{\mathcal{M}}^{h,0}_r(S^4,X)$. 
\end{prop}
\noindent
{\bf Proof: } 
The trouble is to deal with  $N'_{\mathfrak{t}}$.
Without  loss of generality, we assume that $N=0$. Now 
$N_{\mathfrak{t}}$ is an orbifold bundle 
with fiber isomorphic
to $R^4/\mathbb{Z}_2$. 
We first make some assumptions so that
 the computation can be done
under the best situation:
(1)
in general, this orbifold bundle is not a global quotient
of some $R^4$-bundle. Let us assume that 
it can be identified as  a global quotient for a moment. 
(2) Suppose both $SO(3)$-bundles $P_\alpha$,
 $\underline{\mathcal{M}}^{b,0}_r(S^4,X)$ can be lifted to 
$SU(2)$-bundles, say $Y_1,Y_2$.
(3) Let $Z_i=Y_i\times _{SU(2)} \mathbb{C}^2$. 
As we know, $Z_1$ can be split as $Z_1= L_\alpha^{1/2}
\oplus L_\alpha^{-1/2}$.  We assume that 
$Z_2$ can be also split as $Z_2= L_2^{1/2}\oplus
L^{-1/2}_2 $ where $p_1(r)=L_2^2$. With these assumptions,
we
consider the cone
of $Hom_{SU(2)}(Z_1,Z_2)$. This is actually 
a $\mathbb{C}^2$ bundle, we denote it by $Hom(Z_1,Z_2)$.
Then
$$
Hom(Z_1,Z_2)\cong (L_\alpha^{-1/2}\otimes L_2^{1/2})
\oplus (L_\alpha^{1/2}\otimes L_2^{1/2}).
$$
$\Gamma_A=S^1$ acts on the bundle by $(e^{i\theta/2}\cdot,
 e^{-
i\theta/2}\cdot)$. One can check that $Hom(Z_1,Z_2)$
gives a double cover of $N_{\mathfrak{t}}$. 
However, this bundle  exists only when our assumptions hold.
In general, we can pass to a quotient bundle of this bundle which 
always exists: let $G_0=\mathbb{Z}_2\times \mathbb{Z}_2$
act on $Hom(Z_1,Z_2)$ by multiplication. Then the quotient
bundle $Hom(Z_1,Z_2)/G_0 = (L_\alpha^{-1}\otimes L_2)
\oplus (L_\alpha\otimes L_2)$ and the $S^1$ action is
given by $(e^{i\theta}\cdot, e^{-i\theta}\cdot)$.
There is a natural map 
$$
\xi: N_{\mathfrak{t}}\to Hom(Z_1,Z_2)/G_0
$$
which is a double cover.
Now, applying the localization formula over $Hom(Z_1,Z_2)/G_0$
is legal.
There is a small problem to deal with:
 $\theta$ is not $G_0$ invariant, however 
we can replace $\theta$ by the average of $g^\ast \theta, g\in G_0$.
This will not change the left hand side of (25). So we assume that
$\theta$ can be reduced to $Hom(Z_1,Z_2)/G_0$.  
The equivariant Euler class of $Hom(Z_1,Z_2)/G_0$ is then
$$
(u -\alpha +  c_1(L_2))(-u+\alpha +c_1(L_2))= -(u-\alpha)^2 + p_1(r).
$$
Counting the multiplicity of covering and apply the localization
formula, we prove (25).
q.e.d.
\vskip 0.1in
\noindent
Now the main problem in the right hand side of (25)
is to compute 
$p_1(r)$.                     To deal with this problem
we use a trick that was essentially used in \cite{L}. 

Suppose $f: X\to BSO(4)$ is the map that induces $Fr(X)$, i.e, 
$f^\ast ESO(4) = Fr(X)$. All associated
bundles of $Fr(X)$ can be obtained by $f^\ast$.
So $\underline{\mathcal{M}}_r^{b}(S^4,X)$ and
$\underline{\mathcal{M}}_r^{b,0}(S^4,X)$
are the pull-backs of bundles
$$
(\underline{\mathcal{M}}_r^{b}(S^4))_{SO(4)}=
ESO(4)\times_{SO(4)} 
\underline{\mathcal{M}}_r^{b}(S^4)
$$
and
$$
(\underline{\mathcal{M}}_r^{b,0}(S^4))_{SO(4)}=
ESO(4)\times_{SO(4)} 
\underline{\mathcal{M}}_r^{b,0}(S^4).
$$
Moreover 
$(\underline{\mathcal{M}}_r^{b,0}(S^4))_{SO(4)}
$ is still an  $SO(3)$-bundle 
over $(\underline{\mathcal{M}}_r^{b}(S^4))_{SO(4)}
$. We denote the $p_1$ of this bundle by 
$\tilde{p}_1$. Then $p_1(r)=(f^\ast)^\ast \tilde{p}_1$.
By definition,
 $\tilde{p}_1$  is 
an equivariant
cohomology of $\underline{\mathcal{M}}^{b}_r(S^4)$.
 As it is suggested in \cite{L}, we work with
$Spin(4)=SU(2)_L\times SU(2)_R$ equivariant cohomology rather 
than with $SO(4)=SU(2)\times_{\mathbb{Z}_2} SU(2)$. Working
with rational coefficients, no information is lost, as evident from
by this lemma.
\begin{lemma}[\cite{HH}]
Let $c_L,c_R$ be the second Chern classes of 
$SU(2)$-
 bundles
$$
ESpin(4)/SU(2)_R
\to BSpin(4)
\mbox{ and }
ESpin(4)/SU(2)_L\to BSpin(4)
$$
 respectively, then 
$H^\ast(BSpin(4),\mathbb{Q})=\mathbb{Q}[c_L,c_R]$.

If $s: BSpin(4)\to BSO(4)$
is the classifying map  for the natural $SO(4)$ bundle over 
$BSPin(4)$ and $e,p_1\in H^4(BSO(4))$ are the universal Euler
 class and
the  universal Pontrjagin class, then 
$s^\ast (p_1+2e)=-4c_R$ and $s^\ast (p_1-2e)= -4c_L$.
\end{lemma}
Since $\tilde{p}_1(r)\in H^4_{Spin(4)}
(\underline{\mathcal{M}}^{b}_r(S^4))$, we have
\begin{corollary}
$\tilde{p}_1(r)= \omega_0 + \omega_L c_L + \omega_R c_R$,
where $\omega_0\in 
H^4(\underline{\mathcal{M}}^b_r(S^4))$
and $\omega_L, \omega_R\in 
H^0(\underline{\mathcal{M}}^b_r(S^4))$.
\end{corollary}
Now we are able to prove our main theorem.
\begin{theorem}
KM-Conjecture  is true.
\end{theorem}
\noindent
{\bf Proof: } 
Suppose $[A]$ is a reducible connection with respect
to $\alpha$ and lies in level-$r$ strata, i.e,
$r=(\alpha^2-p_1(P))/4$.
Using the  same notations as in proposition 5.7, the neighborhood of
reducible family is $\mathcal{U}_{\alpha,r}$. 
Let $\gamma$ be the equivariant extension of
$\bar{\mu}(\Sigma)$. So $\gamma = \bar{\mu}(\Sigma) +fu$.
Given a partition $\mathfrak{t}$ of $r$,
we know $\gamma$ restricted to the  fixed loci $F_\mathfrak{t}$
is  $Sym^r(\omega) +\frac{1}{2} \langle \alpha,\Sigma
\rangle u $. 
Apply theorem 5.4 to this model by letting
$M=\mathcal{U}_{\alpha,r}$. It is not hard to see that 
$$
r(\gamma) = \bar{\mu}(\Sigma) + \frac{1}{2}
\langle \alpha, \Sigma \rangle d\theta.
$$
So use    this information and plug $r(\gamma^d)$
  into the left hand side of (25).
We get exactly  $\delta_P(\alpha)$. 
By the Stokes theorem, one can easily show that  $\delta_P(\alpha)$
is the sum of $\delta_{\alpha,\mathfrak{t}}$ which is defined by
$$
\delta_P(\alpha)= \sum_{\mathfrak{t}} \int_{\partial U_\mathfrak{t}/S^1}
r(\gamma^d)=: \sum_{\mathfrak{t}} \delta_{\alpha,\mathfrak{t}}. 
$$
So it is sufficient to prove the following statement about
 $\delta_{\alpha,\mathfrak{t}}$:
 $\delta_{\alpha,\mathfrak{t}}$ is in the form
$$
\delta_{\alpha,\mathfrak{t}}:
= \sum_{i=0}^r a_i (r,d ,X ,\mathfrak{t}) q^{r-i}\alpha
^{d-2r-2i},
$$
and as it  indicates  $a_i$ depend 
on $r,d \mathfrak{t}$ and homotopy invariants $e,\sigma$
of $X$.

First assume $\mathfrak{t}=(r)$.
Therefore, the right hand side of (25) is the 
constant term of $u$-series
\begin{equation}
C(r)\frac{1}{u^N}\int_{\underline{\mathcal{M}}_r(S^4,X)}
\frac{u^{4r-2}(r\cdot \sigma + \langle \alpha,\Sigma\rangle)^{d-1}}
{-(u-\alpha)^2+p_1(r)},
\end{equation}
where $C(r)$ is some constant depending on $r$ due to 
the orbifold structure. 
Note that $\pi:\underline{\mathcal{M}}_r(S^4,X)
\to X$ is a bundle 
with fiber
$\underline{\mathcal{M}}^b_r$,
we can apply the integration along fiber $\pi_\ast$
 first. The contributions
of forms in fiber direction come from $p_1(r)$. More precisely,
they are terms containing $p_1^{2r-2}(r)$  and $p_1^{2r-1}(r)$.
Since $\pi_\ast$ commutes with the pull-back $(f^\ast)^\ast$, 
we can apply integration along fiber 
at $(\underline{\mathcal{M}}^b_r(S^4))_{Spin(4)}$
first and pass down to $BSpin(4)$, then pull back 
by $f^\ast$. By corollary 5.11,
the integration over fiber operation $\pi_\ast$
at $(\underline{\mathcal{M}}^b_r(S^4))_{Spin(4)}$
implies
$$
\pi_*(\tilde{p}_1(r))^{2(r-1)}= A ,
\mbox{ and } 
\pi_*(\tilde{p}_1(r))^{2r-1}= B_Rc_R +B_Lc_L 
$$
in 
$H^4(BSpin(4), \mathbb{Q})$,
where $A,B_R, B_L$ are universal numbers depending
on $\omega_0,\omega_L,\omega_R$. By pulling back,
$$\pi_*(p_1(r))^{2(r-1)}=A \mbox{ and
}\pi_*(p_1(r))^{2r-1}=B_R(2e+3\sigma) +B_L(2e-3\sigma).
$$
With this fact, applying (26) we can verify the statement 
for $\delta_{(r)}(\alpha)$. 
For general $\mathfrak{t}$, the statement  
follows from the product property. 
Thus we prove the conjecture.
q.e.d.

\end{document}